\documentclass[preprint,10pt]{elsarticle}

\usepackage{amsmath,amssymb,amsthm}
\usepackage{graphicx}
\usepackage{caption}
\usepackage{subcaption}
\usepackage[lined,boxed]{algorithm2e}
\usepackage{gensymb}
\usepackage{fancybox}
\usepackage{listings}
\usepackage{color}
\usepackage{enumerate}
\usepackage{marginnote}
\usepackage{breqn}

\usepackage[margin=1.in]{geometry}

\makeatletter \newcommand{\figcaption}{\def\@captype{figure}\caption} \newcommand{\tabcaption}{\def\@captype{table}\caption} \makeatother

\biboptions{numbers,sort&compress}
\usepackage{mathrsfs}

\usepackage{pinlabel}
\usepackage{caption}
\usepackage{seceqn}

\newcommand{\complexi}{\mathrm{i}}

\usepackage{hyperref}

\newtheorem*{theorem*}{Theorem}
\theoremstyle{definition}

\theoremstyle{remark}

\def\dsst{\displaystyle}
\def\L{\mathscr{L}}

\newcommand{\nv}{\mathbf{n}}
\newcommand{\biharm}{\nabla^4}
\newcommand{\G}{\mathcal{G}} 
\newcommand{\xv}{\mathbf{x}}
\newcommand{\yv}{\mathbf{y}}
\newcommand{\rv}{\mathbf{r}}
\newcommand{\Xv}{\mathbf{X}}
\newcommand{\iv}{\mathbf{i}}
\newcommand{\Phiv}{\mathbf{\Phi}}
\newcommand{\Wv}{\mathbf{W}}
\newcommand{\Fv}{\mathbf{F}}
\newcommand{\Uv}{\mathbf{U}}
\newcommand{\Vv}{\mathbf{V}}
\newcommand{\Fcal}{\mathcal{F}}
\newcommand{\diag}{\mathrm{diag}}
\newcommand{\fsolve}{\text{\bf fsolve}}
\newcommand{\tv}{\mathbf{t}}
\newcommand{\av}{\mathbf{a}}
\newcommand{\bv}{\mathbf{b}}

\usepackage{tikz}
\newlength{\tfwidth}
\newlength{\tfheight}
\newlength{\tfxa}
\newlength{\tfxb}
\newlength{\tfya}
\newlength{\tfyb}
\newcommand{\trimFigWithBox}[6]{%
\setlength\fboxsep{0pt}%
\setlength\fboxrule{1.0pt}
\fbox{\includegraphics[width=#2, clip, trim=#3 #4 #5 #6]{#1}}%
}
\newcommand{\trimFigNoBox}[6]{%
\setlength\fboxsep{1pt}
\setlength\fboxrule{0.0pt}
\fbox{\includegraphics[width=#2, clip, trim=#3 #4 #5 #6]{#1}}%
}
\newcommand{\trimFigHeightWithBox}[6]{%
\setlength\fboxsep{0pt}%
\setlength\fboxrule{1.0pt}
\fbox{\includegraphics[height=#2, clip, trim=#3 #4 #5 #6]{#1}}%
}
\newcommand{\trimFigHeightNoBox}[6]{%
\setlength\fboxsep{1pt}
\setlength\fboxrule{0.0pt}
\fbox{\includegraphics[height=#2, clip, trim=#3 #4 #5 #6]{#1}}%
}

\newsavebox\figBox

\newcommand{\trimw}[6]{%
\sbox\figBox{\includegraphics{#1}}
\setlength{\tfwidth}{\the\wd\figBox}
\setlength{\tfheight}{\the\ht\figBox}
\setlength{\tfxa}{\tfwidth*\real{#3}}%
\setlength{\tfxb}{\tfwidth*\real{#4}}%
\setlength{\tfya}{\tfheight*\real{#5}}%
\setlength{\tfyb}{\tfheight*\real{#6}}%
\trimFigNoBox{#1}{#2}{\tfxa}{\tfya}{\tfxb}{\tfyb}%
}

\newcommand{\trimwb}[6]{%

\sbox\figBox{\includegraphics{#1}}
\setlength{\tfwidth}{\the\wd\figBox}
\setlength{\tfheight}{\the\ht\figBox}
\setlength{\tfxa}{\tfwidth*\real{#3}}%
\setlength{\tfxb}{\tfwidth*\real{#4}}%
\setlength{\tfya}{\tfheight*\real{#5}}%
\setlength{\tfyb}{\tfheight*\real{#6}}%
\trimFigWithBox{#1}{#2}{\tfxa}{\tfya}{\tfxb}{\tfyb}%
}

\newcommand{\trimh}[6]{%
\sbox\figBox{\includegraphics{#1}}
\setlength{\tfwidth}{\the\wd\figBox}
\setlength{\tfheight}{\the\ht\figBox}
\setlength{\tfxa}{\tfwidth*\real{#3}}%
\setlength{\tfxb}{\tfwidth*\real{#4}}%
\setlength{\tfya}{\tfheight*\real{#5}}%
\setlength{\tfyb}{\tfheight*\real{#6}}%
\trimFigHeightNoBox{#1}{#2}{\tfxa}{\tfya}{\tfxb}{\tfyb}%
}

\newcommand{\trimhb}[6]{%

\sbox\figBox{\includegraphics{#1}}
\setlength{\tfwidth}{\the\wd\figBox}
\setlength{\tfheight}{\the\ht\figBox}
\setlength{\tfxa}{\tfwidth*\real{#3}}%
\setlength{\tfxb}{\tfwidth*\real{#4}}%
\setlength{\tfya}{\tfheight*\real{#5}}%
\setlength{\tfyb}{\tfheight*\real{#6}}%
\trimFigHeightWithBox{#1}{#2}{\tfxa}{\tfya}{\tfxb}{\tfyb}%
}
\begin{document}


\begin{frontmatter}


\title{
Numerical methods for thermally stressed shallow shell equations
}

%

\author[UCLA]{Hangjie~Ji}
 \ead{hangjie@math.ucla.edu}
  
\author[UL]{Longfei~Li}
 \ead{longfei.li@louisiana.edu}


\address[UCLA]{Department of Mathematics, University of California Los Angeles, Los Angeles, CA 90095, USA.}
\address[UL]{Department of Mathematics, University of Louisiana at Lafayette, Lafayette, LA 70504, USA.}

%

\begin{abstract}

We develop efficient and accurate numerical methods to solve a class of shallow shell problems of the von Karman type. The governing equations form a fourth-order coupled system of nonlinear biharnomic equations for the transverse deflection and Airy's stress function.
A second-order finite difference discretization with three iterative methods (Picard, Newton and Trust-Region Dogleg) are proposed for the  numerical solution of the nonlinear PDE system. 
Three simple boundary conditions and two application-motivated mixed boundary conditions are considered. 
Along with the nonlinearity of the system, boundary singularities that appear when mixed boundary conditions are specified are the main  numerical challenges.
Two approaches  that use either  a transition function or local corrections are developed to deal with these boundary singularities.  All the proposed numerical methods are validated using carefully designed numerical tests, where expected orders of accuracy and rates of convergence are observed. A rough run-time performance comparison is also conducted to illustrate the efficiency of our methods. 
As an application of the methods, a snap-through thermal buckling problem is considered. The critical thermal loads of shell buckling with various boundary conditions are numerically calculated, and snap-through bifurcation curves 
are also obtained using our numerical methods together with a pseudo-arclength continuation method.  Our  results are consistent with previous studies.

\end{abstract}
\end{frontmatter}


\section{Introduction}

High quality thin  glass sheets are ubiquitously used  in  modern  electronic devices such as smart phone screens and large TV displays. In order to  maintain the quality of glass sheets, the manufacturing processes (e.g., Corning's  revolutionary  ``Fusion''  process  \cite{fussionProcessCorningWebsite})  need to avoid any imperfections that can cause
deviations from a desired shape. Small non-idealities in manufacturing may produce non-uniform stress-free deflections that fails to meet  the  tightening specifications in glass industry.   For large thin glass sheets, variations can be introduced during the
glass forming process, and the  subsequent   cooling   and  transporting processes.  The cooling process can  yield heterogeneous ``frozen-in'' thermal stresses \cite{abbottmethods}, and  transporting  by partially  holding or clamping edges of a glass sheet \cite{fussionProcessCorningWebsite} can introduce various boundary stresses. These
thermal and boundary stresses  can generate further deflections to the products.
Therefore,  there is a pressing need for further investigation of thermal-elastic deformations in shallow shells so that improved manufacturing procedures can be designed to minimize defects during cooling and transporting. 

Over the years,   numerous theoretical works    have been developed  on related  areas of elasticity and solid mechanics; for example, theories on the thermal stability of regular-shaped structures like doubly curved, conical, spherical and cylindrical shells are developed in  \cite{mahayni1966thermal, thornton1993thermal}.  Many mathematical models have also been formulated to capture various aspects of  shallow shells,
among which models of von Karman type \cite{howell2009applied} provide a solid foundation for the characterization of  shallow shells. The nonlinear governing equations concerning  the transverse deflections of the shell and the Airy's stress function are able to characterize  large shell deformations that are of primary interest in industrial applications.
For a review of the geometrically nonlinear  theory of shallow shells exhibiting large displacement, we refer the readers to \cite{vlasov1964general,ventsel2001thin} and the reference therein.

Incorporating thermal stresses into shallow shell models is not trivial, and many literatures study the nonlinear loaded shell problems and the thermoelastic problems separately \cite{TimoshenkoWoinowsky59,donnell1976beams,hetnarski2009thermal}. The governing equations  for a flat thin isotropic plate under a thermal stress can be easily derived from the von Karman theory \cite{tauchert1986thermal}. 
Thermal buckling of plates and regular-shaped shells has also  been considered long time ago in \cite{mahayni1966thermal}. We note that the difference between a plate and a shell lies in the precast shape, which is flat for plates and curved for a shell in the stress-free stage. It was only until recently when Abbott et al. developed the thermoelastic theory for nonlinear thin shells of general shapes subject to thermal stresses \cite{abbottmethods};
the  model  is a system of two biharmonic PDEs nonlinearly   coupled together. 
Analytical solutions can hardly be obtained for this type of PDEs; therefore, numerical approaches are normally applied to investigate the solutions. 

To the best of our knowledge, there are no existing numerical studies on solving the nonlinear shallow shell  equations under thermal stresses developed in \cite{abbottmethods}. Nevertheless, 
a great number of numerical methods \cite{palsev1966expansion,dang2009iterative,hadjidimos1971numerical} have been proposed to solve the  biharmonic  equation, which is a fundamental part of the nonlinear shallow shell model. 
Common numerical methods can be applied 
successfully to solve the biharmonic equation with ordinary boundary conditions, such as  Dirichlet and Neumann boundary conditions. However, when mixed boundary conditions are involved, it is well-known that standard numerical methods perform poorly for elliptic PDEs around boundary singularities, which are  introduced by jump discontinuities in the mixed boundary conditions. In this case, both global methods (series-type method, Ritz method, etc.) and local methods (finite differences, finite elements, strip elements, etc.) suffer from  loss of accuracy. 
A related benchmark problem,  Motz's problem \cite{motz1946treatment}, that considers the Laplace's equation with Neumann-Dirichlet mixed boundary  conditions in a rectangular domain can be used to reveal the loss of accuracy due to boundary singularities; interesting readers are referred to \cite{li2000singularities} for an extensive survey on this topic.

To maintain the desired accuracy for the  biharmonic equation  with mixed boundary conditions, global methods usually require extremely high order approximations around the singularities; see for example   the series-based method  introduced   in \cite{narita1981application} to solve the biharmonic problem with mixed boundary conditions.   Meanwhile,  local methods need to be implemented with adaptive mesh refinement around the singularities or combined with singular function approximations, such as the numerical methods developed in \cite{eastep1982natural, mizusawa1987vibration, fan1984flexural} to study the vibration and buckling of plates.  
 Even though standard local methods combined with local mesh refinement can be applied to a large variety of singular problems with fewer requirements, singular function methods are preferred since it is  generally more efficient provided appropriate functions are chosen to fit the singularities. 
 Incorporating these function approximations requires the understanding of analytic forms of the boundary singularities. Among several special methods that take into consideration of local corrections of solutions, we in particular mention the methods developed by Richardson  \cite{richardson1970stick}  and  Poullikkas et al. \cite{poullikkas1998methods}.
Richardson \cite{richardson1970stick} applied the Wiener-Hopf method to obtain solutions to the biharmonic equation that involves clamped and simply supported mixed boundary conditions. Poullikkas et al. \cite{poullikkas1998methods} combined
the knowledge of fundamental solutions near the singularity with a least square routine to determine unknown coefficients in the numerical approximation.
 There have been numerous   other numerical approaches designed to  prevent the loss of accuracy for mixed boundary conditions, such as Galerkin method \cite{chia1985non}, Rayleigh-Ritz variational method \cite{leissa1980vibrations, liew1992vibration} and domain decomposition method \cite{liew1993use,liew1994use} to name just a few. 

For coupled nonlinear problems of thin plates with large deflections similar to the shallow shell system that we are interested in solving, several finite difference and finite element techniques \cite{bilbao2008family, leung1995symplectic, ribeiro1999geometrical}, boundary element methods \cite{wang2000dual} and Picard iterations \cite{uscilowska2011implementation} have been developed. In the study of nonlinear dynamics of shallow shells, different methods \cite{dowell1968modal, chia1988nonlinear, chuen1987non, abe2000non, kurpa2007nonlinear} have been applied to numerically solve the system of nonlinear PDEs with various boundary conditions.

In this paper,  we  focus on developing new efficient and accurate numerical methods to solve  the type of nonlinear biharmonic PDEs developed in \cite{abbottmethods} that incorporates the thermal stresses. The boundary conditions we consider for the coupled system 
are  both the standard simple boundary conditions derived from preserving the energy of the shell  and the mixed boundary conditions motivated by engineering applications.  We are particularly interested in the partially clamped mixed boundary conditions not only because it is more numerically challenging, but also because it is closely related to the  glass manufacturing applications; for example, the glass sheets are partially clamped  during the cooling and transporting process in  Corning's Fusion technique \cite{fussionProcessCorningWebsite}.   
 We develop and compare three numerical techniques to solve the nonlinear biharmonic system iteratively, and propose two approaches to address the boundary singularities of the partially clamped mixed boundary conditions. In addition,  strategies of regularizing the singular system with free boundary conditions are also proposed,  noting that the biharmonic system is singular with free boundary conditions since the displacement is only determined up to an arbitrary plane. All our numerical methods are carefully validated with numerical convergence studies, and the various methods are compared with a rough run-time performance comparison.
As an application of the proposed numerical methods, we solve  a snap-through thermal buckling problem to numerically obtain the critical thermal loads for several boundary conditions. In conjunction with a pseudo-arclength continuation method \cite{Keller87}, we are able to obtain snap-through bifurcation curves for those boundary conditions as well. 

The remainder of the paper is organized as follows. In section \ref{sec:formulation}, the model for elastic shallow shells subject to thermal stresses is formulated. Three types of simple boundary conditions and two types of mixed boundary conditions are introduced for the problem. In section \ref{sec:numericalScheme}, we propose  three iterative schemes (Picard, Newton and Trust-Region Dogleg) based on a common finite difference discretization of the coupled biharmonic system  to solve the  governing equations. In particular, a transition function approach and a local asymptotic solution approach are developed for special treatments for boundary singularities of the mixed boundary conditions. In section \ref{sec:numericalResults}, numerical results of the proposed approaches are presented, and the influences of thermal stresses, mixed boundary conditions and geometric nonlinearity to shallow shells are investigated via an example application problem.

\section{Formulation}
\label{sec:formulation}
We consider an elastic thin shallow shell defined on a rectangular domain $\Omega'=[x'_a,x'_b]\times[y'_a,y'_b]$ with a precast shape $w'_0(x',y')$ under the influence of  a temperature field $T'$; all the primed variables are dimensional quantities.
The governing equations of this problem consist of a coupled system of two biharmonic equations for the transverse deflection function $w'(x',y')$ and Airy stress function $\phi'(x',y')$  \cite{ventsel2001thin}: 
\begin{equation*}
\nabla^4\phi'  =  -\frac{1}{2}Eh\L{[}w',w'{]}-Eh\L{[}w'_0,w'{]}- \nabla^2 N(T'),
\end{equation*}
\begin{equation*}
D\nabla^4 w' =  \L{[}\phi',w'{]}+\L{[}\phi',w'_0{]}-\frac{1}{1-\nu}\nabla^2 M(T')+P',
\end{equation*}
where the bilinear operator $\L$ is defined as
$$
\L{[}u,v{]} \equiv u_{xx}v_{yy}+u_{yy}v_{xx}-2u_{xy}v_{xy}.
$$
Here $h$ is the thickness of the shell, $E$ is the Young's modulus, $\nu$ is the Poisson's ratio,  $D = {Eh^3}/{12(1-\nu^2)}$ is the bending stiffness and $P'$ accounts for any external load. In addition,
the resultant thermal force $N(T')$ and thermal moment $M(T')$ are given by 
\begin{equation*}
N(T') = {E \alpha}\int^{h/2}_{-h/2}T'~dz \quad\text{and} \quad M(T') = {E\alpha}\int^{h/2}_{-h/2}T'z~dz,
\end{equation*}
where $\alpha$ is the coefficient of thermal expansion, and $T'$ is the temperature distribution \cite{tauchert1986thermal}. We note that, for the special case with $w_0\equiv 0$,
the governing equations reduce to a classical nonlinear model for thin plates.

For a thin shallow shell,  the thickness $h$ is assumed to be small compared to other dimensions;  thus the temperature variations through the thickness can be ignored, namely $T'=T'(x',y')$. The thermal force and moment are then reduced to
\begin{equation*}
N(T')=Eh\alpha T'(x',y') \quad\text{and}\quad M(T') = 0,
\end{equation*}
and therefore the model can be simplified to
\begin{subequations}
\begin{equation}
\frac{1}{Eh}\nabla^4\phi'  =  -\frac{1}{2}\L{[}w',w'{]}-\L{[}w'_0,w'{]}- \alpha\nabla^2 T',\label{eq:dim_phi}
\end{equation}
\begin{equation}
D\nabla^4 w'  =  \L{[}\phi',w'{]}+\L{[}\phi',w'_0{]}+P'.\label{eq:dim_w}
\end{equation}\label{eq:dimEqns}
\end{subequations}
The biharmonic-type coupled system \eqref{eq:dimEqns} is in a form of the von Karman nonlinear static shallow shell equations \cite{ventsel2001thin}, and is valid for shallow shell with large transverse displacements. 

To non-dimensionalize the nonlinear shell model, we follow the similar scalings for the nonlinear von Karman plate equations used in \cite{dowell1974aeroelasticity, lyman2014application}; i.e.,
\begin{equation*}
x' = L {x}, \qquad y' = L{y}, \qquad w'_0 = \sqrt{\frac{D}{Eh}}{w_0}, \qquad w' = \sqrt{\frac{D}{Eh}}{w}, 
\end{equation*}
\begin{equation*}
\phi' = D{\phi},\qquad T' = \frac{D}{\alpha E h L^2}{T}, \qquad P' = \frac{D}{L^4}\sqrt{\frac{D}{Eh}}{P}.
\end{equation*}
Substituting the scales into the model \eqref{eq:dimEqns} leads to the dimensionless coupled system governing the displacement $w$ and the Airy stress function $\phi$,
\begin{subequations} \label{eq:coupledSystemNonlinear}
\begin{align}
\biharm \phi &=-\frac{1}{2}\L[w,w]  -\L[w_0,w]-f_\phi, \label{eq:coupledSystemNonlinear_phi}\\
\biharm w &=  \L[w,\phi] +\L[w_0,\phi]+f_w, \label{eq:coupledSystemNonlinear_w}
\end{align}
\end{subequations}
where the forcing terms are given by $f_\phi= \nabla^2T$ and $f_w=P$.

For a shallow shell with small deformations, the linear shallow shell theory is applicable and leads to a coupled system of linear partial differential equations,
\begin{subequations}\label{eq:coupledSystemLinear}
\begin{equation}\label{eq:coupledSystemLinear_phi}
\nabla^4\phi  = -\L{[}w_0,w{]}-f_\phi,
\end{equation}
\begin{equation}\label{eq:coupledSystemLinear_w}
\nabla^4 w  =  \L{[}w_0,\phi{]}+f_w.
\end{equation}
\end{subequations}
Numerical solutions to this linear system can serve as an initial guess for the iterative methods of solving the nonlinear shell equations \eqref{eq:coupledSystemNonlinear}.

For the shallow shell models \eqref{eq:coupledSystemLinear} and \eqref{eq:coupledSystemNonlinear}, we first consider three types of commonly used boundary conditions; these simple boundary conditions are normally derived from  conservation of energy  \cite{bilbao2008family}. Motivated by industrial applications \cite{fussionProcessCorningWebsite,abbottmethods}, we are interested in understanding the effects of mixed boundary conditions on the static behavior of the shallow shell. To this end, two partially clamped mixed boundary conditions are also considered.

\subsection{Simple boundary conditions}
To be specific, the three simple boundary conditions considered are
\begin{itemize}
\item Clamped boundary conditions: 
\begin{equation}
w = \frac{\partial w}{\partial \nv} = 0,\qquad \phi = \frac{\partial \phi}{\partial \nv} = 0,
\label{clampedbc}
\end{equation}
\item Simply supported boundary conditions: 
\begin{equation}
w = \frac{\partial^2 w}{\partial \nv^2}= 0, \qquad \phi = \frac{\partial^2 \phi}{\partial \nv^2} = 0,
\label{ssbc}
\end{equation}
\item Free boundary conditions:
\begin{equation}
\frac{\partial^2 w}{\partial{\nv}^2}+\nu\frac{\partial^2w}{\partial \tv^2}  =0,\qquad
\frac{\partial}{\partial \nv }\left[\frac{\partial^2 w}{\partial \nv^2} + (2-\nu)\frac{\partial^2 w}{\partial \tv^2}\right] = 0,\qquad \phi = \frac{\partial \phi}{\partial \nv} = 0,
\label{freebc}
\end{equation}
\end{itemize} 
where $\nv$ and $\tv$ are normal and tangential vectors to the boundary of the domain. The free boundary conditions must be complemented by a corner condition that impose zero forcing at corners of the rectangular region \cite{bilbao2008family}:
\begin{equation}\label{cornerforcing}
\frac{\partial^2 w}{\partial{x}\partial y} = 0.
\end{equation}

Noting that similar to the Possion equation, the forcing term $f$ of a biharmonic equation $\nabla^4 w=f$ with  free boundary conditions has to satisfy a compatibility condition, $\int_\Omega f dX=0.$
It is also important to point out that with the assumption of no temperature variation through the shell thickness, the resultant thermal moment $M$ does not affect the boundary conditions; while the influences from other boundary constraints have been investigated in \cite{arnold1989edge, qatu1992effects}.

\subsection{Mixed boundary conditions}

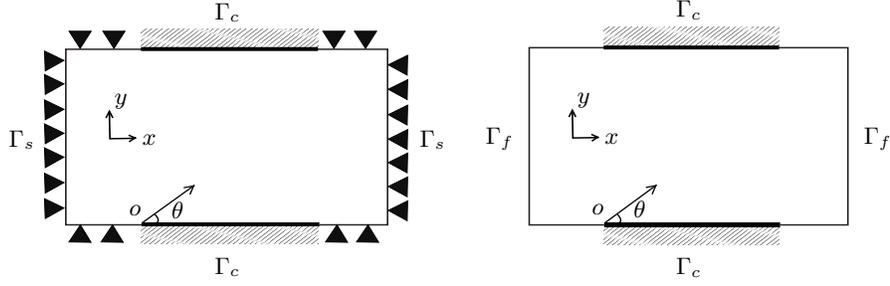
\begin{figure}[h]
\begin{center}
\labellist
\small
\pinlabel $\Gamma_{s}$ at 10, 400
\pinlabel $\Gamma_{s}$ at 590, 400
\pinlabel $\Gamma_{c}$ at 300, 600
\pinlabel $\Gamma_{c}$ at 300, 200
\pinlabel $x$ at 190, 400
\pinlabel $y$ at 150, 460
\pinlabel $\theta$ at 230, 290
\pinlabel $o$ at 170, 290
\endlabellist
\includegraphics[width=5cm,height=3cm]{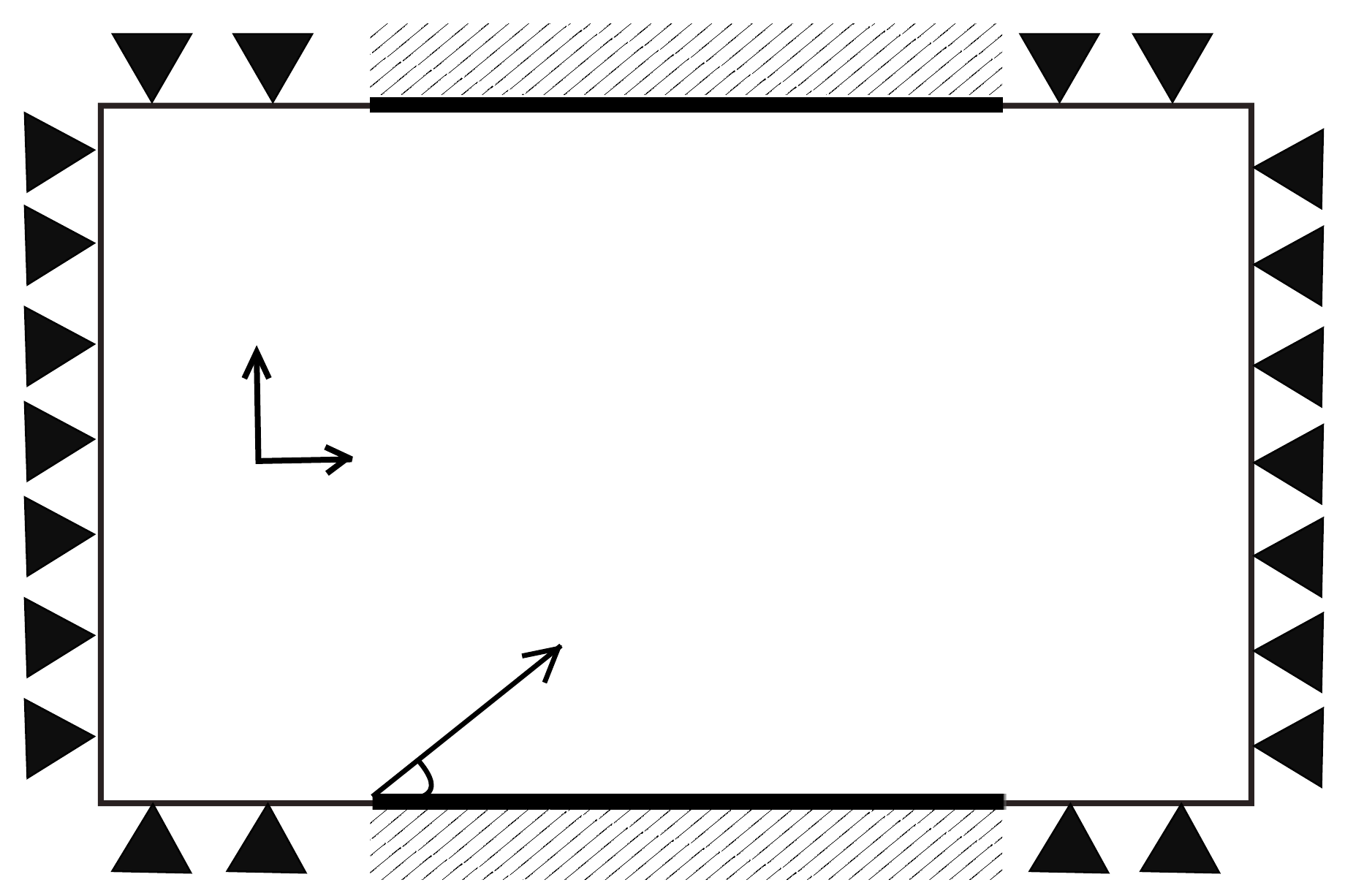}\hspace{0.4in}
\labellist
\small
\pinlabel $\Gamma_{f}$ at 30, 400
\pinlabel $\Gamma_{f}$ at 570, 400
\pinlabel $\Gamma_{c}$ at 300, 600
\pinlabel $\Gamma_{c}$ at 300, 200
\pinlabel $x$ at 190, 400
\pinlabel $y$ at 150, 460
\pinlabel $\theta$ at 230, 290
\pinlabel $o$ at 170, 290
\endlabellist
\includegraphics[width=5cm,height=3cm]{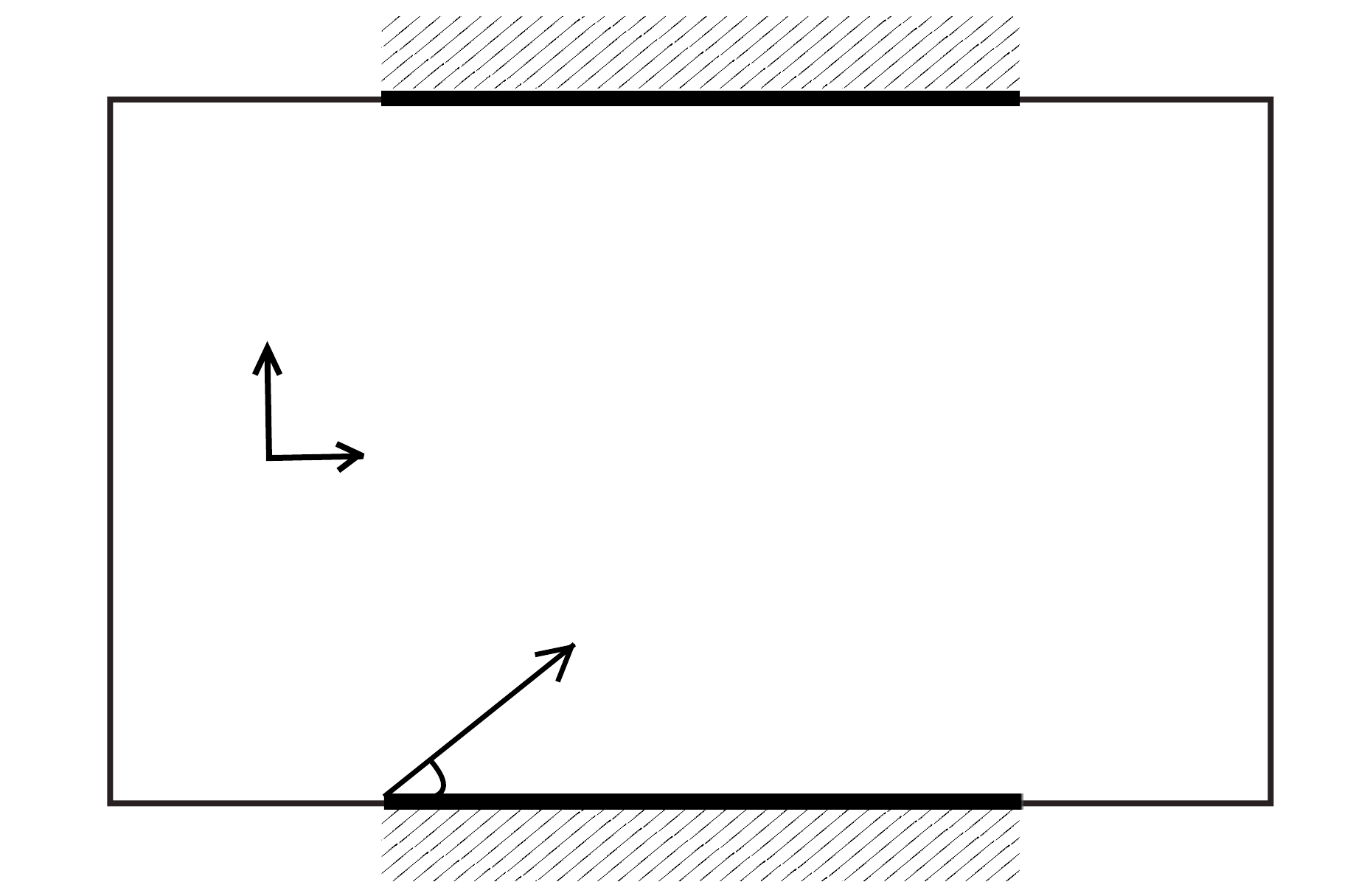}
\end{center}
\caption{Illustration of the mixed boundary conditions. Left:  clamped-supported boundary condition. Right:  clamped-free boundary condition.}
\label{domain}
\end{figure}

Rectangular shells with various mixed boundary conditions involving several combinations of the aforementioned simple boundary conditions have been studied analytically in \cite{TimoshenkoWoinowsky59}.
Shells with partially clamped edges have a great number of applications in industry. For example, in glass industry, large glass panels are sometimes partially clamped during cooling and transport processes. 

In this work, we focus on rectangular shells partially clamped on two opposite edges with the rest of the boundary being either simply supported or free. 
Specifically, we divide the boundary $\Gamma$ of the region $\Omega$ into two parts, $\Gamma = \Gamma_c \cup \overline{\Gamma}_{c}$, where $\Gamma_c$ represents the collection of the center portion of two opposite edges that are clamped. The rest of the boundary   $\overline{\Gamma}_{c}$ is  denoted by either $\Gamma_s$ if  simply supported  or $\Gamma_f$ if  free.  As illustrated in Figrue~\ref{domain}, the two partially clamped mixed boundary conditions considered are
\begin{itemize}
\item Clamped-Supported (CS): 
\begin{subequations}\label{eq:CSBC}
\begin{align}
w = \frac{\partial w}{\partial \nv} = 0,\qquad \phi = \frac{\partial \phi}{\partial \nv} = 0 &\quad\text{on}\quad \Gamma_c\\
w = \frac{\partial^2 w}{\partial \nv^2}= 0, \qquad \phi = \frac{\partial^2 \phi}{\partial \nv^2} = 0 &\quad\text{on}\quad \Gamma_s
\end{align}
\end{subequations}
\item Clamped-Free (CF): 
\begin{subequations}\label{eq:CFBC}
\begin{align}
w = \frac{\partial w}{\partial \nv} = 0,\qquad \phi = \frac{\partial \phi}{\partial \nv} = 0 &\quad\text{on}\quad \Gamma_c\\
\frac{\partial^2 w}{\partial{\nv}^2}+\nu\frac{\partial^2w}{\partial \tv^2}  =0,\qquad
\frac{\partial}{\partial \nv }\left[\frac{\partial^2 w}{\partial \nv^2} + (2-\nu)\frac{\partial^2 w}{\partial \tv^2}\right] = 0,\qquad \phi = \frac{\partial \phi}{\partial \nv} = 0  &\quad\text{on}\quad \Gamma_f
\end{align}
\end{subequations}

\end{itemize}
We note again that,  for the CF boundary condition,  the corner condition \eqref{cornerforcing} is needed at each of the free corners for completion. Mixed boundary conditions  of this type are sometimes referred to as strongly mixed boundary conditions  since the boundary conditions change at  inner points of the domain edges rather than  on the domain vertices \cite{dang2014simple}.  The sudden switch of  boundary conditions on an interior point of the boundary edge  introduces a jump discontinuity (singularity), as the two  boundary conditions cannot be both satisfied at the point of discontinuity. 

%

The effects of all the five boundary conditions on shallow shells will be demonstrated in  numerical examples. The difficulties of the numerical computation  lie in  the nonlinearity of the shallow shell equations, and in  
the singularities induced by discontinuities in the mixed boundary conditions. Numerical approaches proposed below  address both difficulties satisfactorily.

\section{Numerical schemes}\label{sec:numericalScheme}

In this section, three iterative methods (Picard, Newton and Trust-Region Dogleg) will be proposed for the numerical solution of the coupled  nonlinear system \eqref{eq:coupledSystemNonlinear} on a rectangular domain $\Omega=\left[x_a,x_b\right]\times\left[y_a,y_b\right]$.  All three iterative methods are based on a common spatial discretization of the coupled system utilizing a second-order accurate centered finite-difference scheme. 

\subsection{Spatial discretization}\label{sec:spatialDiscretization}
To be specific, the equations are solved on a Cartesian grid $\G_N$, with grid spacings $h_x=(x_b-x_a)/N$ and $h_y=(y_b-y_a)/N$, for a positive integer $N$:
\begin{equation}\label{eq:gridGN}
\G_N = \left\{\xv_\iv =\left(x_i,y_j\right)=(x_a+ih_x,y_a+jh_y): \quad  i, j = -2,-1,0,1, \dots,N+2\right\}.
\end{equation}
Here $\iv=(i,j)$ is a multi-index. We note that two layers of ghost points  are also included  at each boundary  to aid the discretization, and  we use $h=\min\{h_x,h_y\}$ to characterize the grid size.

Let $\Phi_\iv$ and  $W_{\iv}$ be the numerical approximation to $\phi(\xv_\iv)$ and $w(\xv_\iv)$, and denote $W_{0 \iv}=w_0(\xv_\iv)$, $F_{\phi\iv}=f_\phi(\xv_\iv) $ and  $F_{w\iv}=f_w(\xv_\iv)$ for notational brevity. For each grid  index $\iv$, the spatial discretized approximation to the coupled system reads
\begin{subequations}
\label{eq:discretizedCoupledSystemNonlinear}
\begin{equation}
\biharm_h \Phi_\iv =-\frac{1}{2}\L_h[W_\iv,W_\iv]  -\L_h[W_{0\iv},W_{\iv}]-F_{\phi\iv} \label{eq:discretizedCoupledSystemNonlinear_w},
\end{equation}
\begin{equation}
\biharm_h W_\iv =  \L_h[W_\iv,\Phi_\iv] +\L_h[W_{0\iv},\Phi_\iv]+F_{w\iv} \label{eq:discretizedCoupledSystemNonlinear_phi}.
\end{equation}
\end{subequations}
The discrete operators $\biharm_h$ and $\L_h$ are the standard centered finite-difference approximation to $\biharm$ and $\L$:
\begin{align*}
&\biharm_hU_{\iv} = \left(D_{xx} D_{xx}+2D_{xx}D_{yy}+D_{yy}D_{yy}\right)U_{\iv},\\
&\L_h[U_\iv,V_\iv]= D_{xx}U_{\iv}\,D_{yy}V_{\iv}+D_{yy}U_{\iv}\,D_{xx}V_{\iv}-2D_{xy}U_{\iv}\,D_{xy}V_{\iv},\\
\end{align*}
where 
\begin{align*}
 &D_{xx}U_{\iv} = \frac{U_{i+1,j}-2U_{i,j}+U_{i-1,j}}{h_x^2}, \quad D_{yy}U_{\iv} =  \frac{U_{i,j+1}-2U_{i,j}+U_{i,j-1}}{h_y^2},\\
 &D_{xy}U_{\iv} = \frac{U_{i+1,j+1}-U_{i-1,j+1}-U_{i+1,j-1}+U_{i-1,j-1}}{4h_xh_y}.
\end{align*}
The discretized system of equations \eqref{eq:discretizedCoupledSystemNonlinear} can be denoted as  two  matrix equations:
\begin{subequations}
\label{eq:matrixCoupledSystemNonlinear}
\begin{equation}
M_{\biharm_h} \Phiv =-\frac{1}{2}L_h[\Wv,\Wv]  -L_h[\Wv_0,\Wv]-\Fv_\phi \label{eq:matrixCoupledSystemNonlinear_phi},
\end{equation}
\begin{equation}
M_{\biharm_h} \Wv =  L_h[\Wv,\Phiv] +L_h[\Wv_0,\Phiv]+\Fv_w \label{eq:matrixCoupledSystemNonlinear_w}.
\end{equation}
\end{subequations}
Here $\Phiv$ denotes the column vector obtained by reshaping the grid function $\Phi_{\iv}$ and similar for all the other grid functions such as $\Wv$, $\Wv_0$, $\Fv_\phi$, $\Fv_w$, etc. Let  $M_{xx}$, $M_{yy}$ and $M_{xy}$ denote the matrices associated with the difference operators $D_{xx}$, $D_{yy}$ and $D_{xy}$, respectively;  the vector operators $M_{\biharm_h}$ and $L_h$ can then be written as 
\begin{align*}
&L_h[\Uv,\Vv] = M_{{xx}}\Uv\circ M_{{yy}}\Vv+M_{{yy}}\Uv\circ M_{{xx}}\Vv-2M_{{xy}}\Uv\circ M_{{xy}}\Vv,\\
&M_{\biharm_h} = M_{{xx}}M_{{xx}}+2M_{{xx}}M_{{yy}}+M_{{yy}}M_{{yy}},
\end{align*}
where $\Uv$ and $\Vv$ denote any column vectors that  are of the same size as $\Phiv$ (and $\Wv$). 
Here $A\circ B$ represents the Hadamard (entrywise) product of two matrices (or vectors) of the same dimensions, and $AB$ is the standard matrix multiplication.

To complete the statement of the discretized problem, appropriate discrete boundary conditions need to be applied to the matrix equation system \eqref{eq:matrixCoupledSystemNonlinear}. As is already noted in  section \ref{sec:formulation}, there are three simple boundary conditions and two partially clamped mixed boundary conditions that are considered for the shallow shell equations.

\subsection{Numerical implementation of simple boundary conditions}

The discretizations of the three simple boundary conditions are straightforward and their formulations are given by
\begin{itemize}
\item {Supported}
\begin{equation}
W_{\iv_b}=0,~ D^2_{\nv_{\iv_b}}W_{\iv_b}=0,~\Phi_{\iv_b}=0,~ D^2_{\nv_{\iv_b}}\Phi_{\iv_b}=0,\label{eq:discreteClampedBC}
\end{equation}
\item {Clamped}
\begin{equation}
W_{\iv_b}=0,~ D_{\nv_{\iv_b}}W_{\iv_b}=0,~\Phi_{\iv_b}=0,~ D_{\nv_{\iv_b}}\Phi_{\iv_b}=0,\label{eq:discreteSupportedBC}
\end{equation}
\item {Free} 
\begin{equation}
\left( D^2_{\nv_{\iv_b}}+\nu D^2_{\tv_{\iv_b}}\right)W_{\iv_b}=0,~D_{\nv_{\iv_b}}\left( D^2_{\nv_{\iv_b}}+(2-\nu) D^2_{\tv_{\iv_b}}\right)W_{\iv_b}=0,~\Phi_{\iv_b}=0,~ D_{\nv_{\iv_b}}\Phi_{\iv_b}=0.\label{eq:discreteFreeBC}
\end{equation}
\end{itemize}
Here $\iv_b=(i_b,j_b)$ denotes the index of a boundary node, $\nv_{\iv_b}$ and $\tv_{\iv_b}$ are the normal and tangential vectors at the boundary node, respectively. The directional difference operator is defined by 
$$
D_{\mathbf{a}}=a_1D_x+a_2D_y,
$$
where $\mathbf{a}=(a_1,a_2)$ is any given direction, and $D_x$ and $D_y$ are given by 
$$
D_{x}U_\iv = \frac{U_{i+1,j}-U_{i-1,j}}{2h_x} ~~\text{and}~~ D_{y}U_\iv = \frac{U_{i,j+1}-U_{i,j-1}}{2h_y}.
$$

\subsection{Numerical implementation of mixed boundary conditions}
The partially clamped mixed boundary conditions are less straightforward to implement, as the sudden switch of boundary conditions poses a boundary singularity at the point of jump  discontinuity. 
To address this singularity and achieve second order spatial accuracy, we explore the following two approaches to remove the boundary discontinuity.
\subsubsection{Transition function approach}
\begin{figure}[h!]
\centering
\includegraphics[width=6cm]{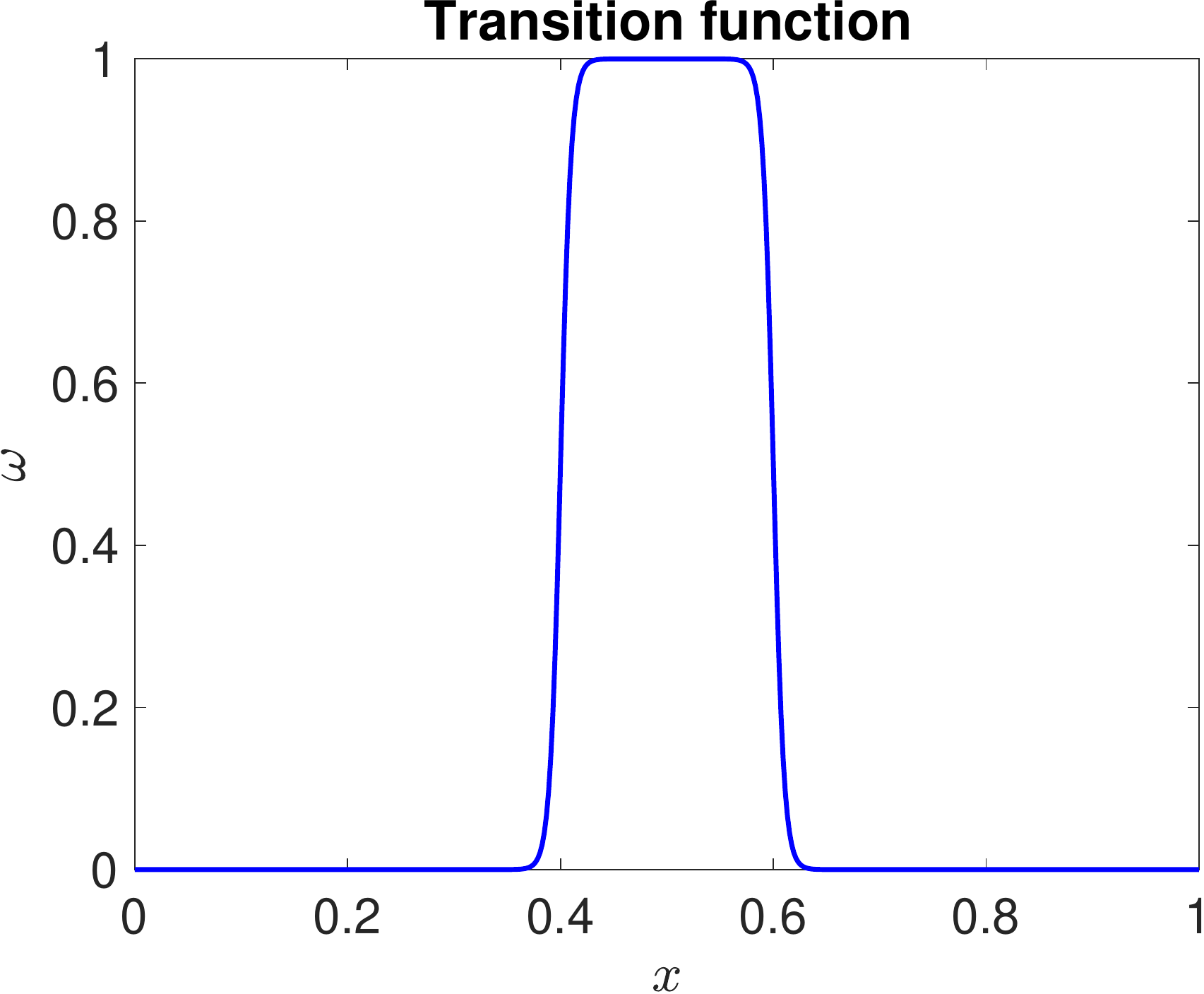}
\caption{The plot of the smooth transition function defined in \eqref{eqn:transitionFunction} for $0\leq x\leq 1$. It   provides a smooth transition from non-clamped region to the clamped region ($0.4<x<0.6$). The parameters are specified as $x_c=0.5$ and $r_c=0.1$.}\label{fig:smoothTrans}
\end{figure}

 The first approach considered  removes the discontinuity at the continuous level by introducing a transition function that  enables the boundary conditions to change  smoothly from  one to the other.  For simplicity, we restrict  the partially clamped region to the top and bottom boundaries of the square domain $\Omega=[x_a,x_b]\times[y_a,y_b]$. The clamped region is defined by  $$\Gamma_c=\{(x,y)~:y=y_a ~\text{or}~ y_b,  ~x_c-r_c<x<x_c+r_c\},$$
where $x_c$ and $r_c$  denote the $x$-coordinate of the center   and  the radius of the clamped region, respectively.
 The top-bottom  and left-right boundaries are given respectively by
\begin{align*}
\Gamma_{t,b}=\{(x,y)~:y=y_a ~\text{or}~ y_b,  ~x_a\leq x\leq x_b\}
~~\text{and}~~
\Gamma_{l,r}=\{(x,y)~:x=x_a ~\text{or}~ x_b,  ~y_a\leq y\leq y_b\}.
\end{align*}  
We then define a  transition function as following,
\begin{equation}
\omega(x) =  1-\frac{1}{2}\left[\tanh\left(\frac{|x-x_c|-r_c}{\epsilon}\right)+1\right],
\label{eqn:transitionFunction}
\end{equation}
where  $\epsilon$ is a parameter that controls the width of the transition region, and we set $\epsilon=0.01$ for the rest of this paper. We note that $\omega(x) = 1$ when $x$ is in the clamped region and $\omega(x) = 0$ otherwise (see Fig.~\ref{fig:smoothTrans}).

 With the introduction of $\omega(x)$, we redefine the two partially clamped mixed boundary conditions such that there is no discontinuity at the point where boundary conditions switch. The mixed clamped-supported boundary conditions are redefined to be
\begin{subequations}
\begin{align}
 w=0,~\frac{\partial^2 w}{\partial \nv^2}=0, \phi=0,~\frac{\partial^2 \phi}{\partial \nv^2}=0 ~\text{on}~\Gamma_{l,r},\\
 w=0,~\left(1-\omega(x)\right)\frac{\partial^2 w}{\partial \nv^2}+\omega(x)\frac{\partial w}{\partial \nv}=0 ~\text{on}~\Gamma_{t,b},\\
 \phi=0,~(1-\omega(x))\frac{\partial^2 \phi}{\partial \nv^2}+\omega(x)\frac{\partial \phi}{\partial \nv}=0 ~\text{on}~\Gamma_{t,b}.
\end{align}
\end{subequations}
Meanwhile the  mixed clamped-free boundary conditions are redefined to be
\begin{subequations}
\begin{align}
 \frac{\partial^2 w}{\partial \nv^2}+\nu \frac{\partial^2 w}{\partial \tv^2}=0,~ \frac{\partial }{\partial \nv}\left(\frac{\partial^2 w}{\partial \nv^2}+(\nu-2) \frac{\partial^2 w}{\partial \tv^2}\right)=0 ~\text{on}~\Gamma_{l,r},\\
 \left(1-\omega(x)\right)\left(\frac{\partial^2 w}{\partial \nv^2}+\nu \frac{\partial^2 w}{\partial \tv^2}\right)+\omega(x)w=0 ~\text{on}~\Gamma_{t,b},\\
\left(1-\omega(x)\right) \frac{\partial }{\partial \nv}\left(\frac{\partial^2 w}{\partial \nv^2}+(\nu-2) \frac{\partial^2 w}{\partial \tv^2}\right)+\omega(x)\frac{\partial w}{\partial \nv}=0 ~\text{on}~\Gamma_{t,b},\\
\phi=0,~\frac{\partial \phi}{\partial \nv}=0 ~\text{on}~\Gamma.
\end{align}
\end{subequations}
The discrete version of the redefined mixed boundary conditions are  readily obtained:
\begin{itemize}
\item {Clamped-Supported (CS)}
\begin{subequations}\label{eq:discreteCSBC}
\begin{align}
 W_{\iv_b}=0,~D_{\nv_{\iv_b}}^2W_{\iv_b}=0, ~\Phi_{\iv_b}=0,~D_{\nv_{\iv_b}}^2\Phi_{\iv_b}=0 ~\text{on}~\Gamma_{l,r},\\
 W_{\iv_b}=0,~\left(1-\omega(x_{i_b})\right)D_{\nv_{\iv_b}}^2W_{\iv_b}+\omega(x_{i_b})D_{\nv_{\iv_b}}W_{\iv_b}=0 ~\text{on}~\Gamma_{t,b},\\
 \Phi_{\iv_b}=0,~\left(1-\omega(x_{i_b})\right)D_{\nv_{\iv_b}}^2\Phi_{\iv_b}+\omega(x_{i_b})D_{\nv_{\iv_b}}\Phi_{\iv_b}=0  ~\text{on}~\Gamma_{t,b}.
\end{align}
\end{subequations}
\item {Clamped-Free (CF)}
\begin{subequations}\label{eq:discreteCFBC}
\begin{align}
\left( D^2_{\nv_{\iv_b}}+\nu D^2_{\tv_{\iv_b}}\right)W_{\iv_b}=0,~D_{\nv_{\iv_b}}\left( D^2_{\nv_{\iv_b}}+(2-\nu) D^2_{\tv_{\iv_b}}\right)W_{\iv_b}=0 ~\text{on}~\Gamma_{l,r},\\
\left(1-\omega(x_{i_b})\right)\left( D^2_{\nv_{\iv_b}}+\nu D^2_{\tv_{\iv_b}}\right)W_{\iv_b}+\omega(x_{i_b})W_{\iv_b}=0 ~\text{on}~\Gamma_{t,b},\\
\left(1-\omega(x_{i_b})\right)D_{\nv_{\iv_b}}\left( D^2_{\nv_{\iv_b}}+(2-\nu) D^2_{\tv_{\iv_b}}\right)W_{\iv_b}+\omega(x_{i_b})D_{\nv_{\iv_b}}W_{\iv_b}=0 ~\text{on}~\Gamma_{t,b},\\
\Phi_{\iv_b}=0,~D_{\nv_{\iv_b}}\Phi_{\iv_b}=0 ~\text{on}~\Gamma.
\end{align}
\end{subequations}

\end{itemize}

\subsubsection{Local asymptotic solution approach}
 The second approach considered removes the discontinuity at the discrete level by using local asymptotic analytical solutions. This approach is related to the Wiener-Hopf technique which has been applied to determine exact solutions to many problems including biharmonic equations with mixed boundary conditions \cite{abrahams2008matrix, richardson1970stick}.
Here we briefly describe this approach by demonstrating the numerical schemes for a simple biharmonic plate equation for the transverse deflection $w(x,y)$,
\begin{equation}\label{eq:WienerHopfBiharmEqn}
\nabla^4 w(x,y) = f_w(x,y), \quad (x,y)\in\Omega,
\end{equation}
which is subject to either CS or CF mixed boundary conditions. Since the equation for the Airy stress  $\phi$ is not considered in \eqref{eq:WienerHopfBiharmEqn}, only the $w$ parts of the boundary conditions \eqref{eq:CSBC} and \eqref{eq:CFBC} are needed.  

The key point of this approach is to locally construct an analytical solution in a small neighborhood of the singular point that approximately satisfies the boundary conditions. The analytical solution consists of a solution to the homogeneous version of \eqref{eq:WienerHopfBiharmEqn}, $\nabla^4 w = 0$, that satisfies the boundary conditions exactly and a solution of the inhomogeneous equation \eqref{eq:WienerHopfBiharmEqn} that asymptotically satisfies the boundary conditions. Then a leading order approximation of the analytical solution is used to design a special numerical scheme that bypasses the singular point with the assumption that the singular point lies on a grid point.

To be specific, we seek solutions of \eqref{eq:WienerHopfBiharmEqn} in a half-disk domain $\mathcal{B}(O,r_\epsilon)$ that is centered at the singular point $O$ with a small radius $r_\epsilon$. The equation \eqref{eq:WienerHopfBiharmEqn} is then converted to polar coordinates $(r,\theta)$ for convenience.
To simplify the discussion for boundary conditions, we assume that the boundary is clamped at $\theta=0$, and is either  supported or free  at $\theta=\pi$.
Note that  the normal derivatives  in  the  boundary conditions   become $\theta$ derivatives at  $\theta=0$ and $\pi$ since the domain $\mathcal{B}(O,r_\epsilon)$ is a half disk.
Using separation of variables and the fact that $w$ satisfies the clamped boundary conditions $w = \frac{\partial w}{\partial \theta}=0$ at $\theta=0$, we write the solution to the homogeneous biharmonic equation $\nabla^4 w = 0$ for an eigenvalue $\lambda$ in polar coordinates as,
\begin{equation}
w_\lambda(r,\theta) = r^{\lambda+1} f_{\lambda}(\theta),
\label{polarcoordinates}
\end{equation} 
where the  corresponding eigenfunction  $f_{\lambda}(\theta)$  is given by 
\begin{equation}
f_{\lambda}(\theta)=A\left[\cos\left((\lambda+1)\theta\right) - \cos\left((\lambda-1)\theta\right)\right]+B\left[ \frac{\sin\left((\lambda+1)\theta\right)}{\lambda+1}- \frac{\sin\left((\lambda-1)\theta\right)}{\lambda-1}\right].
\label{biharmonic_polarexpansion}
\end{equation} 
The unknown coefficients $A$ and $B$  will be determined by boundary conditions at  $\theta=\pi$. 


For the CS boundary conditions,  we have
\begin{equation}
\label{CSBCPolar}
 w=\frac{\partial^2 w}{\partial \theta^2}=0 \quad \text{on}\quad \theta=\pi.
\end{equation}
With \eqref{CSBCPolar} applied to $w_{\lambda}$ defined in \eqref{polarcoordinates}, 
the unknown coefficients in \eqref{biharmonic_polarexpansion} are determined,  and the eigenvalues and their eigenfunctions are found to be
\begin{equation}
f_{\lambda}(\theta)=
\begin{cases}
\cos(\lambda+1)\theta-\cos(\lambda-1)\theta \quad &\text{for} \quad \lambda = \frac{1}{2},\frac{3}{2}, \frac{5}{2} \cdots\\
(\lambda-1)\sin(\lambda+1)\theta-(\lambda+1)\sin(\lambda-1)\theta \quad &\text{for} \quad \lambda=2,3,4\cdots
\end{cases}.
\label{eigenfcn_clamped}
\end{equation}

Motivated by the modified method of fundamental solutions (MFS) proposed in \cite{poullikkas1998methods} where singular radial basis functions are integrated to approximate the biharmonic solution, we construct an approximate solution to the equation \eqref{eq:WienerHopfBiharmEqn}:
\begin{align}
\hat{w}_{\text{cs}}(r,\theta) = & \alpha_1 r^{\frac{3}{2}}f_{\frac{1}{2}}(\theta)+ \alpha_2 r^{\frac{5}{2}}f_{\frac{3}{2}}(\theta) +\alpha_3 r^{\frac{7}{2}}f_{\frac{5}{2}}(\theta)+\alpha_4 r^{3}f_{2}(\theta)+\alpha_5 r^{4}f_{3}(\theta) \nonumber\\ 
& + b_0 \epsilon^2 r^2 \ln r+ (a_0 \epsilon^2- b_0 \epsilon^2 \ln \epsilon) r^2.
\label{approx_cs_inhomo}
\end{align}
We note that the terms with coefficients $\alpha_i$'s form a  local approximation to the solution of the homogenous version of \eqref{eq:WienerHopfBiharmEqn} that satisfies the CS boundary conditions exactly, while the rest terms represent an approximation  to a solution of the inhomogeneous equation \eqref{eq:WienerHopfBiharmEqn} that satisfies the CS boundary conditions asymptotically as $r\to 0$.

\begin{figure}
\centering
\begin{subfigure}[t]{0.49\textwidth}
\centering
\labellist
\small
\pinlabel $\theta$ at 330, 300
\pinlabel $O$ at 260, 300
\pinlabel $1$ at 311, 370
\pinlabel $3$ at 350, 370
\pinlabel $2$ at 237, 370
\pinlabel $4$ at 308, 400
\pinlabel $13$ at 217, 250
\pinlabel $12$ at 154, 250
\pinlabel $14$ at 370, 250
\pinlabel $15$ at 430, 250
\pinlabel $5$ at 155, 370
\pinlabel $6$ at 427, 370
\pinlabel $7$ at 238, 450
\pinlabel $8$ at 345, 450
\pinlabel $9$ at 155, 450
\pinlabel $10$ at 427, 450
\pinlabel $11$ at 293, 530
\endlabellist
\includegraphics[width=7cm]{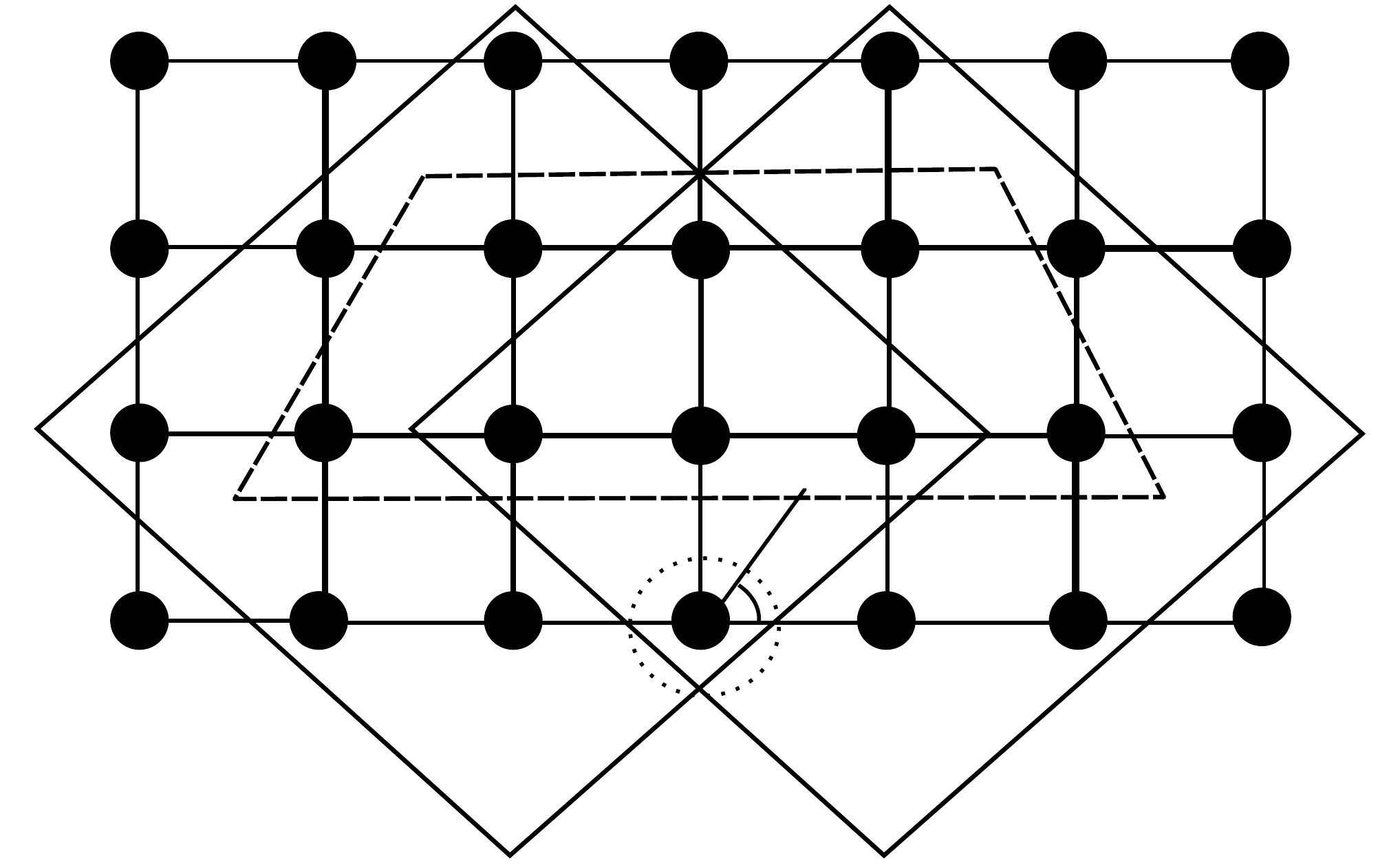}
\caption{CS boundary conditions}
\label{boundarygrid_cs}
\end{subfigure}
\begin{subfigure}[t]{0.49\textwidth}
\centering
\labellist
\small
\pinlabel $\theta$ at 330, 300
\pinlabel $O$ at 260, 300
\pinlabel $1$ at 311, 370
\pinlabel $3$ at 350, 370
\pinlabel $2$ at 237, 370
\pinlabel $4$ at 308, 400
\pinlabel $13$ at 217, 250
\pinlabel $12$ at 154, 250
\pinlabel $14$ at 370, 250
\pinlabel $15$ at 430, 250
\pinlabel $5$ at 155, 370
\pinlabel $6$ at 427, 370
\pinlabel $7$ at 238, 450
\pinlabel $8$ at 345, 450
\pinlabel $9$ at 155, 450
\pinlabel $10$ at 427, 450
\pinlabel $11$ at 293, 530
\endlabellist
\includegraphics[width=7cm]{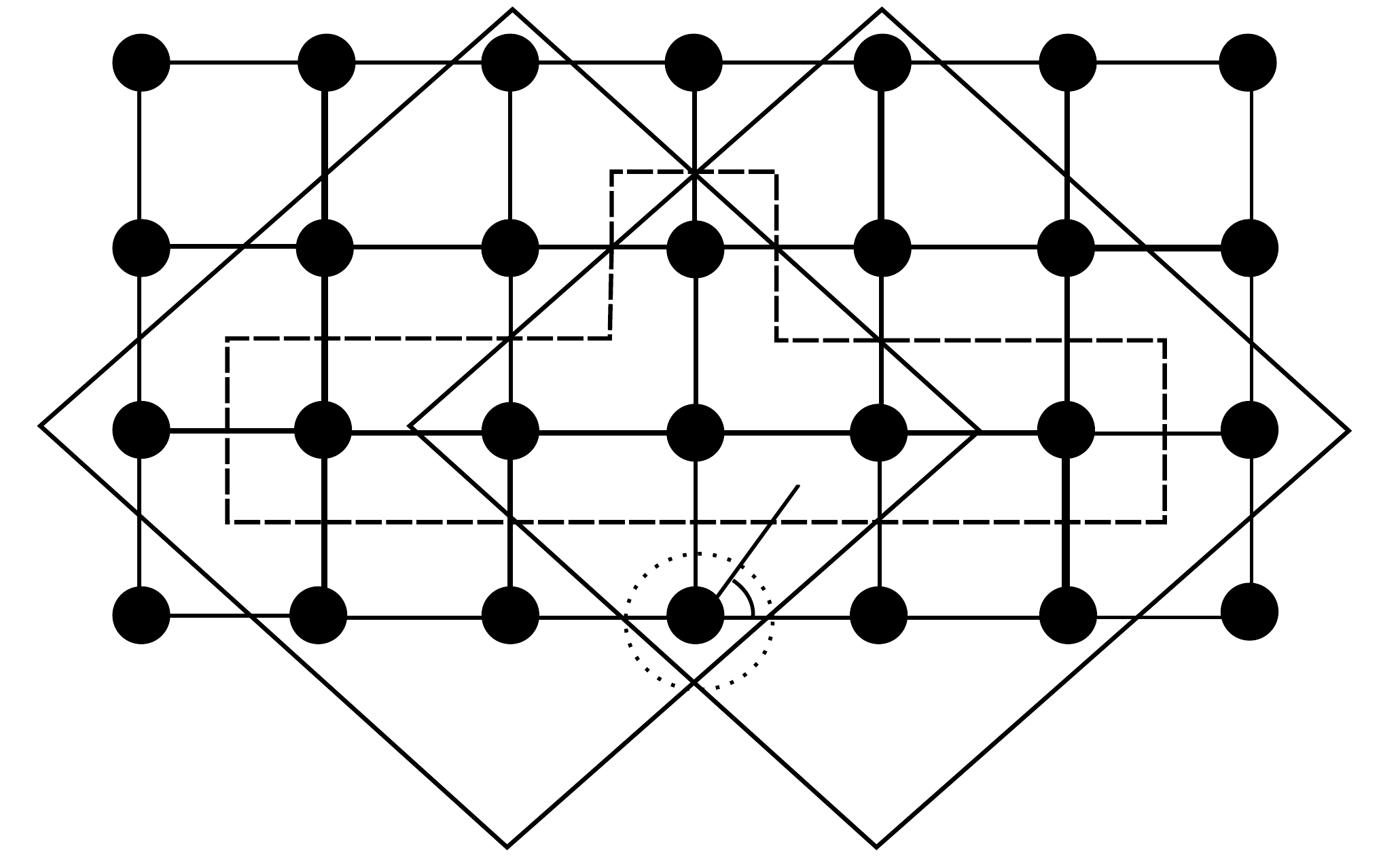}
\caption{CF boundary conditions}
\label{boundarygrid_cf}
\end{subfigure}
\caption{Schematic illustration of the stencil of the  numerical scheme near the singular point $O$.  The boundary condition \eqref{constraint_general}  derived from the local asymptotic solution approach for the corresponding boundary conditions is implemented at node 1 with its stencil  enclosed by the dashed line. Dirichlet boundary conditions are imposed at $O$. The standard centered-difference schemes or the standard  numerical boundary conditions are implemented at all the other nodes.}
\end{figure}
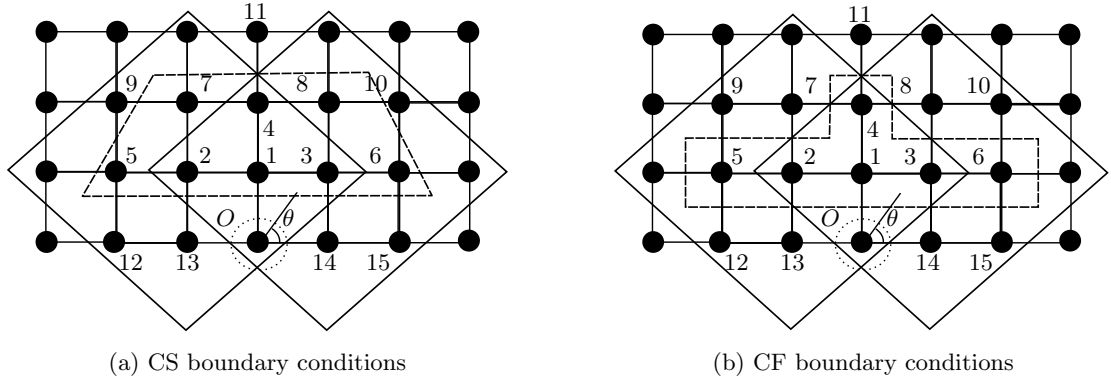
Similarly,  for the CF mixed boundary conditions, applying the following free boundary conditions
\begin{equation}
\frac{\partial^2 w}{\partial \theta^2}+\nu \frac{\partial^2 w}{\partial r^2}=\frac{\partial }{\partial\theta}\left[\frac{\partial^2 w}{\partial \theta^2}+(2-\nu)\frac{\partial^2 w}{\partial r^2}\right]=0 \quad \text{on}\quad \theta=\pi
\label{CFBCPolar}
\end{equation}
to $w_{\lambda}$ defined in \eqref{polarcoordinates}, we determine the unknown coefficients $A$ and $B$ in \eqref{biharmonic_polarexpansion}, and thus find all the possible eigenvalues,
\begin{equation}
\lambda=(n-\frac{1}{2})\pm \complexi K \mbox{ for } n\ge 1, \qquad \mbox{ where } \tanh(K\pi)=\frac{1+\nu}{2}.
\end{equation}

Then we derive a local approximation to the solution of \eqref{eq:WienerHopfBiharmEqn} associated with CF boundary conditions,
\begin{equation}
\hat{w}_{\text{cf}}(r,\theta) = \alpha_1 r^{\frac{3}{2}+\complexi K} f_{\frac{1}{2}+\complexi K}(\theta) + \alpha_2 r^{\frac{3}{2}-\complexi K} f_{\frac{1}{2}-\complexi K}(\theta) + b_0 r^2 \ln r + a_0 r^2 + a_0' r^2 \theta. 
\label{approx_cf_inhomo}
\end{equation}
In this expansion, the first two terms with coefficients $\alpha_i$'s are the first complex conjugate pairs of the solutions that exactly satisfy the homogeneous biharmonic equation with the mixed CF boundary conditions.  The rest terms form a local approximation to the inhomogeneous equation that satisfies the CF boundary conditions approximately.  In particular, the biharmonic terms $r^2$ and $r^2 \theta$ asymptotically satisfy the boundary conditions  at $\theta=0$ as $r\to 0$ with a truncation error of $O(r^2)$; those terms would satisfy the free boundary conditions exactly at $\theta=\pi$ if  proper coefficients are chosen. The term $r^2 \ln r$ is included in the expansion as it is important for matching the local approximation in the singular subdomain to the outer solution at $r = O(1)$. 

Based on the idea that overlapping ``patches'' with different stencils can be  used to approximate local solutions from the Flexible Local Approximation MEthod (FLAME) \cite{tsukerman2006class}, the analytical approximate solutions $\hat{w}_{\text{cs}}(r,\theta)$ given by \eqref{approx_cs_inhomo} and $\hat{w}_{\text{cf}}(r,\theta)$ given by  \eqref{approx_cf_inhomo}   are incorporated to remove the boundary singularity from the numerical scheme for the CS and CF mixed boundary conditions, respectively. 
For illustration purpose, we label the grid points near the boundary singularity $O$
 as shown in Figure~\ref{boundarygrid_cs} and Figure~\ref{boundarygrid_cf} for each case.  Applying  \eqref{approx_cs_inhomo} for the CS boundary conditions or \eqref{approx_cf_inhomo} for the CF boundary conditions to the grid points adjacent to the singular point $O$, we obtain a set of algebraic equations
\begin{equation}\label{eq:algebraicEqm}
w_i = \hat{w}_{\text{bc}}(r_i,\theta_i), \quad \text{bc}=\text{cs or  cf},
\end{equation}
where    $w_i$ represents the numerical approximation of the solution $w(x,y)$ at node $i$  and $(r_i,\theta_i)$ denote its polar coordinates (see Fig.~\ref{boundarygrid_cs} and Fig.~\ref{boundarygrid_cf}). Here  $i=1,\dots,8$ for the CS case and $i=1,\dots,6$ for the CF case. 
Eliminating the unknown coefficients in \eqref{approx_cs_inhomo} and \eqref{approx_cf_inhomo} using the  algebraic equations \eqref{eq:algebraicEqm} leads to a relationship 
\begin{equation}
\dsst \sum_{j=1}^p C_jw_j=0.
\label{constraint_general}
\end{equation}
where $p=8$ for the CS boundary conditions, and $p=6$ for CF boundary conditions.
This numerical  boundary condition is then implemented at  node $1$ in both cases (see Fig.~\ref{boundarygrid_cs} and Fig.~\ref{boundarygrid_cf}).
Note that Dirichlet boundary conditions are imposed at node $O$ and standard centered-difference schemes with appropriate numerical boundary conditions as discussed in section \ref{sec:spatialDiscretization} are used for all the other grid points. For example, the diamond shape stencils in Fig.~\ref{boundarygrid_cs} and Fig.~\ref{boundarygrid_cf} are used to indicate all the nodes involved in the standard finite difference approximation for the biharmonic operator.

It is important to point out that, without using the local asymptotic solution approach to remove the boundary singularity, traditional methods simply enforce one of the boundary conditions involved in  the mixed boundary conditions at the singular point and proceed with  the standard numerical  schemes and boundary conditions at all the grid points. In the traditional methods, the singularity, which is left untreated in the discretized system,  deteriorates the order of accuracy of the whole system. For comparison purpose, a traditional approach that  imposes the clamped boundary conditions at the singular point for both the CS and CF cases will be implemented and its numerical results will be compared with that of the asymptotic analytical solution approach.

%

\subsection{Initial guess}\label{sec:initialGuess}

The system \eqref{eq:matrixCoupledSystemNonlinear} together with  the discrete boundary conditions are solved iteratively using one of the following algorithms. We note that  all the iterative methods  start with a given initial guess denoted by $(\Phiv^0,\Wv^0)$. Unless otherwise noted, we use the precast shell shape $\Wv_0$ as the initial guess for $\Wv^0$, and use the solution to the following Airy stress $\Phiv$ equation 
$$
M_{\biharm_h} \Phiv^0 =-\frac{1}{2}L_h[\Wv^0,\Wv^0]  -L_h[\Wv_0,\Wv^0]-\Fv_\phi
$$
as the initial guess for $\Phiv$.

\subsection{Picard method}
 Motivated by \cite{uscilowska2011implementation}, we propose a Picard-type iterative method to solve the matrix equations described in Alogirthm \ref{alg:picardSolve}. It is important to note that the Picard method decouples the shallow shell equations by solving two biharmonic equations \eqref{eq:picarSolve_phi} \& \eqref{eq:picarSolve_w}  independently at each iteration step. Each of the biharmonic equations  has a matrix dimension that is four times smaller than the original coupled system; therefore it could potentially be more efficient in  overall performance than the other two algorithms that solve the coupled system as a single matrix equation. The efficiency of Picard method is confirmed by the numerical tests presented in  \S\ref{sec:numericalResults}. We also have an option ($\delta\in[0,1]$) to treat the $\Wv$ equation semi-implicitly; the scheme is  explicit for  $\delta=0$ and  implicit   for $\delta=1$.

\begin{algorithm}[h]
 \KwData{given initial guess: $(\Phiv^0,\Wv^0)$}
 \KwResult{numerical solutions to the shallow shell equations \eqref{eq:matrixCoupledSystemNonlinear_phi} \& \eqref{eq:matrixCoupledSystemNonlinear_w}: $(\Phiv,\Wv)$}
initialization: set $\Phiv^k=\Phiv^0,\,\Wv^k=\Wv^0$, $converged=\text{false}$ and $step=0$\;   
 \While{not converged {\bf and} step $<$ maxIter}{
  solve $\Phiv$ equation: 
  \begin{equation}\label{eq:picarSolve_phi}
  M_{\biharm_h} \Phiv^{k+1} =-\frac{1}{2}L_h[\Wv^k,\Wv^k]  -L_h[\Wv_0,\Wv^k]-\Fv_\phi;
  \end{equation}\textcolor{white}{\;}
  solve $\Wv$ equation: 
  \begin{equation}\label{eq:picarSolve_w}
  M_{\biharm_h} \Wv^{k+1}=\delta L_h[\Wv^{k+1},\Phiv^{k+1}]+(1-\delta) L_h[\Wv^k,\Phiv^{k+1}] +L_h[\Wv_0,\Phiv^{k+1}]+\Fv_w;
  \end{equation}\textcolor{white}{\;}
\If{$||\Phiv||_{\infty}+||\Wv||_{\infty} < \text{tol}$}
   {$converged=$true\;}
   prepare for next iteration step: $\Phiv^k=\Phiv^{k+1}$, $\Wv^k=\Wv^{k+1}$, $step$++\;
}
\eIf{converged}{
 solutions obtained: $\Phiv=\Phiv^{k}$, $\Wv=\Wv^{k}$\;}
 {iteration failed after max number of iteration steps reached\;}
 \caption{A Picard-type iterative method for the coupled system \eqref{eq:matrixCoupledSystemNonlinear}, where
 $tol$ is the tolerance and $maxIter$ is the maximum number of iterations allowed. The nonlinear term in  the $\Wv$ equation \eqref{eq:picarSolve_w} is treated semi-implicitly with $\delta\in[0,1]$ representing the degree of implicity.} \label{alg:picardSolve}
\end{algorithm}

\subsection{Newton's method}
The most obvious approach of solving a  nonlinear system  of equations is Newton's method. To this end, we also develop a Newton solver to numerically solve the shallow shell equations.  We rewrite the equations  \eqref{eq:matrixCoupledSystemNonlinear} as 
\begin{equation}
\Fcal(\Xv) = \mathbf{0},
\end{equation}
where 
$$
\Xv=\begin{bmatrix}
\Phiv\\
\Wv
\end{bmatrix}
\quad \text{and} \quad
\Fcal(\Xv) =
\begin{bmatrix}
M_{\biharm_h} \Phiv +\frac{1}{2}L_h[\Wv,\Wv]  +L_h[\Wv_0,\Wv]+\Fv_\phi \\
M_{\biharm_h} \Wv -  L_h[\Wv,\Phiv] -L_h[\Wv_0,\Phiv]-\Fv_w 
\end{bmatrix}.
$$
The key for the Newton's method, as well as the Trust-Region Dogleg Method discussed below in \S\ref{sec:fsolve}, to work efficiently for problems with a large number of grid points is to  find a way to efficiently evaluate the Jacobian matrix of $\Fcal(\Xv)$.   Fortunately, in our case, we are able to determine  the analytical expression of the Jacobian matrix.

With the introduction of a matrix function
$$
M_{L_h}(\Uv) = \diag(M_{xx}\Uv) M_{yy}+\diag(M_{yy}\Uv) M_{xx}-2\diag(M_{xy}\Uv) M_{xy},
$$
the bilinear operator can be written in terms of a matrix product
$
L_h[\Uv,\Vv]  = M_{L_h}(\Uv)  \Vv.
$
Here $\diag(\Vv)$ represents  the diagonal matrix with the elements of vector $\Vv$ on the main diagonal.
The Jacobian matrix of $\Fcal(\Xv)$ is therefore  readily obtained:
\begin{equation}\label{eq:JacobeamMatrix}
J(\Xv) = \frac{\partial \Fcal(\Xv)}{\partial \Xv} = 
\begin{bmatrix}
M_{\biharm_h} & M_{L_h}(\Wv)+ M_{L_h}(\Wv_0) \\
-M_{L_h}(\Wv)- M_{L_h}(\Wv_0) & M_{\biharm_h} -M_{L_h}(\Phiv)
\end{bmatrix}.
\end{equation}
Our Newton solver is summarized in Algorithm~\ref{alg:newtonSolve}.
\begin{algorithm}[h]
 \KwData{given initial guess: $\Xv^0=\left[{\Phiv^0}^T,{\Wv^0}^T\right]^T$}
 \KwResult{numerical solutions to the shallow shell equations \eqref{eq:matrixCoupledSystemNonlinear_phi} \& \eqref{eq:matrixCoupledSystemNonlinear_w}: $\Xv=\left[{\Phiv}^T,{\Wv}^T\right]^T$}
initialization: set $\Xv^k=\Xv^0$, $converged=\text{false}$ and $step=0$\;   
 \While{not converged {\bf and} step $<$ maxIter}{
  \begin{align}
 	&\Delta \Xv =-J(\Xv^k)^{-1}\Fcal(\Xv^{k});\label{eq:newtonStep}\\
	&\Xv^{k+1} =\Xv^{k}+\Delta \Xv;
  \end{align}\textcolor{white}{\;}
\If{$||\Phiv||_{\infty}+||\Wv||_{\infty} < \text{tol}$}
   {$converged=$true\;}
   prepare for next iteration step: $\Xv^k=\Xv^{k+1}$, $step$++\;
}
\eIf{converged}{
 solutions obtained: $\Xv=\Xv^{k}$\;}
 {iteration failed after max number of iteration steps reached\;}
 \caption{Newton's method for the coupled system \eqref{eq:matrixCoupledSystemNonlinear} where $tol$ is the tolerance and $maxIter$ is the maximum number of iterations allowed.} \label{alg:newtonSolve}
\end{algorithm}

\subsection{Trust-Region Dogleg method}\label{sec:fsolve}
For comparison purpose, we also solve the discretized system \eqref{eq:discretizedCoupledSystemNonlinear_w}--\eqref{eq:discretizedCoupledSystemNonlinear_phi} using the built-in function $\fsolve$ of  MATLAB (The MathWorks, Inc., Natick, MA).
The underlining algorithm that we choose when using $\fsolve$ is  Trust-Region Dogleg Method, which is a variant of the Powell dogleg method \cite{powell70}. The key difference between this method and the Newton's method lies in the  procedure for computing the step $\Delta \Xv$. In contrast to the Newton's method that computes the step as in equation \eqref{eq:newtonStep}, the Trust-Region Dogleg Method constructs steps from a convex combination of a Cauchy step (a step along the steepest descent direction) and a Gauss-Newton step. The trust-region technique improves the robustness and is able to handle the case when the Jacobian matrix is singular. Further details about this method can be found in the documentation of MATLAB's Optimization Toolbox \cite{optToolbox17a}.

\subsection{Displacement regularization for free boundary conditions}\label{sec:removeSingularityFreeBC}
We note that the displacement equation (a biharmonic equation) subject to free boundary conditions is singular since the displacement is only determined up to an arbitrary plane $c_1x+c_2y+c_3$ ($c_i$'s are arbitrary constants). In addition, similar to a Poisson equation with Neumann boundary condition,  the biharmonic equation is  solvable only if the right hand side satisfies a compatibility condition. In order to solve this singular system one need to eliminate three equations and replace them with equations that set the values of $w$ at three points. Instead of picking the equations to be replaced, we prefer to use a different approach which is better conditioned. This approach is motivated by the method used by Henshaw and Petersson  \cite{splitStep2003} to regularize the  pressure Poisson equation with Neumann boundary condition; it is a crucial step to solve  this singular pressure Poisson equation  for the split-step scheme proposed by the authors to solve   the incompressible Navier-Stokes equations with no-slip wall boundary conditions. Let the biharmonic equation for the displacement with free boundary conditions   be denoted as a matrix equation
\begin{equation}\label{eq:WeqFree}
A\Wv=\bv,
\end{equation}
where the matrix $A$ is singular with the dimension of its null space being $3$. Since the solution is determined up to an arbitrary plane, the right null of $A$ is found to be  $Q=[\xv, \yv, \rv]$, where $\xv$ and $\yv$ are column vectors obtained by reshaping the $x$ and $y$ coordinates of all the grid points and $\rv$ is the vector with all components equal to one. Instead of solving the singular equation \eqref{eq:WeqFree}, we seek solutions of the augmented system
 
\begin{equation}\label{eq:WeqFreeAug}
\begin{bmatrix}
A & Q\\
Q^T & \mathbf{0}_{3\times3}
\end{bmatrix}
\begin{bmatrix}
\Wv\\
\av
\end{bmatrix}
=
\begin{bmatrix}
\bv\\
 \mathbf{0}_{3\times1}
\end{bmatrix}.
\end{equation}
It is well-known that the saddle point problem \eqref{eq:WeqFreeAug} is non-singular and has an unique solution \cite{BenziEtal05}. The last three equations ($Q^T\Wv= \mathbf{0}_{3\times1}$) will set the mean values of $xw$, $yw$ and $w$ to be zero.

\section{Numerical results}\label{sec:numericalResults}
We now present the results for a sequence of simulations to demonstrate the properties of our numerical approaches for the shallow shell equations. We begin with mesh refinement studies to illustrate some basic properties of our approach. Two simple tests for solving a single biharmonic equation are considered first because the accurate solution of a  biharmonic equation  is an essential component for the coupled system. The next set of tests are designed for the coupled system.  We perform mesh refinement studies first for the simplified linear system and then for the nonlinear system using all three of the aforementioned iterative schemes. Efficiency of the iterative schemes are also compared. 
In order to study the effects of boundary conditions,  a numerical example of nonlinear shell  with a precast shell shape and localized thermal forcing is considered with all the proposed simple and mixed boundary conditions \eqref{clampedbc}--\eqref{eq:CFBC}. Finally, as an application of the numerical methods, we study the snap-through thermal buckling problem with an unstressed shell shape. A pseudo-arclength continuation (PAC) method \cite{Keller87} is utilized to find the snap-through bifurcation; one of our iterative methods is used to solve the resulted  system  at the each step of the continuation method.  

For simplicity, unless otherwise noted all the test problems considered are   on a unit square domain, i.e., $\Omega=[x_a,x_b]\times[y_a,y_b]$ with $x_a=y_a=0$ and $x_b=y_b=1$, and the partially clamped region  on the boundary for the two mixed boundary conditions are both assumed to be $\Gamma_c=\{(x,y) :~y=0 ~\text{or}~ 1,  ~0.4<x<0.6\}$.

\subsection{Mesh refinement study}
We use the method of manufactured solutions  to construct exact solutions of test problems by adding forcing
functions to the governing equations. The forcing is specified so that a chosen function becomes an exact solution to the forced equations. Here the approach is used to verify the order of accuracy of the numerical solutions of:  (I) a single biharmonic equation; (II) the linear coupled system and (III) the nonlinear coupled system. Numerical solutions subject to all of the five boundary conditions \eqref{clampedbc}--\eqref{eq:CFBC} are obtained separately.

\subsubsection{Biharmonic equation}
We note that solving biharmonic  equation plays a vital role for all the iterative algorithms proposed for the numerical solution of the  shallow shell equations.  The Picard method essentially solves two biharmonic equations at each iteration and for the Newton and  Trust-Region Dogleg methods the discretized biharmonic operator forms the diagonal blocks of the Jacobian matrix \eqref{eq:JacobeamMatrix} that are inverted at each step. Given its importance, 
the accuracy of the numerical solution of a single biharmonic equation 
\begin{equation}\label{eq:biharmEqn}
\biharm w= f_w
\end{equation}
is verified first. The exact solution $w_e(x,y)$ of the biharmonic equation is chosen to be either of the following:
\begin{enumerate}
\item trigonometric test
$$w_e(x,y)  = \sin^4\left(2\pi\frac{x-x_a}{x_b-x_a}\right)\sin^4\left(2\pi\frac{y-y_a}{y_b-y_a}\right),$$
\item polynomial  test
$$w_e(x,y)  = \frac{1}{100}\left[\frac{(x-x_c)(x-x_a)(x-x_b)(y-y_c)(y-y_a)(y-y_b)}{l_x^3l_y^3}\right]^7,$$
\end{enumerate}
where $x_c=(x_a+x_b)/2$, $y_c=(y_a+y_b)/2$, $l_x=(x_b-x_a)/3$ and $l_y=(y_b-y_a)/3$.
Both exact solutions satisfy all the boundary conditions. The forcing term is then given by
$$
f_w(x,y) = \biharm w_e,
$$
which is plotted in Fig.~\ref{fig:biharmTestForcingContour} for both the trigonometric and polynomial tests.  
{
\newcommand{\figWidth}{8cm}
\newcommand{\trimfig}[2]{\trimw{#1}{#2}{0.}{0.}{0.}{0.0}}
\begin{figure}[h!]
\begin{center}
\begin{tikzpicture}[scale=1]
\useasboundingbox (0.0,0.0) rectangle (17.,6);  
\draw(-0.5,0.0) node[anchor=south west,xshift=0pt,yshift=0pt] {\trimfig{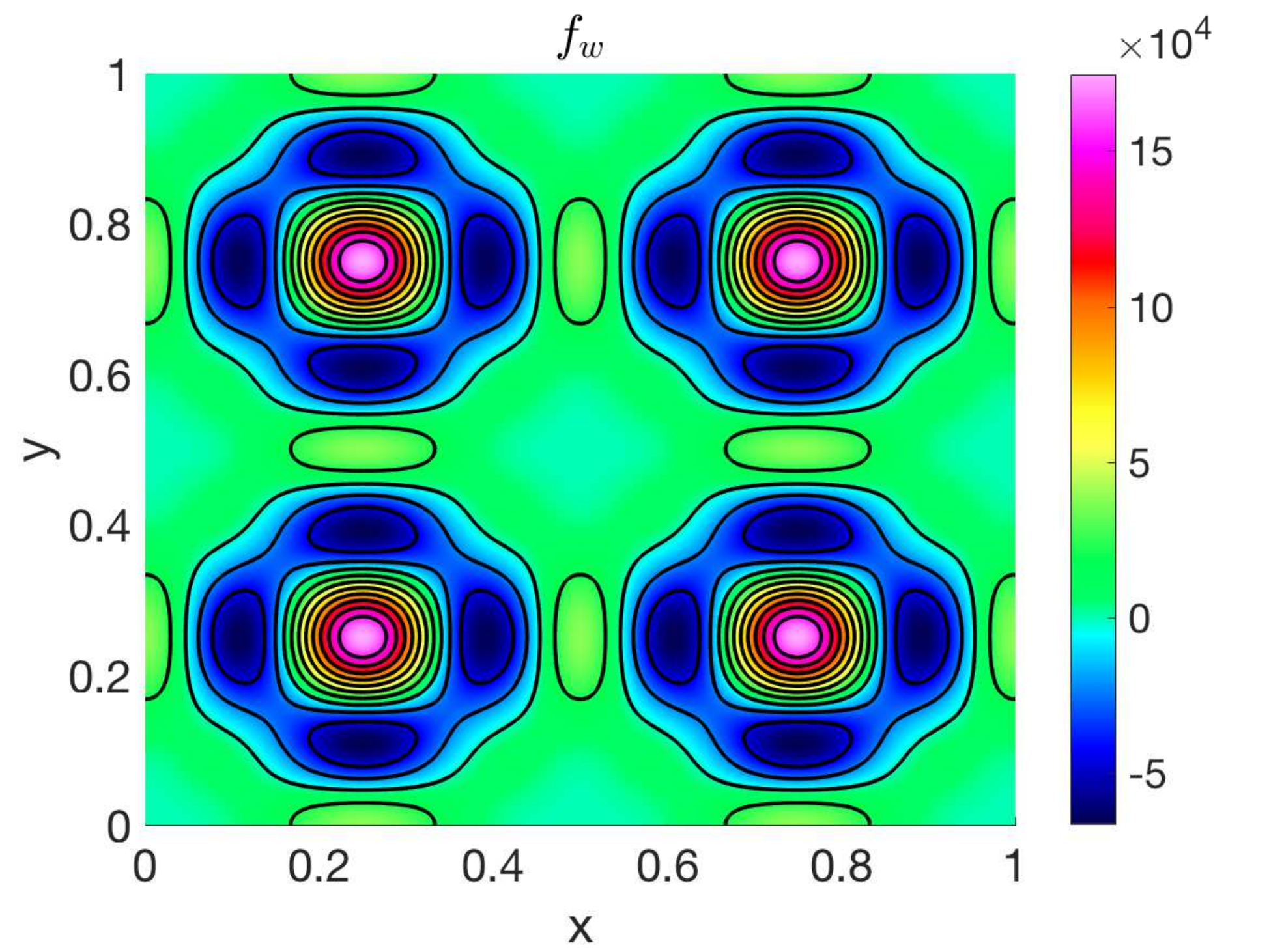}{\figWidth}};
\draw(8.0,0.0) node[anchor=south west,xshift=0pt,yshift=0pt] {\trimfig{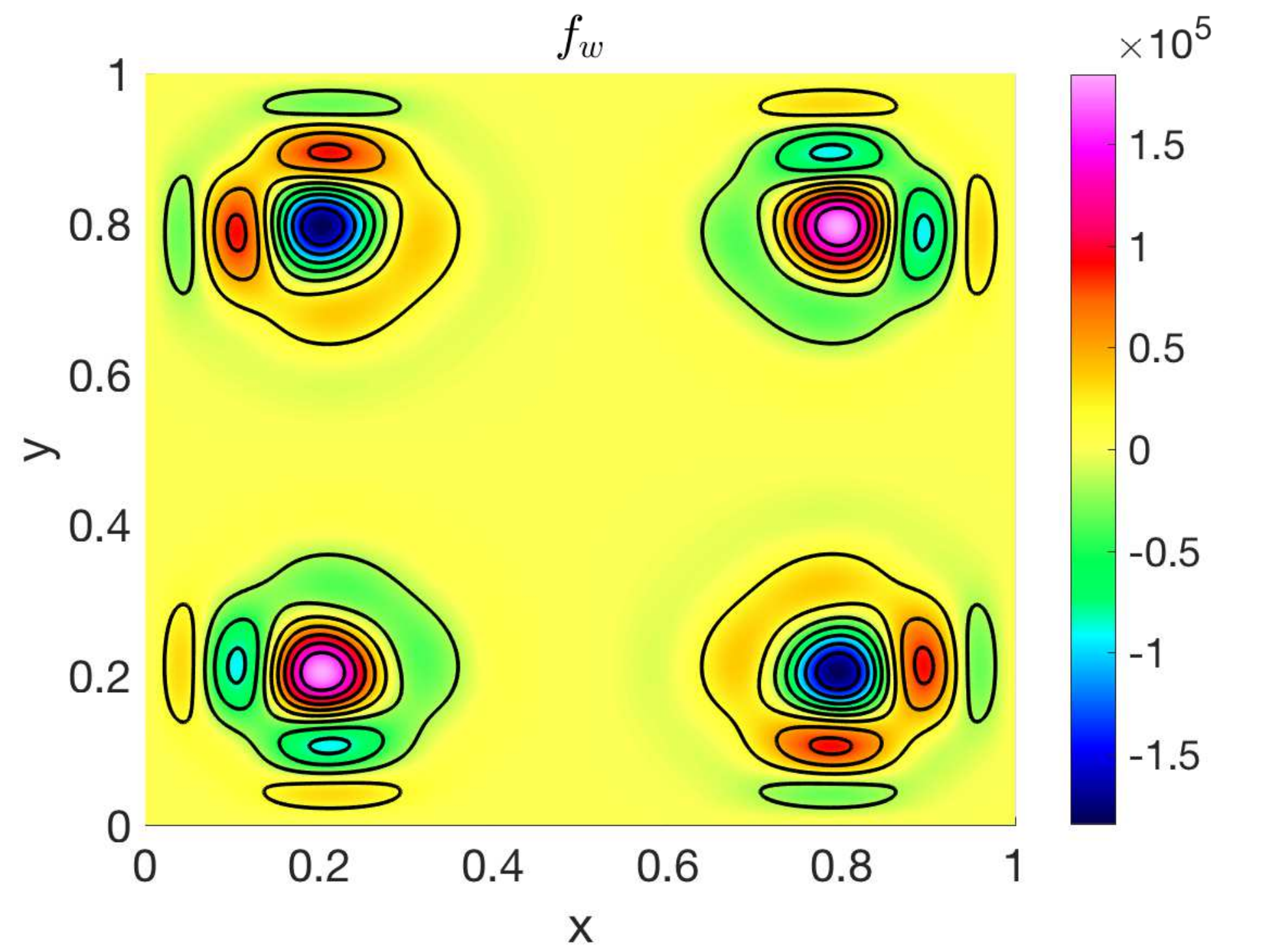}{\figWidth}};
%
\end{tikzpicture}
\caption{Plots of the forcing term $f_w$ in the biharmonic equation \eqref{eq:biharmEqn} for (left) the trigonometric test and (right) the polynomial test.}\label{fig:biharmTestForcingContour}
\end{center}
\end{figure}
}

{
\newcommand{\figWidth}{8cm}
\newcommand{\trimfig}[2]{\trimw{#1}{#2}{0.}{0.}{0.}{0.0}}
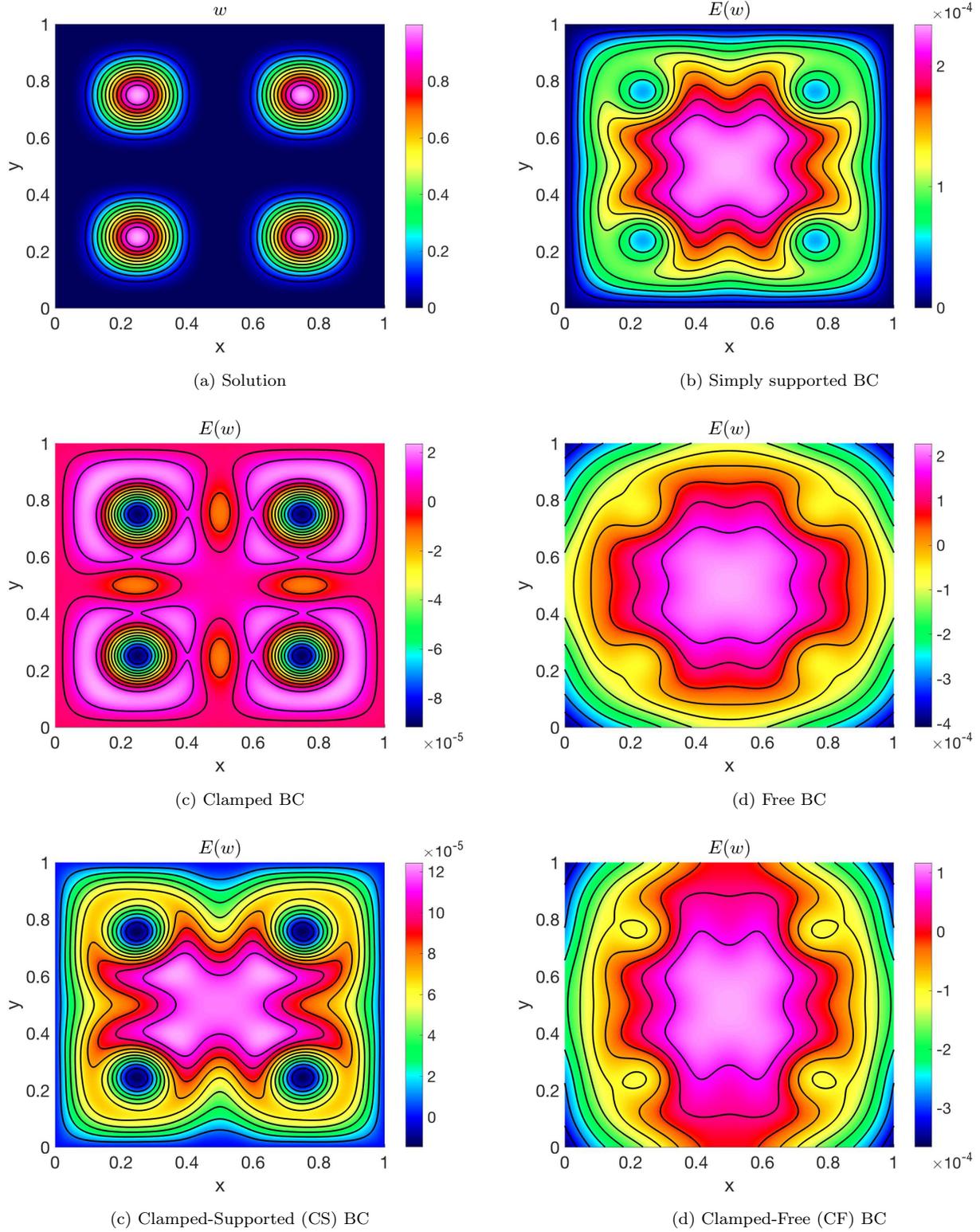
\begin{figure}[hp!]
\begin{center}
\begin{tikzpicture}[scale=1]
\useasboundingbox (0.0,0.0) rectangle (17.,21);  
\draw(-0.5,15.0) node[anchor=south west,xshift=0pt,yshift=0pt] {\trimfig{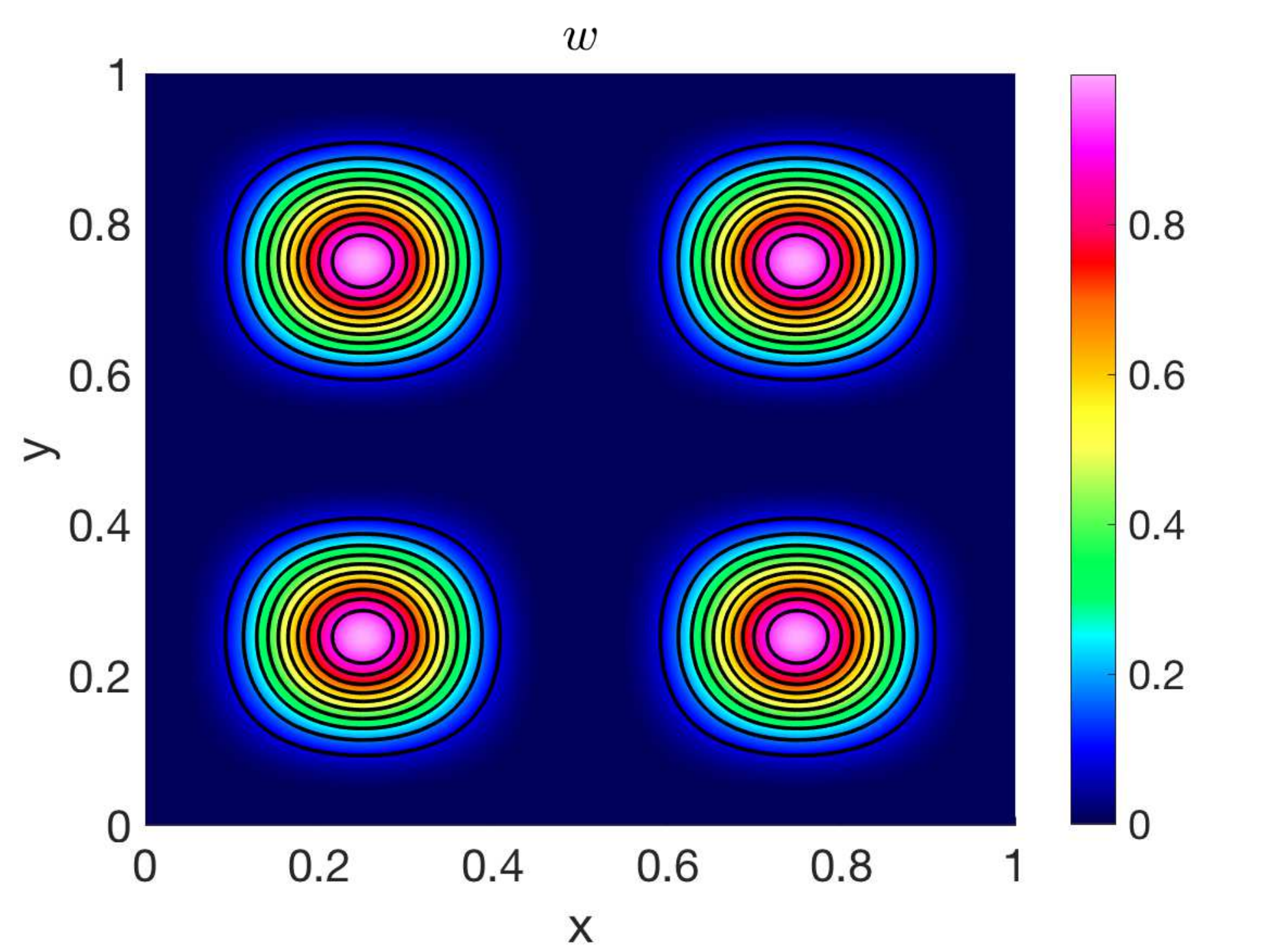}{\figWidth}};
\draw(8.0,15.0) node[anchor=south west,xshift=0pt,yshift=0pt] {\trimfig{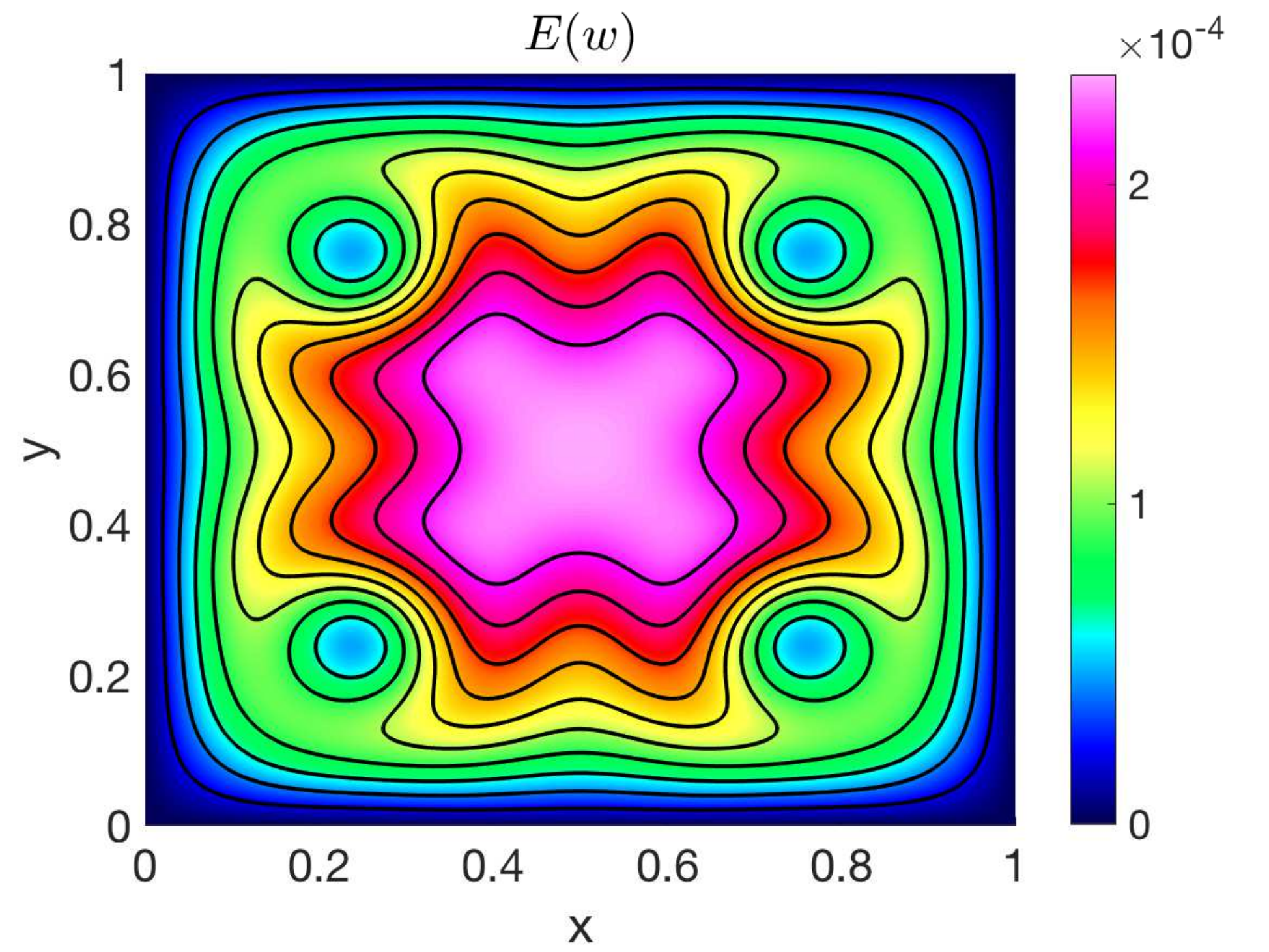}{\figWidth}};
\draw(-0.5,8.0) node[anchor=south west,xshift=0pt,yshift=0pt] {\trimfig{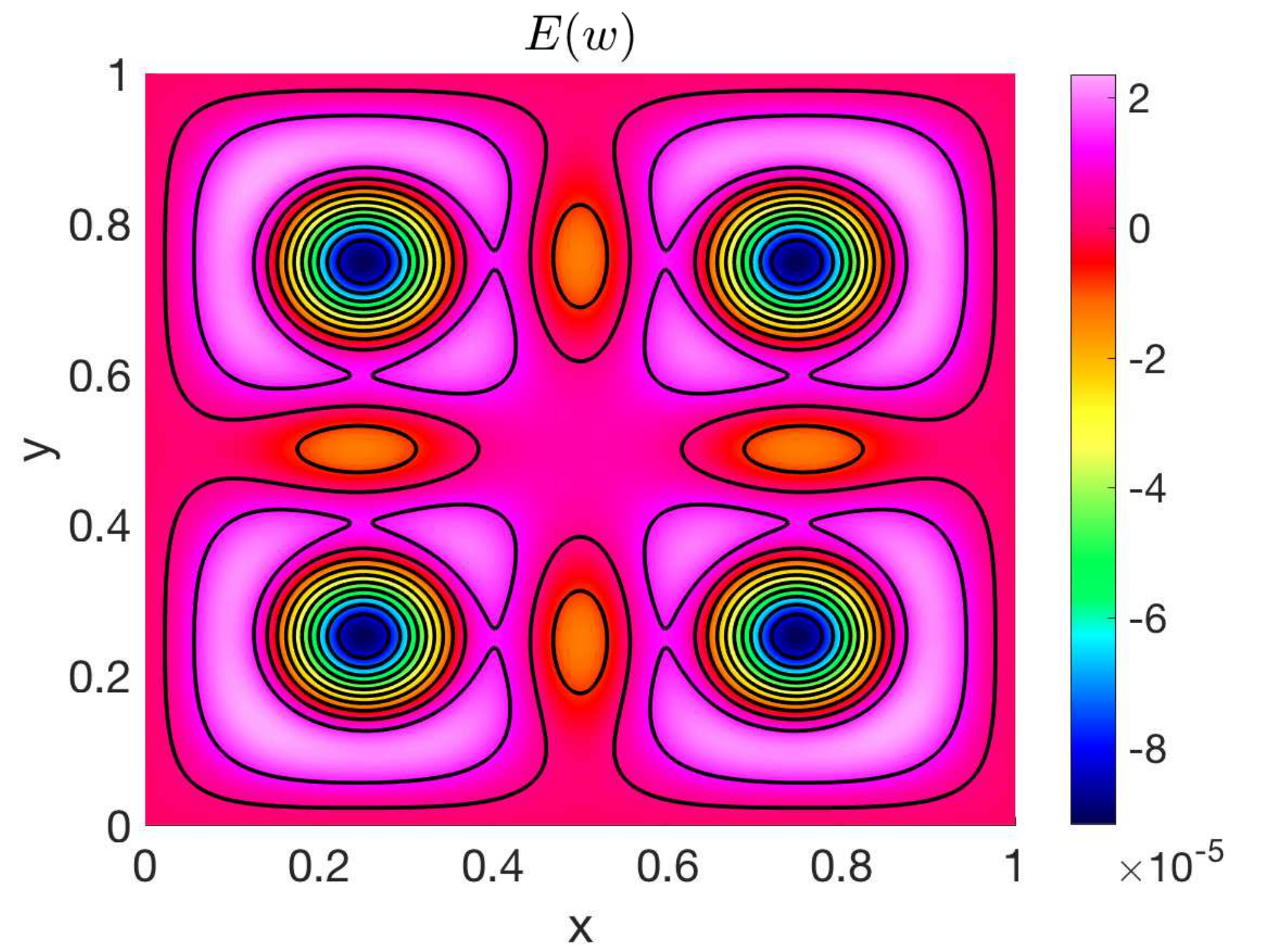}{\figWidth}};
\draw(8.0,8.0) node[anchor=south west,xshift=0pt,yshift=0pt] {\trimfig{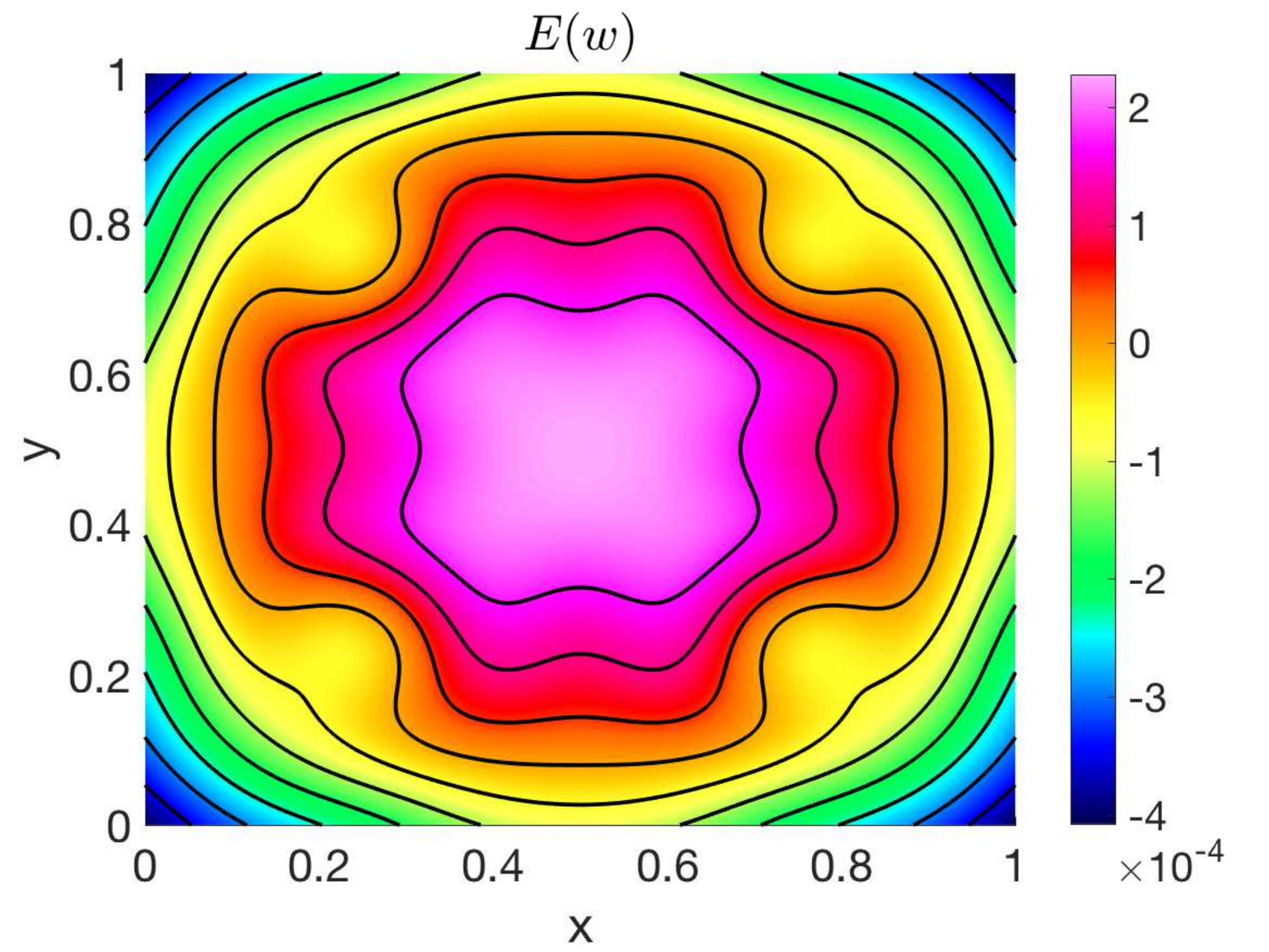}{\figWidth}};
\draw(-0.5,1.0) node[anchor=south west,xshift=0pt,yshift=0pt] {\trimfig{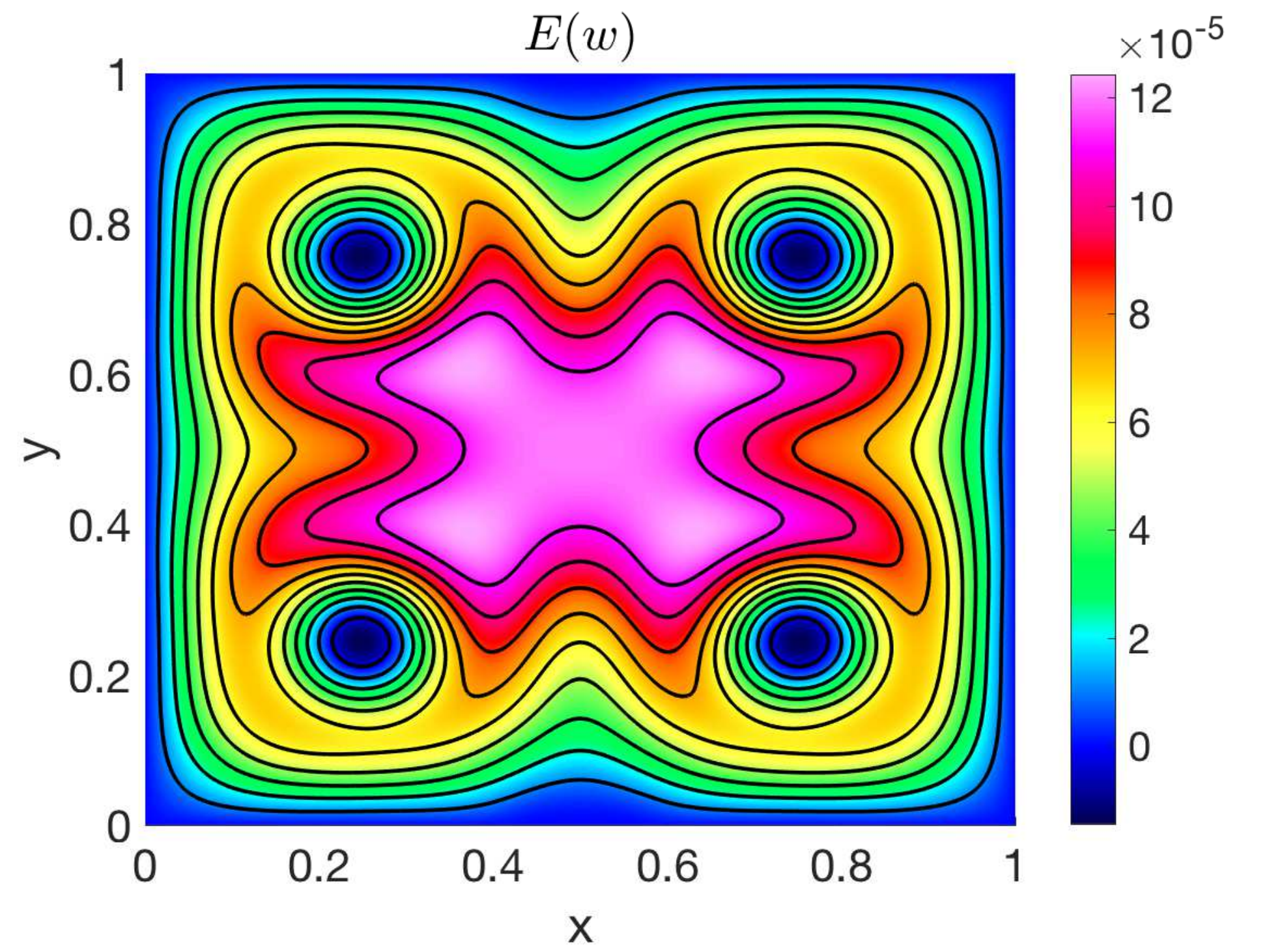}{\figWidth}};
\draw(8.0,1.0) node[anchor=south west,xshift=0pt,yshift=0pt] {\trimfig{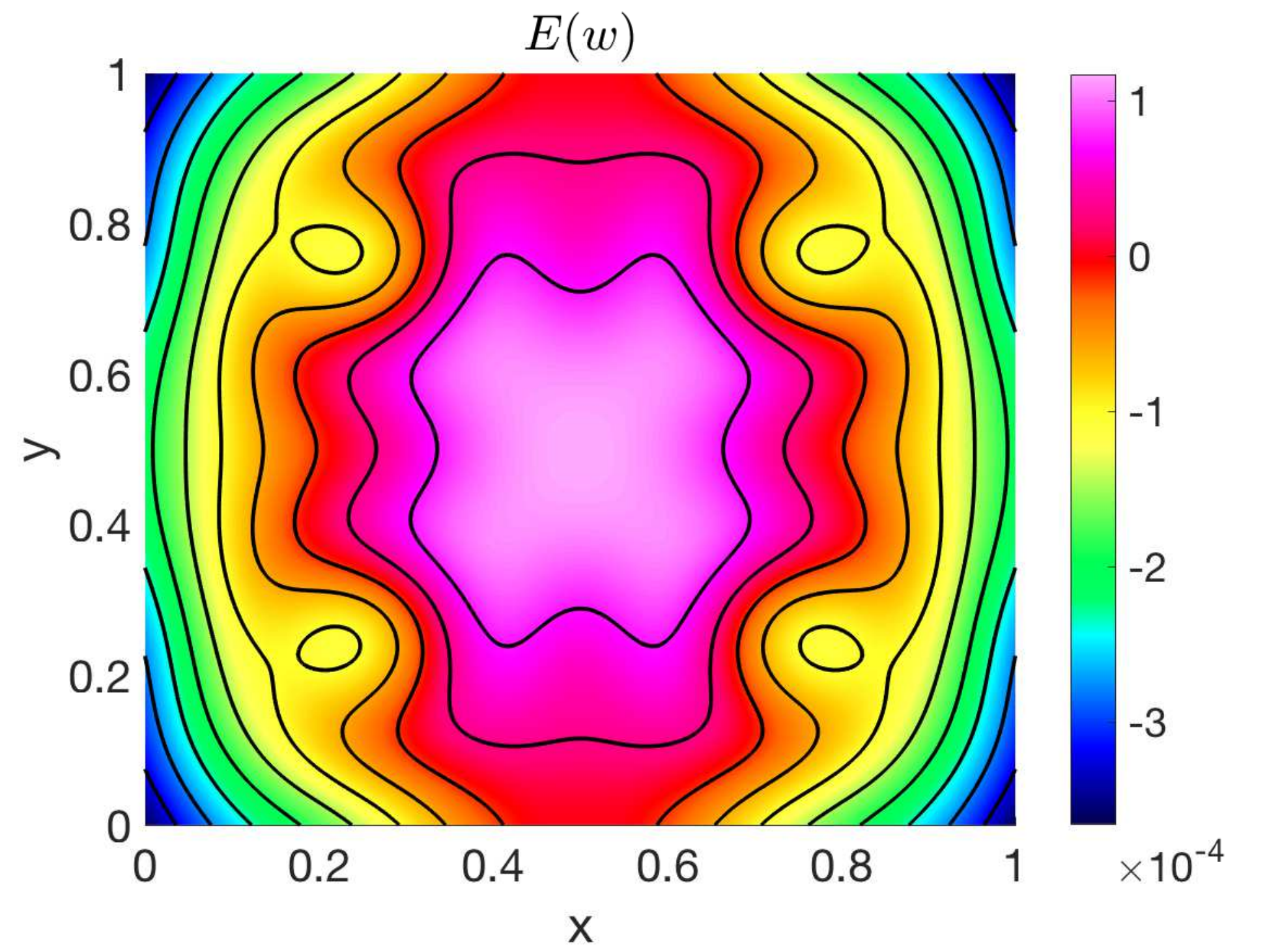}{\figWidth}};

\draw(4,15) node[anchor=north]{\footnotesize(a) Solution};
\draw(13,15) node[anchor=north]{\footnotesize(b) Simply supported BC};
\draw(4,8) node[anchor=north]{\footnotesize(c) Clamped BC};
\draw(13,8) node[anchor=north]{\footnotesize(d) Free BC};
\draw(4,1) node[anchor=north]{\footnotesize(c) Clamped-Supported (CS) BC};
\draw(13,1) node[anchor=north]{\footnotesize(d) Clamped-Free (CF) BC};
%
\end{tikzpicture}
\caption{Contour plots showing the solution and errors of the biharmonic equation \eqref{eq:biharmEqn} with various boundary conditions  on grid $\mathcal{G}_{640}$ for the trigonometric test.}\label{fig:biharmTrigTestResultContour}
\end{center}
\end{figure}
}

The accuracy of the biharmonic solution  is illustrated in Fig.~\ref{fig:biharmTrigTestResultContour} for all the boundary conditions. The first  plot  in the panel shows the numerical solution of the trigonometric test; the rest plots demonstrate the numerical error of  various boundary conditions. The error at grid $\iv$ is given by $E(w_\iv)=w_e(\xv_\iv)-W_{\iv}$. Here we observe that the errors for all the boundary conditions are well behaved in that the magnitude is small and is smooth throughout the domain including the boundaries. The behavior of the errors in the polynomial test are similar.

{
\newcommand{\figWidth}{8cm}
\newcommand{\trimfig}[2]{\trimw{#1}{#2}{0.}{0.}{0.}{0.0}}
\begin{figure}[h!]
\begin{center}
\begin{tikzpicture}[scale=1]
\useasboundingbox (0.0,0.0) rectangle (17.,7);  
\draw(-0.5,0.0) node[anchor=south west,xshift=0pt,yshift=0pt] {\trimfig{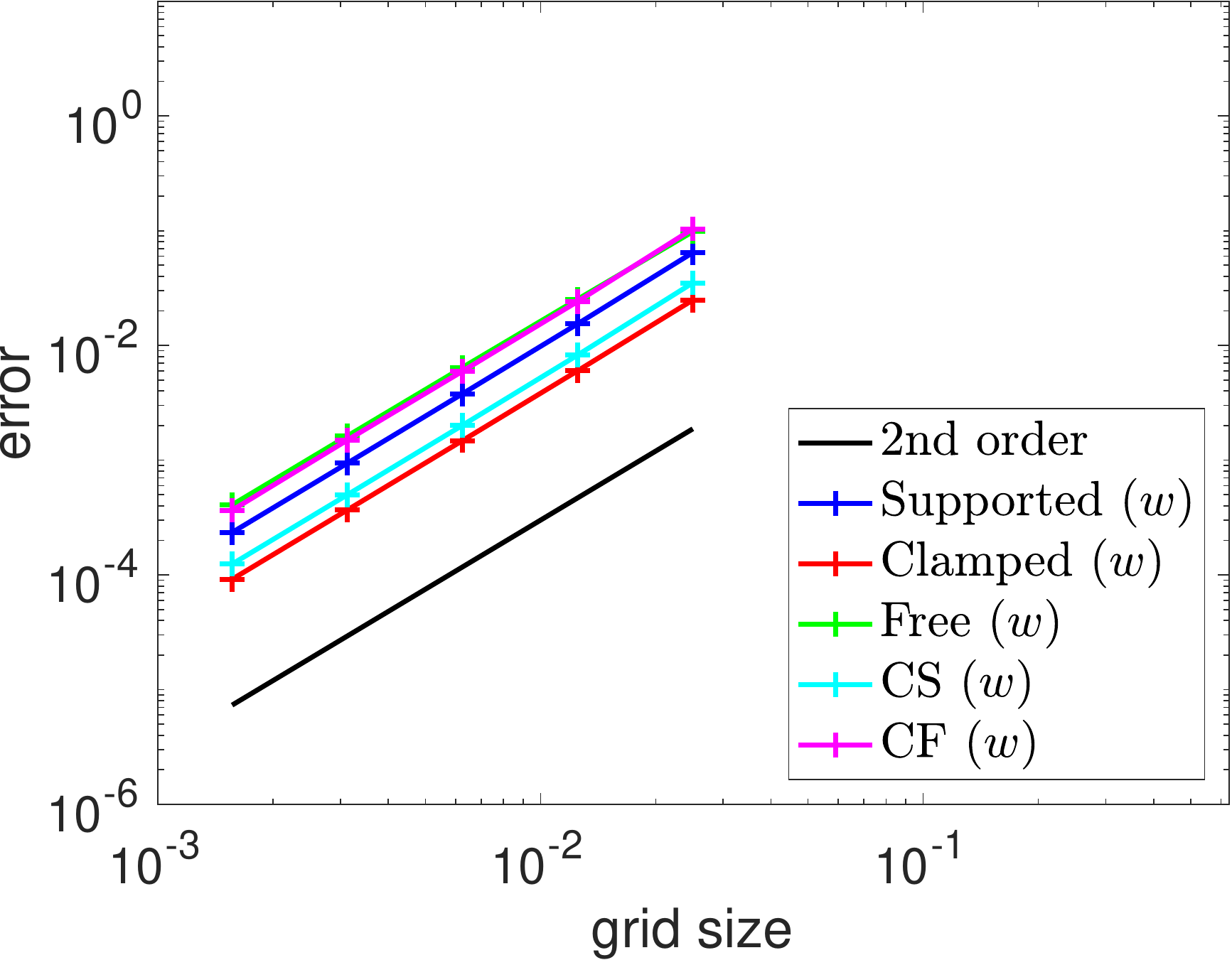}{\figWidth}};
\draw(8.0,0.0) node[anchor=south west,xshift=0pt,yshift=0pt] {\trimfig{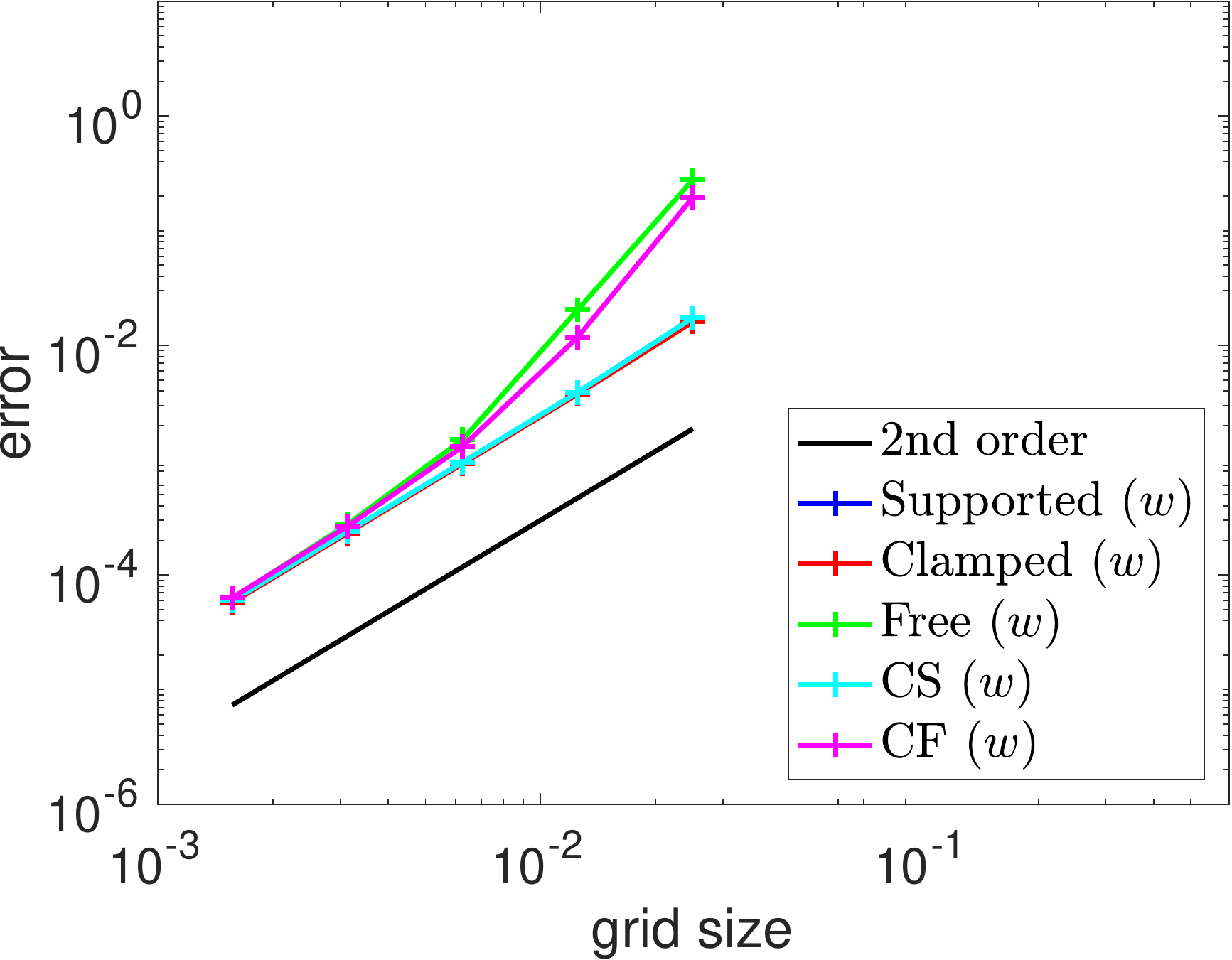}{\figWidth}};
\draw(4,7) node[anchor=north]{Trigonometric test};
\draw(13,7) node[anchor=north]{Polynomial test};
%
\end{tikzpicture}
\caption{A mesh refinement study for the numerical solutions of the biharmonic equation ($L_{\infty}$ norm).}\label{fig:convRateBiharmLinfty}
\end{center}
\end{figure}
}

A mesh refinement study is shown in Fig.~\ref{fig:convRateBiharmLinfty} for all the boundary conditions using the grids $\G_{N}$ defined  in \eqref{eq:gridGN} where $N = 10\times 2^j$ with $j = 1, 2, \cdots, 6$. The maximum-norm errors $||E(w)||_\infty$ against the grid size  for all boundary conditions together with a second-order reference curve are plotted in log-log scale  in Fig.~\ref{fig:convRateBiharmLinfty}. The results in the plot show the expected second-order accuracy for all five choices of  boundary conditions, and in particular for the case of free boundary conditions; the problem would otherwise be singular  without the technique introduced in \S\ref{sec:removeSingularityFreeBC}.  The accuracy result is consistent with the truncation error of our centered finite difference discretization.


\subsubsection{Asymptotic analytical solution approach for mixed boundary conditions}
As an example to show the convergence property of the asymptotic analytical solution approach for mixed boundary conditions, we consider a simple  test problem consisting of the  biharmonic equation \eqref{eq:biharmEqn} with  a constant external forcing $f_w\equiv1$  on the unit square domain   $[0,1]\times[0,1]$. The square shell is assumed to be partially clamped on the  right half of its boundary, i.e., $\Gamma_c = \{{y=0,~0.5\le x\le 1}\}\cup\{y=1,~0.5\le x \le 1\}\cup\{x=1,~0\le y\le 1\}$, and the rest of the boundary is either supported or free. This problem without boundary singularities  corresponds to the plate sagging problem under gravity, which has been well-studied in classic mechanical engineering textbooks \cite{howell2009applied}. 

We first consider the CS mixed boundary conditions; that is, the rest of the boundary is simply supported.  
The reduced rigidity caused by  simply supported boundary conditions enforced on the left half of the boundary should lead to larger displacements at the center of the plate compared to the case associated with fully clamped boundary conditions \cite{howell2009applied}.
The contour of the displacement function $w(x,y)$ shown in the top-left image of  Fig.~\ref{CS_biharmonic_halfdomain} agrees with this result, where the displacement in the clamped half is significantly smaller than that in the simply supported half.  We then consider the CF mixed boundary conditions. The numerical solution for the displacement $w$ is shown in the top-right image of Fig.~\ref{CS_biharmonic_halfdomain}. The free boundary condition on the left half of the boundary has an obvious effect on the sagging of the plate; the  lowest point now locates on the free edge instead of at the center of the domain.


{
\newcommand{\figWidth}{8cm}
\newcommand{\trimfig}[2]{\trimw{#1}{#2}{0.}{0.}{0.}{0.0}}
\def\xa{17.}
\def\ya{13}
\begin{figure}[h!]
\begin{center}
\begin{tikzpicture}[scale=1]
\useasboundingbox (0.0,0.0) rectangle (\xa,\ya);  
\draw(-0.5,6.5) node[anchor=south west,xshift=0pt,yshift=0pt] {\trimfig{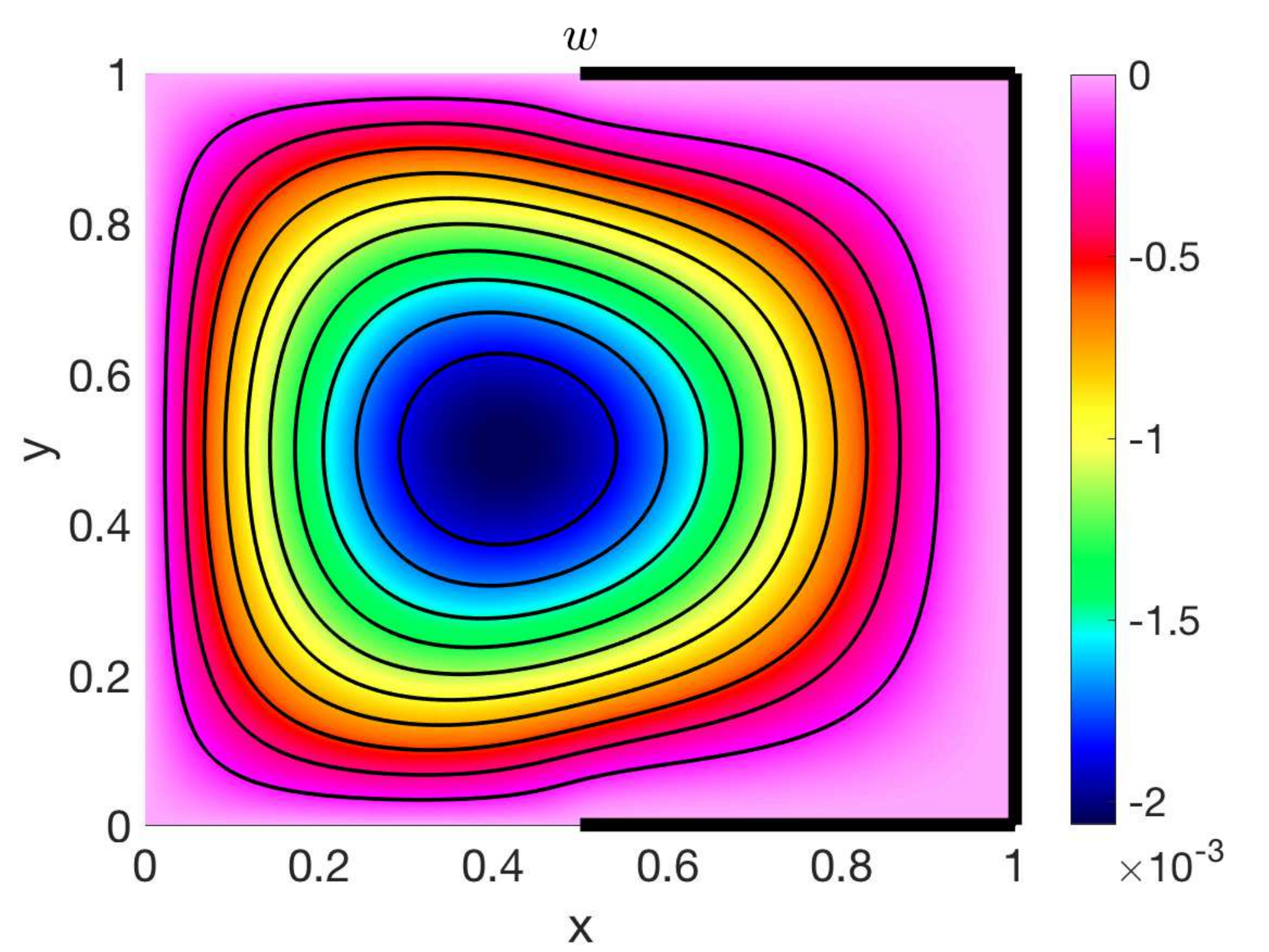}{\figWidth}};
\draw(8.0,6.5) node[anchor=south west,xshift=0pt,yshift=0pt] {\trimfig{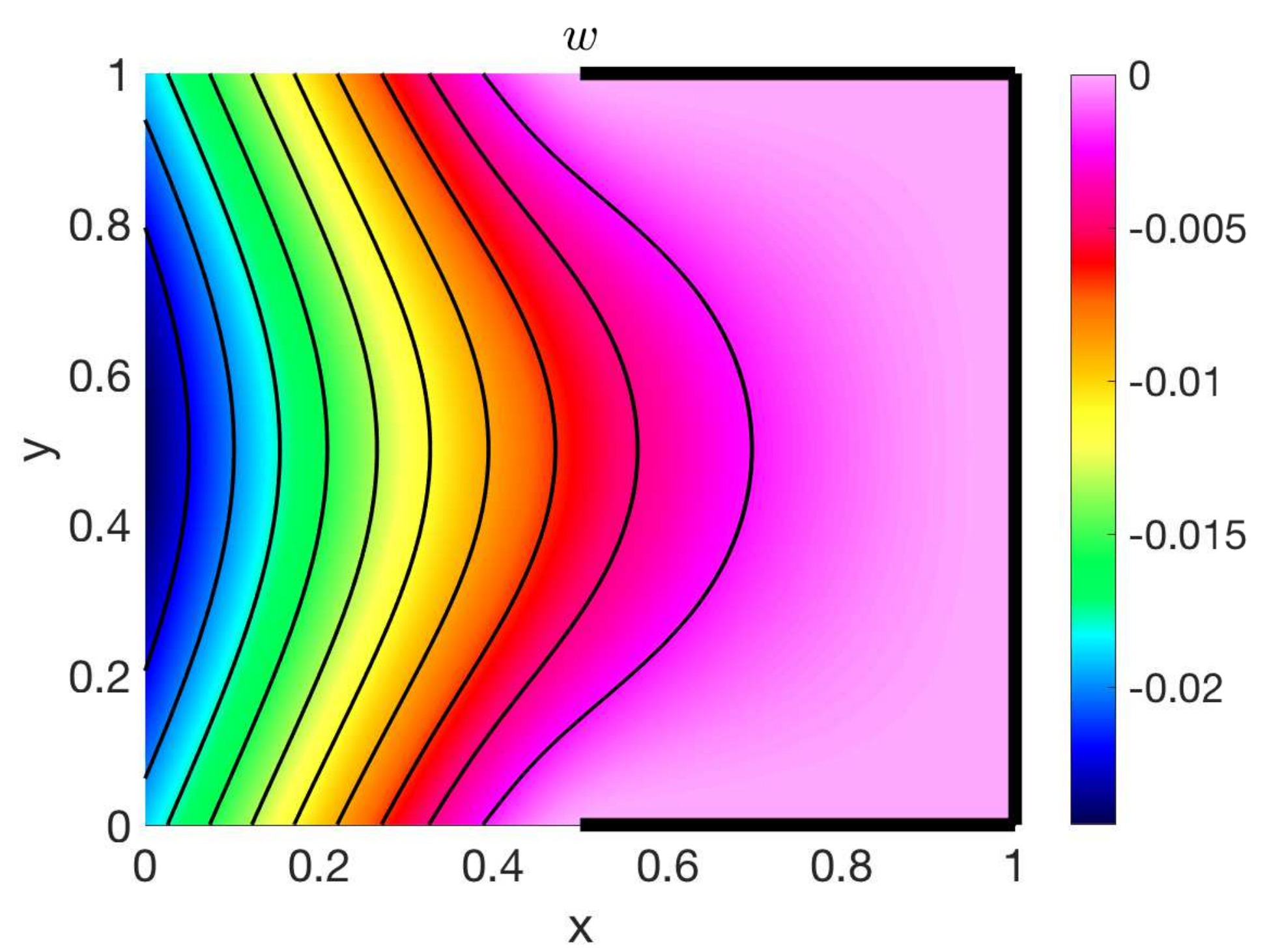}{\figWidth}};
\draw(4.0,0.0) node[anchor=south west,xshift=0pt,yshift=0pt] {\trimfig{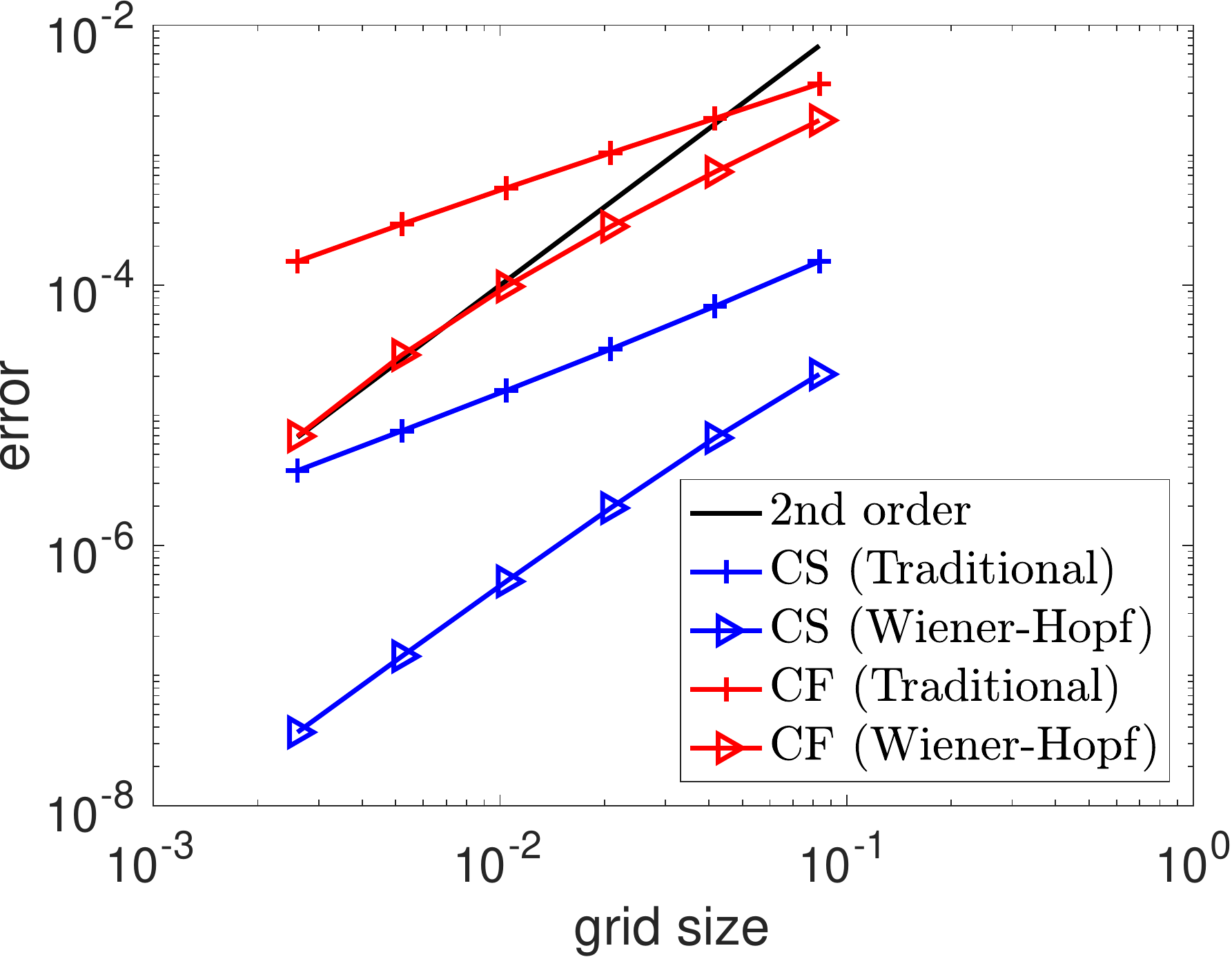}{\figWidth}};
%
\end{tikzpicture}
\caption{Numerical results for the test problem $\nabla^4 w = -1$  subject to partially clamped boundary conditions on the right portion of the edges (marked with solid lines). The solution for the CS mixed boundary conditions is shown in the top-left image, and that for the    CF mixed boundary conditions is shown in the  top-right image. The bottom image shows the  mesh refinement studies for  the  local singularity treatment \eqref{constraint_general} and the traditional finite difference schemes with both CS and CF mixed boundary conditions.}
\label{CS_biharmonic_halfdomain}
\end{center}
\end{figure}
}


Mesh refinement study is also performed for this test problem to reveal the order of accuracy of the local asymptotic solution approach for dealing with mixed boundary conditions. For comparison purposes, we also solve the test problem using traditional finite difference methods; namely, no special treatment are given to the singular point on boundary.   The convergence results for both mixed boundary conditions are shown in the bottom image of  Fig.~\ref{CS_biharmonic_halfdomain}.
We observe for both CS and CF boundary conditions that the traditional method  has a linear rate of convergence; while with the flexible local approximation scheme in \eqref{constraint_general}  applied to the inner node adjacent to the singular point, the resultant numerical scheme exhibits a second-order convergence rate.


We can see that both the transition function approach and the local asymptotic solution approach
are effective in maintaining the second order accuracy for solving the biharmonic equation with the mixed boundary conditions. Given that the implementation is much straightforward for the  transition function approach, we stick with this method for the rest of the paper.


\subsubsection{Linear coupled system}
To illustrate the effectiveness of the iterative schemes for solving the coupled system, we now discuss the problem of linear shallow equations \eqref{eq:coupledSystemLinear} which is obtained by dropping the nonlinear terms in the nonlinear shallow shaw shell equations \eqref{eq:coupledSystemNonlinear}.
This linear coupled system is considered here because it is simple and applicable for shallow shells with small deformations (the linear shallow shell theory). Moreover,  the  solution to this linear system can be utilized as an initial guess to the iterative process of solving the nonlinear shell equations.

Again, the method of manufactured solution is used here. The exact solutions chosen for this test are given by 
\begin{subequations}
\label{eq:exactSolutionCoupled}
\begin{equation}
\phi_e(x,y) = \sin^5\left(2\pi\frac{x-x_a}{x_b-x_a}\right)\sin^5\left(2\pi\frac{y-y_a}{y_b-y_a}\right),\label{eq:exactSolutionCoupled_phi}
\end{equation}
\begin{equation}
w_e(x,y) = \sin^4\left(2\pi\frac{x-x_a}{x_b-x_a}\right)\sin^4\left(2\pi\frac{y-y_a}{y_b-y_a}\right).\label{eq:exactSolutionCoupled_w}
\end{equation}
\end{subequations}
All five boundary conditions are satisfied by the exact solutions. The precast shell shape is specified as
\begin{equation}\label{eq:precastShapeCoupled}
w_0(x,y) = \sin\left(2\pi\frac{x-x_a}{x_b-x_a}\right)\sin\left(2\pi\frac{y-y_a}{y_b-y_a}\right).
\end{equation}
The forcing functions $f_\phi(x,y)$ and $f_w(x,y)$ are obtained accordingly by substituting $\phi_e$, $w_e$ and $w_0$ into the system, and their contour plots  are shown  in Fig.~\ref{fig:linearCoupledSystemForcingsContour}.
{
\newcommand{\figWidth}{8cm}
\newcommand{\trimfig}[2]{\trimw{#1}{#2}{0.}{0.}{0.}{0.0}}
\begin{figure}[h!]
\begin{center}
\begin{tikzpicture}[scale=1]
\useasboundingbox (0.0,0.0) rectangle (17.,6);  
\draw(-0.5,0.0) node[anchor=south west,xshift=0pt,yshift=0pt] {\trimfig{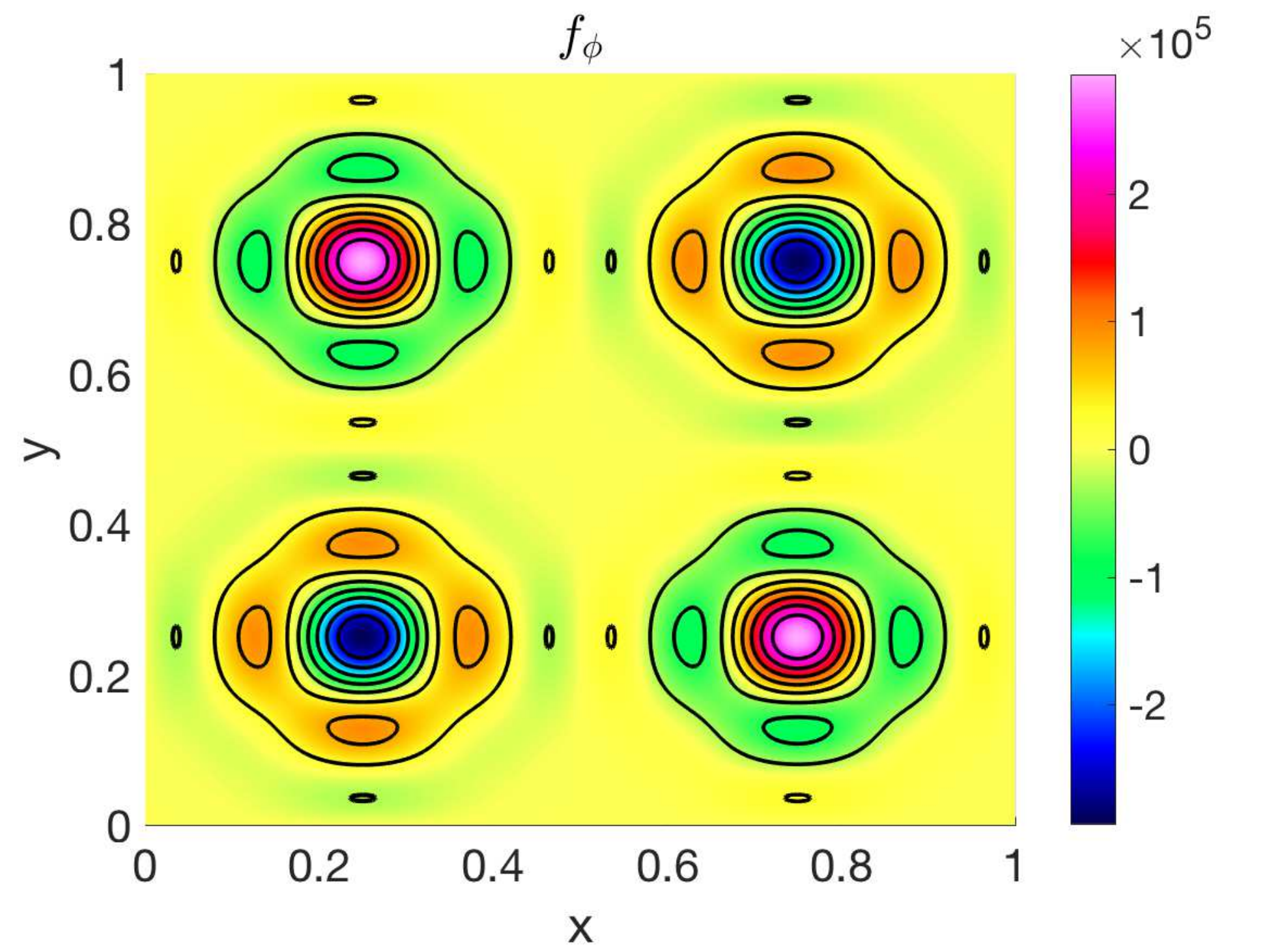}{\figWidth}};
\draw(8.0,0.0) node[anchor=south west,xshift=0pt,yshift=0pt] {\trimfig{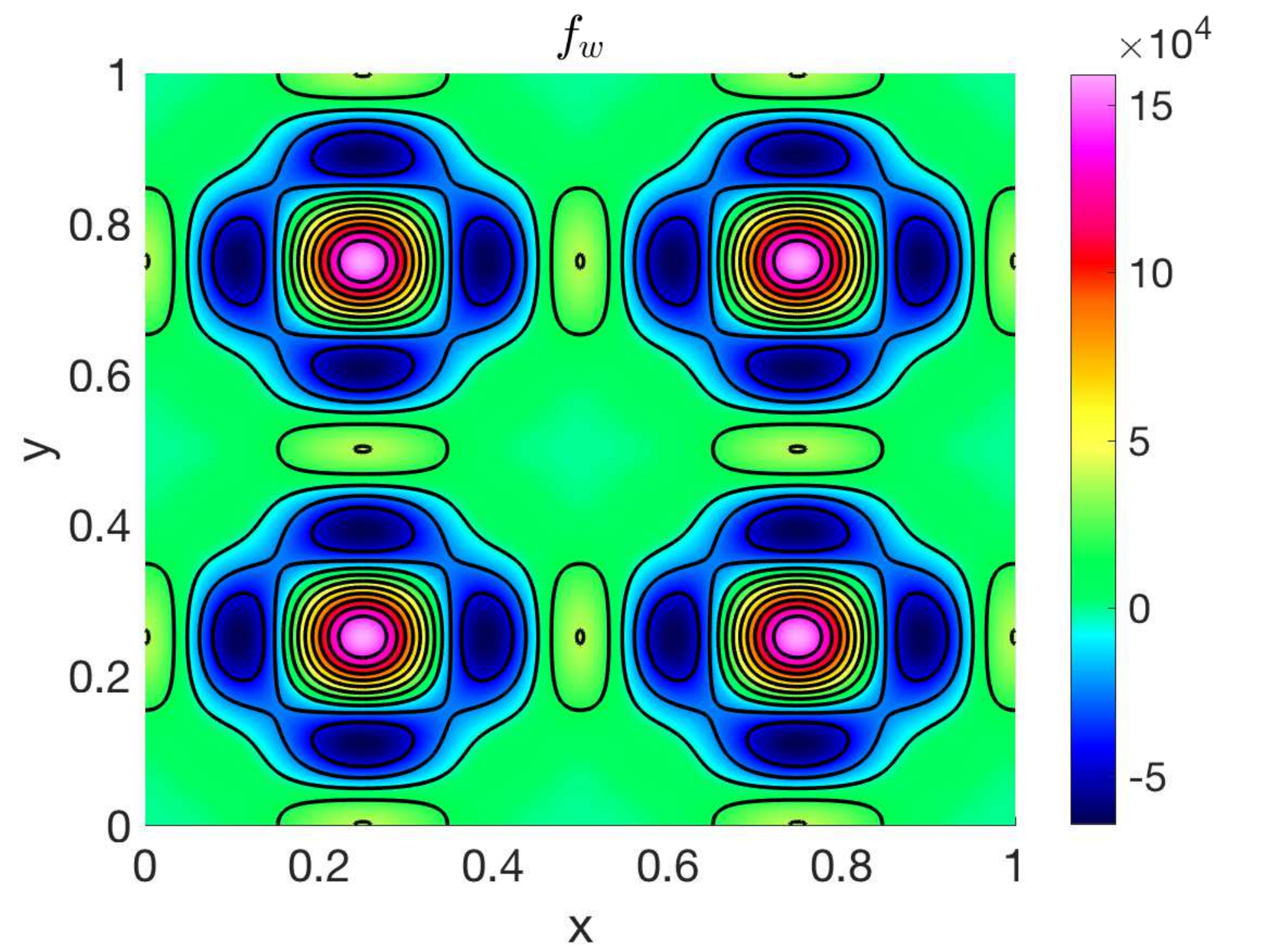}{\figWidth}};
%
\end{tikzpicture}
\caption{Contour plots of the forcing terms $f_{\phi}$ and $f_w$ given by \eqref{eq:exactSolutionCoupled} for the linear coupled system \eqref{eq:coupledSystemLinear}.}\label{fig:linearCoupledSystemForcingsContour}
\end{center}
\end{figure}
}

We solve the linear coupled system using all three iterative schemes together with all five boundary conditions for completeness. The initial conditions used to start the iteration is as discussed in Section \ref{sec:initialGuess}.
The results from the Picard method (Algorithm~\ref{alg:picardSolve}) with the implicit factor $\delta=0$ are presented  in this section.  Solutions obtained using the other iterative algorithms are similar and they are not included here to save space.   The results of the $\phi$ component  are collected in Fig.~\ref{fig:linearCoupledTestResultContour_phi}, and those of the $w$ component are collected in  Fig.~\ref{fig:linearCoupledTestResultContour_w}. We see that  the errors of both $\phi$ and $w$ components subject to all the five choices of boundary conditions  are well behaved; the errors are   small and smooth throughout the domain including the boundaries. 

A careful mesh refinement study is also performed to test the order of accuracy for the linear coupled system. The  series of refined  grids  are $\G_N$'s with  $N=10\times2^j ~ (j=1,2,\dots,6)$.  As expected, second-order spatial accuracy is achieved by all three algorithms as is shown in Fig.~\ref{fig:convRateLinearCoupledSystem}.  It is worth noting that the regularization for the $w$ equation with free boundary conditions works well for the linear coupled system regardless of the numerical methods used for iterations.

{
\newcommand{\figWidth}{8cm}
\newcommand{\trimfig}[2]{\trimw{#1}{#2}{0.}{0.}{0.}{0.0}}
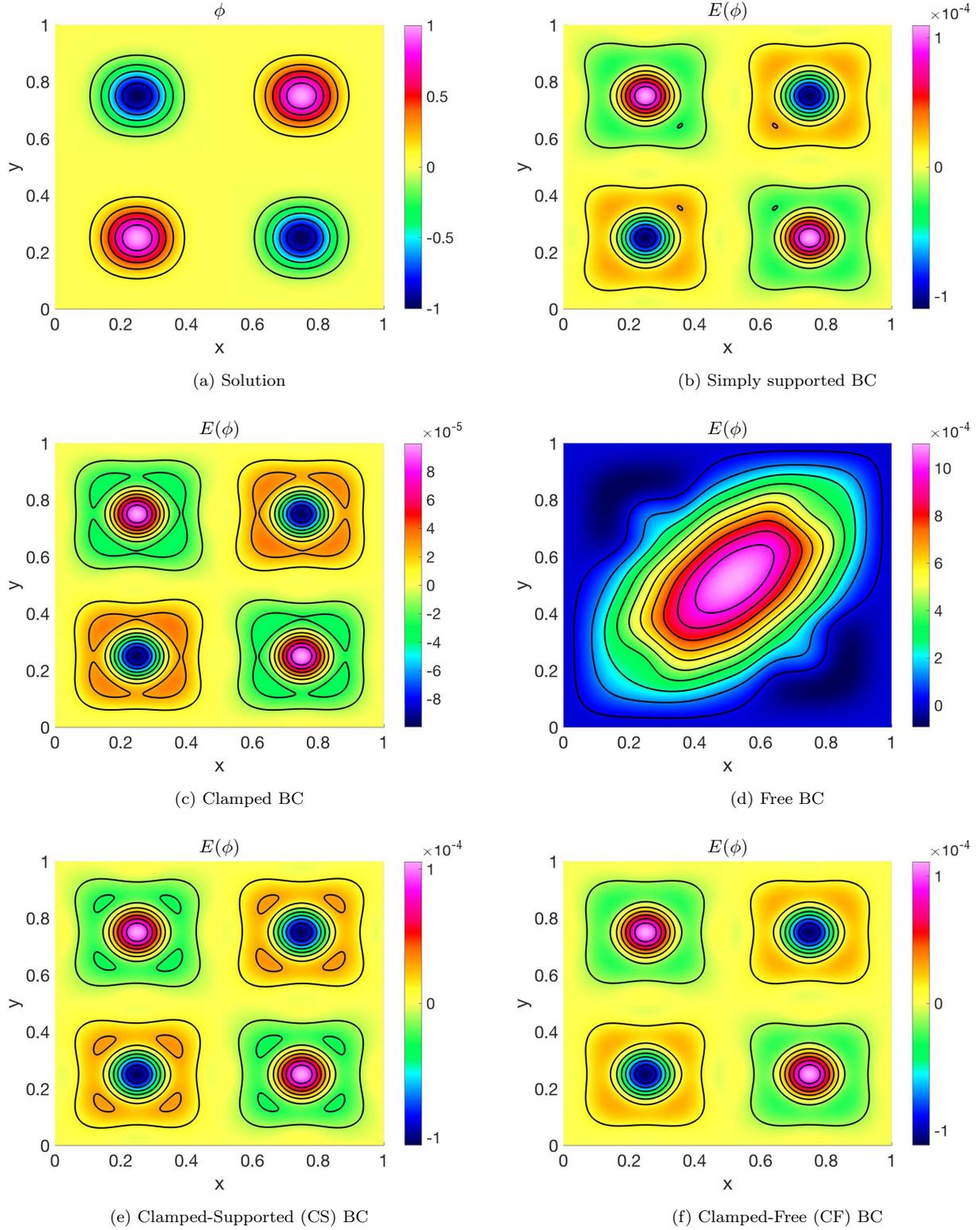
\begin{figure}[hp!]
\begin{center}
\begin{tikzpicture}[scale=1]
\useasboundingbox (0.0,0.0) rectangle (17.,21);  
\draw(-0.5,15.0) node[anchor=south west,xshift=0pt,yshift=0pt] {\trimfig{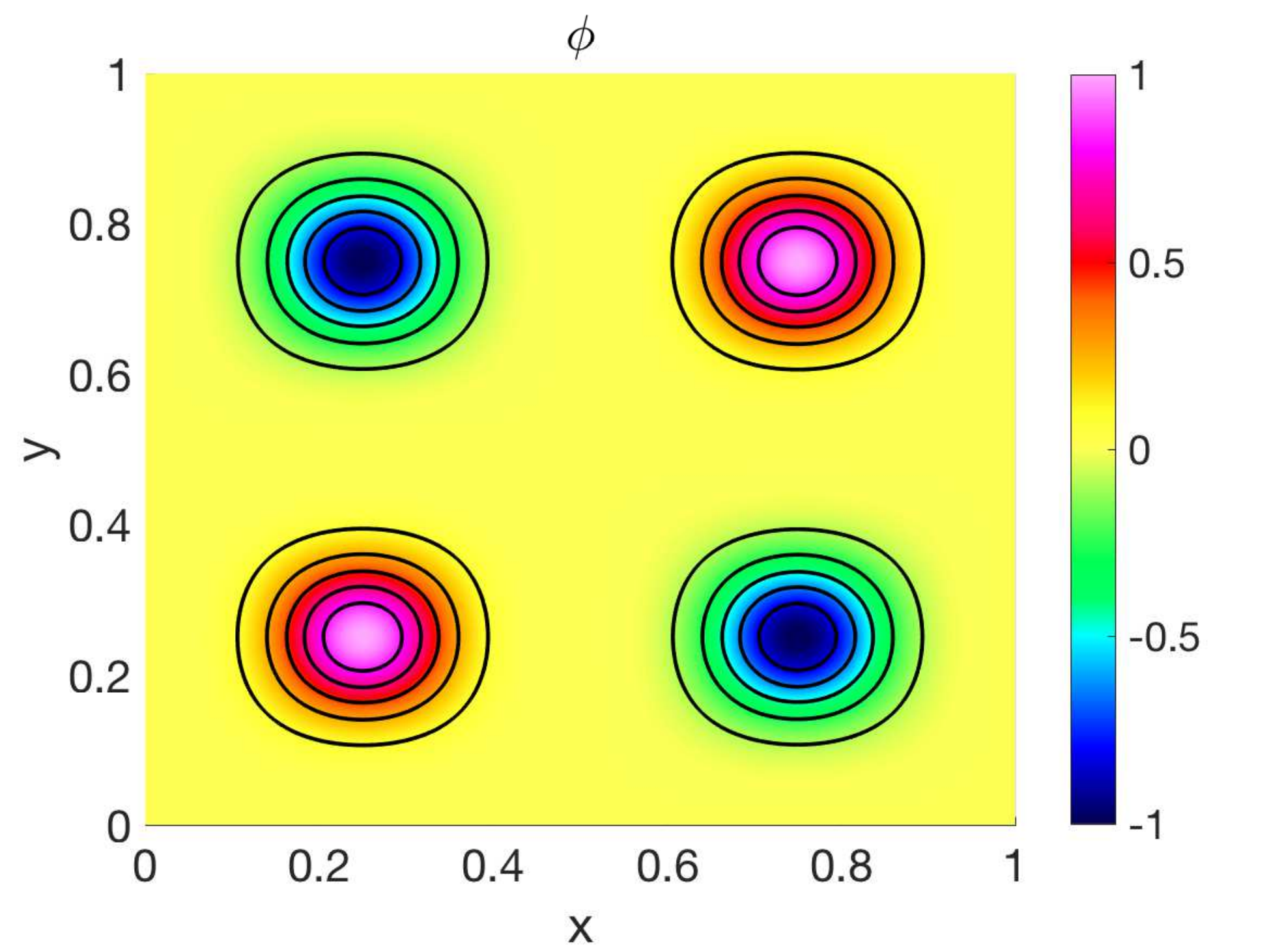}{\figWidth}};
\draw(8.0,15.0) node[anchor=south west,xshift=0pt,yshift=0pt] {\trimfig{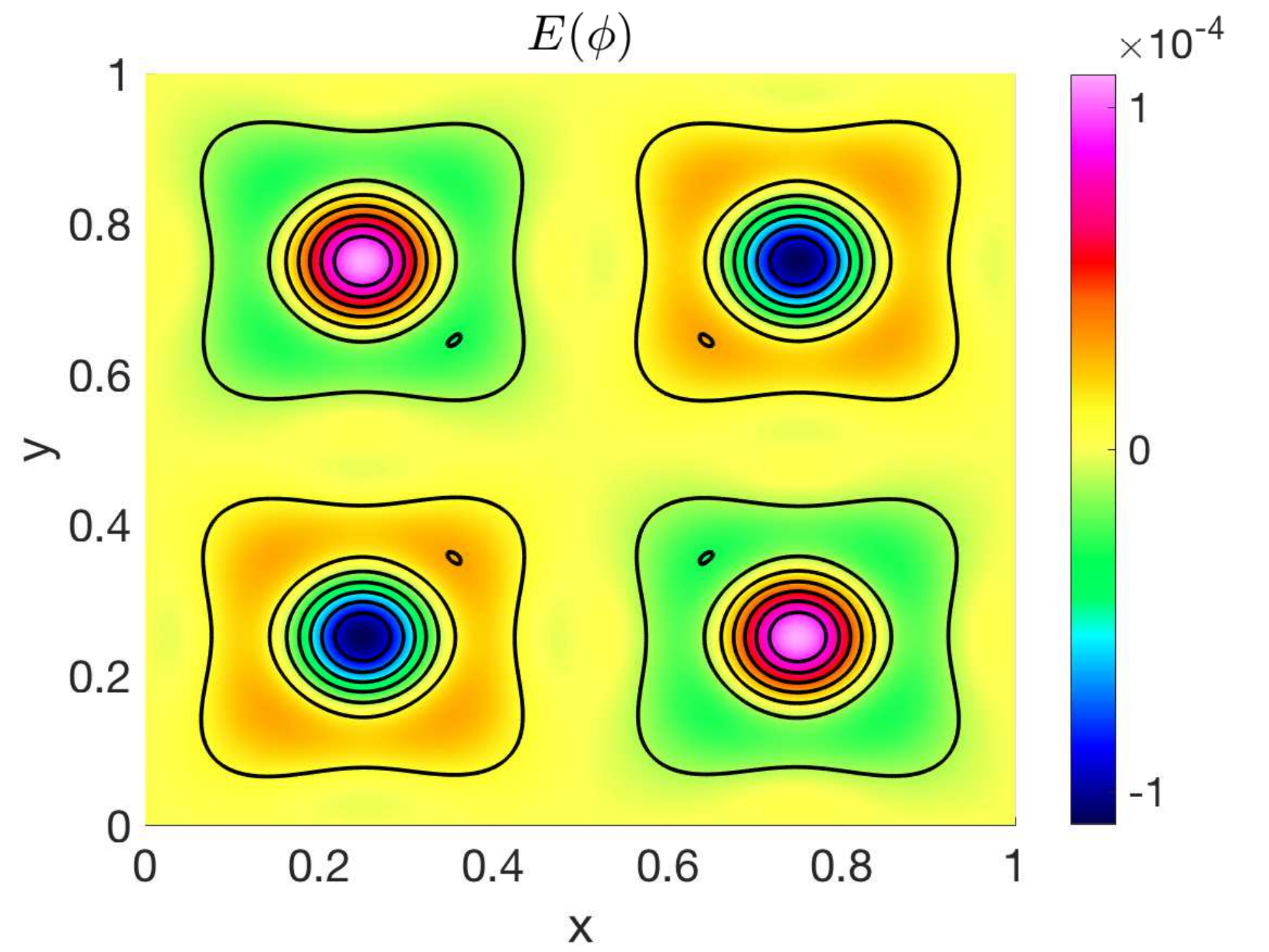}{\figWidth}};
\draw(-0.5,8.0) node[anchor=south west,xshift=0pt,yshift=0pt] {\trimfig{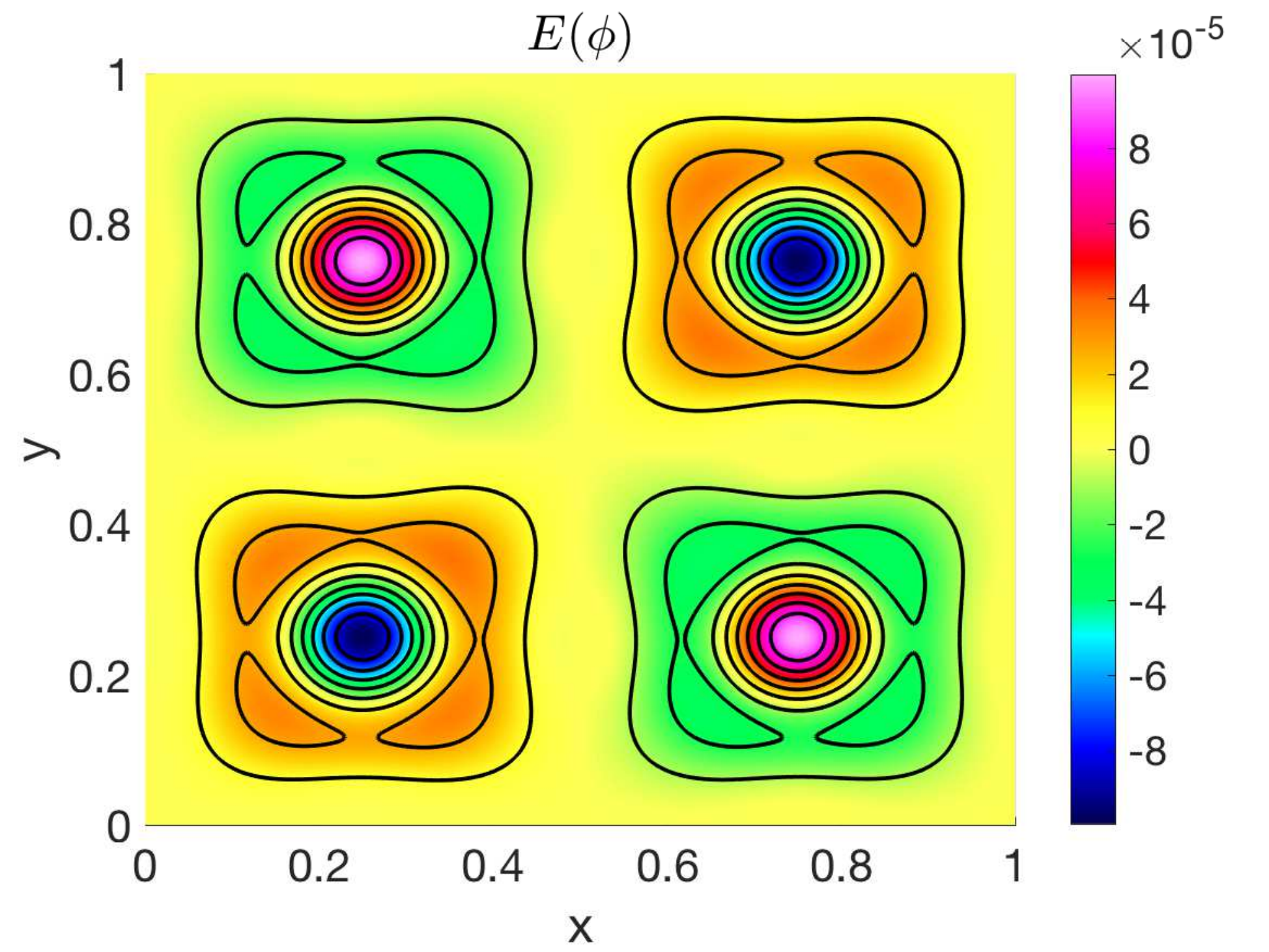}{\figWidth}};
\draw(8.0,8.0) node[anchor=south west,xshift=0pt,yshift=0pt] {\trimfig{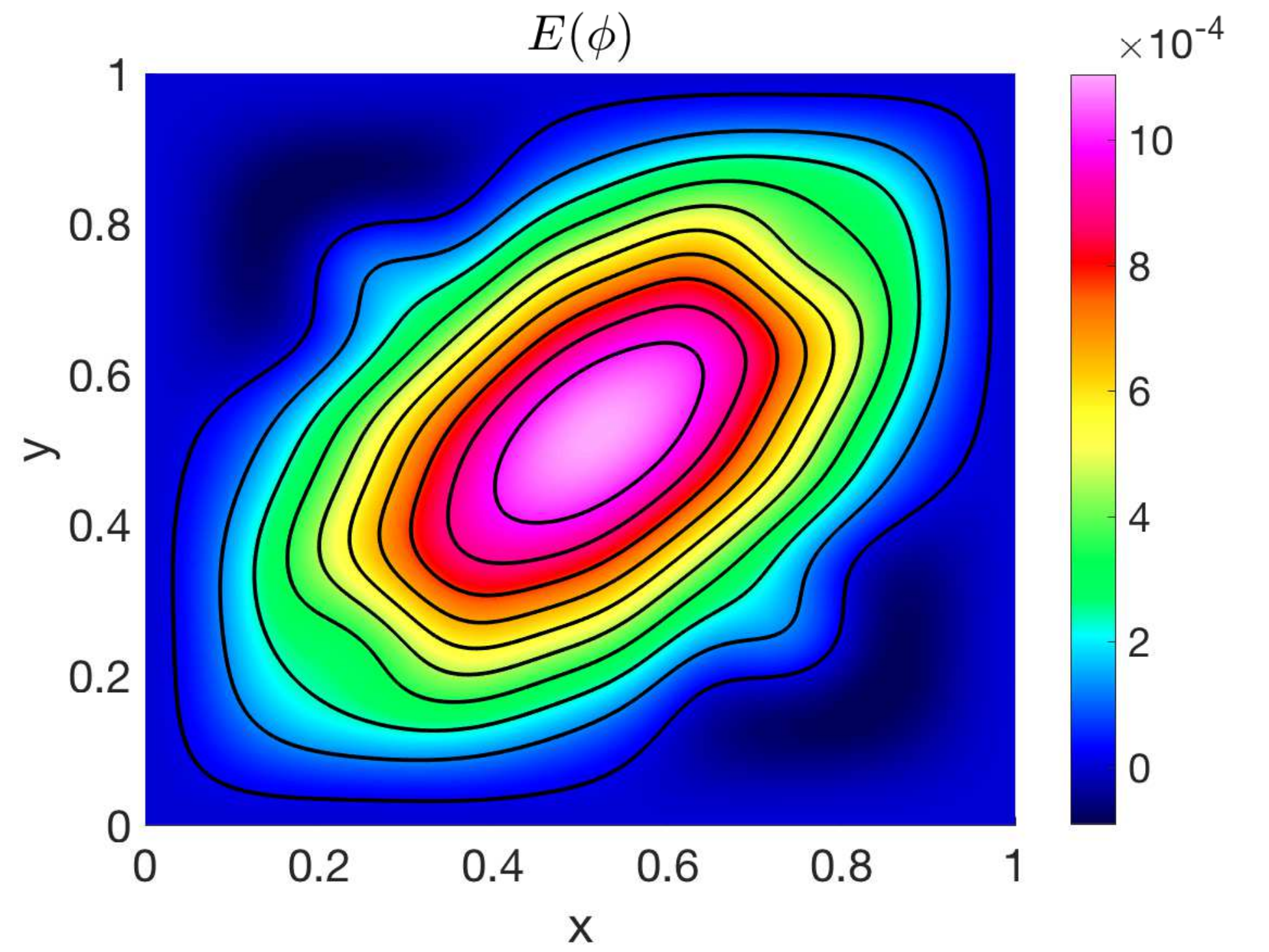}{\figWidth}};
\draw(-0.5,1.0) node[anchor=south west,xshift=0pt,yshift=0pt] {\trimfig{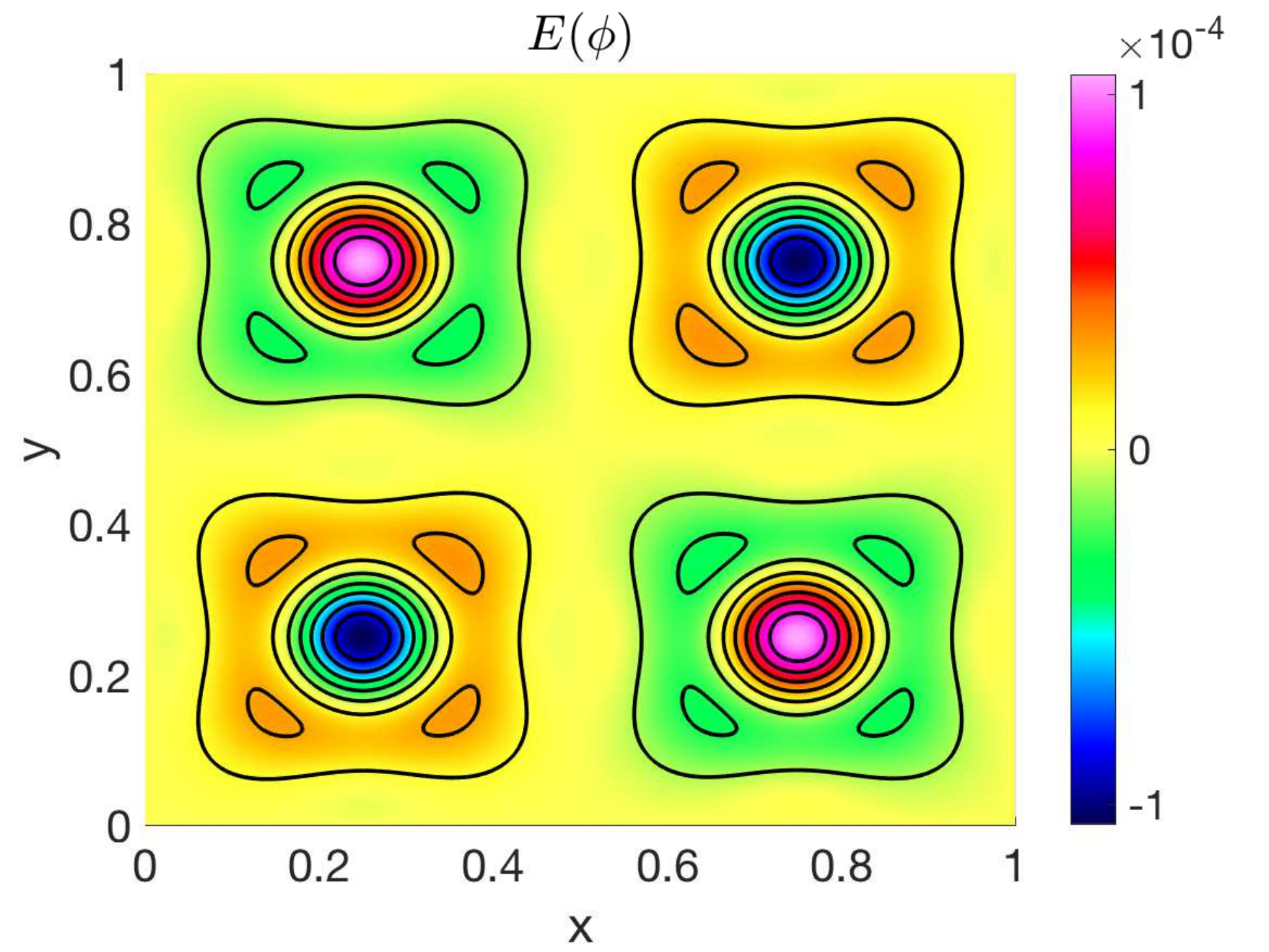}{\figWidth}};
\draw(8.0,1.0) node[anchor=south west,xshift=0pt,yshift=0pt] {\trimfig{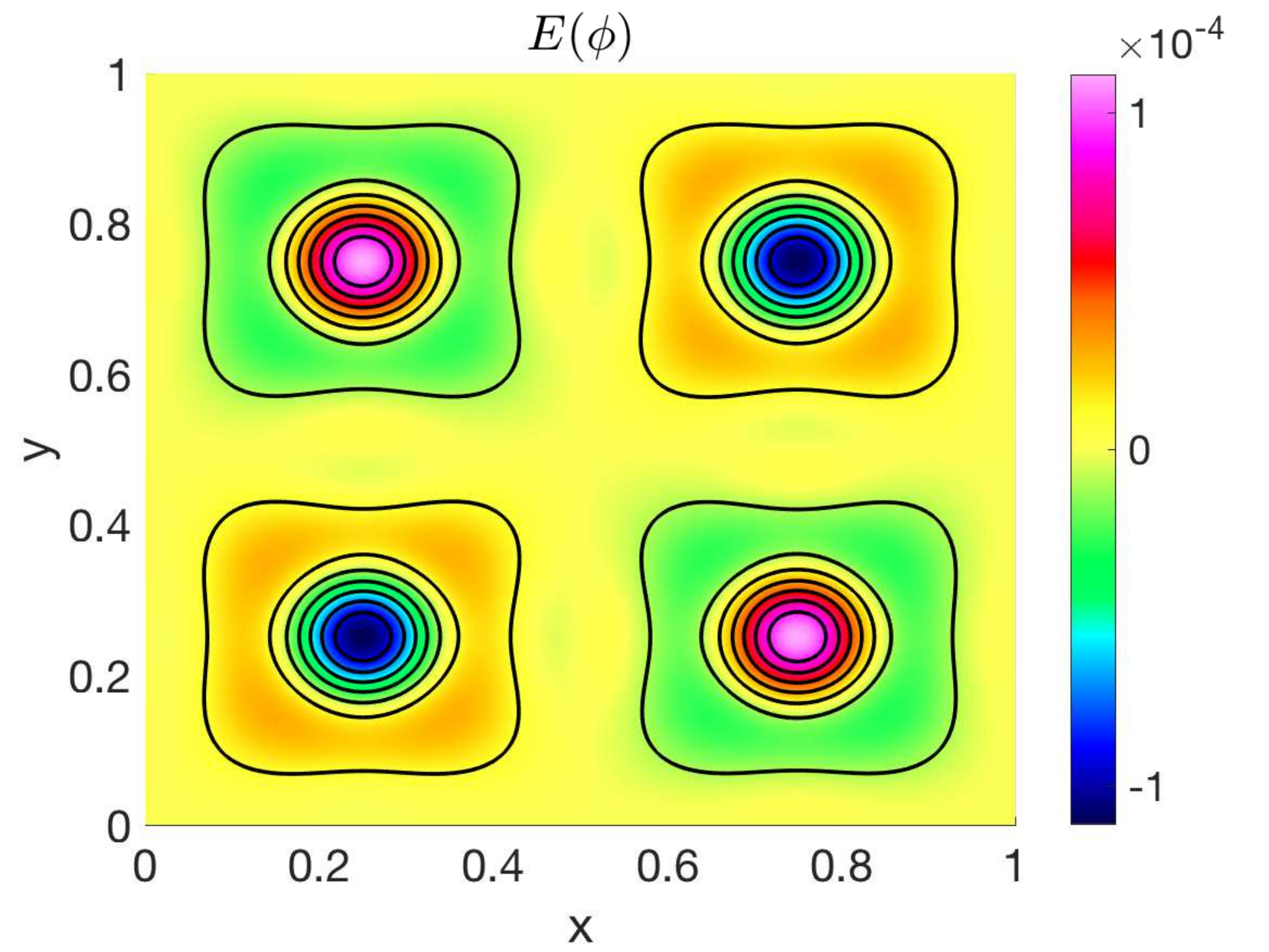}{\figWidth}};

\draw(4,15) node[anchor=north]{\footnotesize(a) Solution};
\draw(13,15) node[anchor=north]{\footnotesize(b) Simply supported BC};
\draw(4,8) node[anchor=north]{\footnotesize(c) Clamped BC};
\draw(13,8) node[anchor=north]{\footnotesize(d) Free BC};
\draw(4,1) node[anchor=north]{\footnotesize(e) Clamped-Supported (CS) BC};
\draw(13,1) node[anchor=north]{\footnotesize(f) Clamped-Free (CF) BC};
%
\end{tikzpicture}
\caption{Contour plots showing the solution and errors of the $\phi$ component  of the linear coupled system with various boundary conditions  on grid $\mathcal{G}_{640}$. The tolerance for this simulation is $tol=10^{-6}$.  Results obtained from  the Picard method are shown here; those of the other two algorithms are similar.   
}\label{fig:linearCoupledTestResultContour_phi}
\end{center}
\end{figure}
}

{
\newcommand{\figWidth}{8cm}
\newcommand{\trimfig}[2]{\trimw{#1}{#2}{0.}{0.}{0.}{0.0}}
\begin{figure}[hp!]
\begin{center}
\begin{tikzpicture}[scale=1]
\useasboundingbox (0.0,0.0) rectangle (17.,21);  
\draw(-0.5,15.0) node[anchor=south west,xshift=0pt,yshift=0pt] {\trimfig{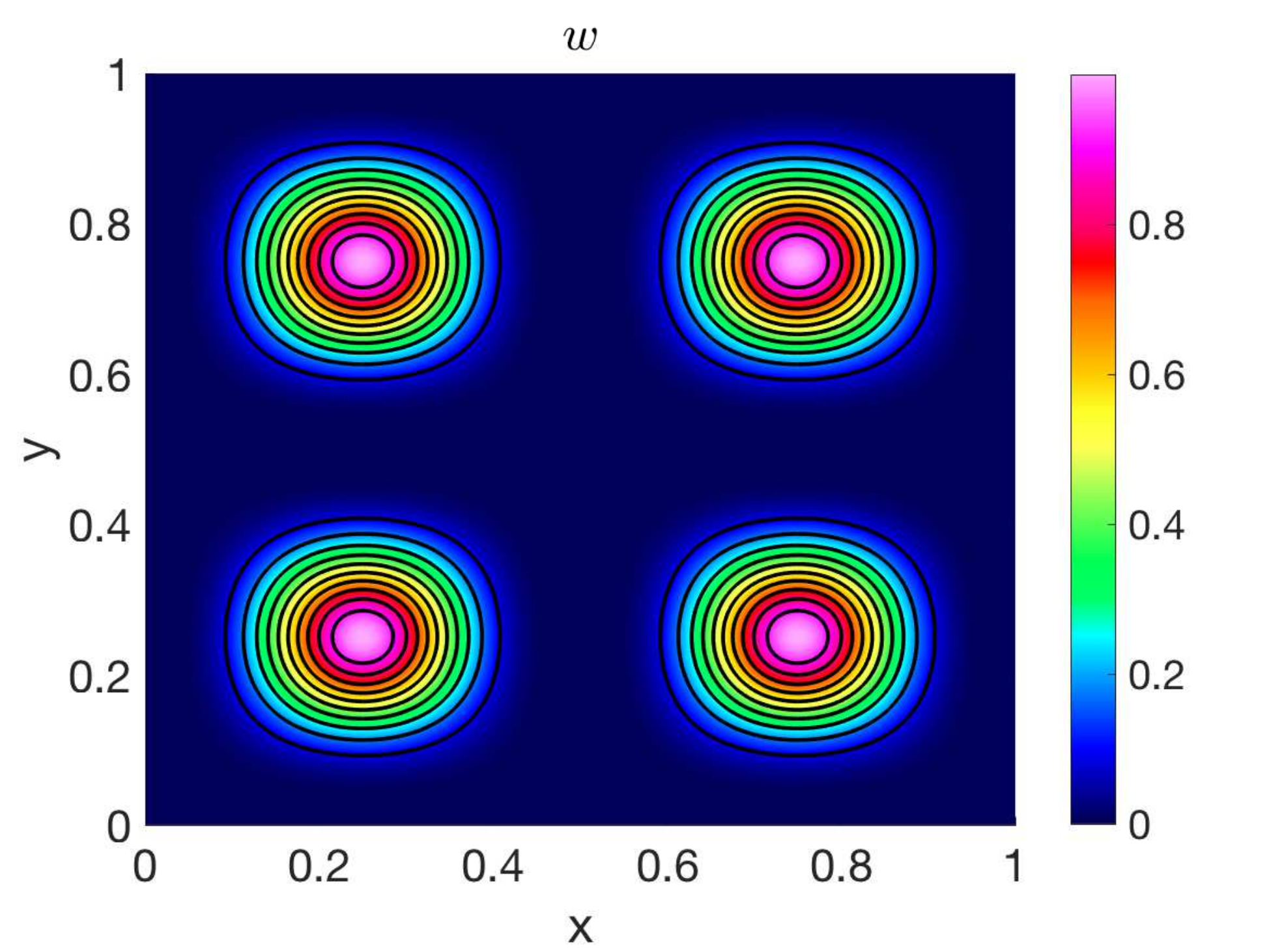}{\figWidth}};
\draw(8.0,15.0) node[anchor=south west,xshift=0pt,yshift=0pt] {\trimfig{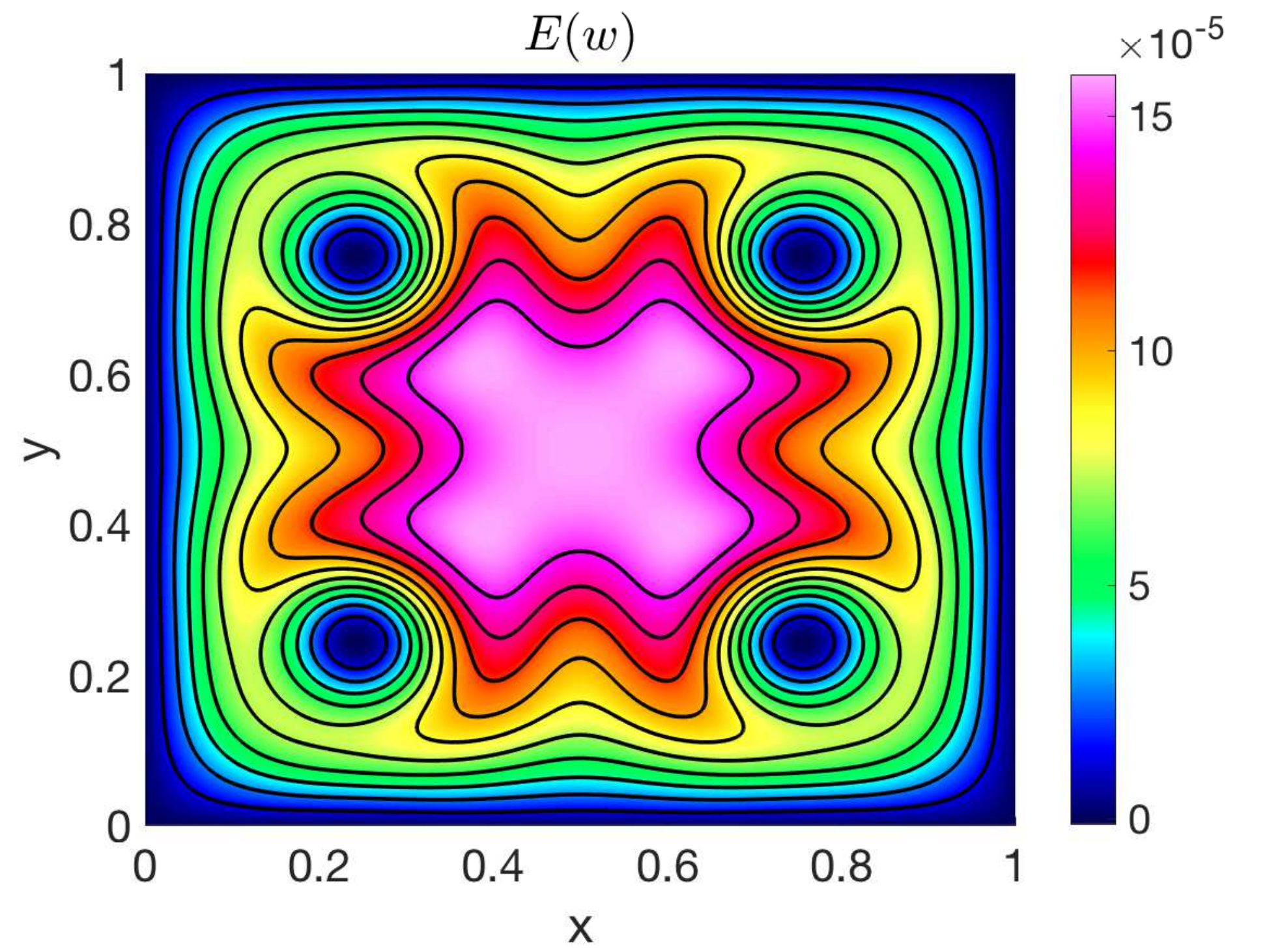}{\figWidth}};
\draw(-0.5,8.0) node[anchor=south west,xshift=0pt,yshift=0pt] {\trimfig{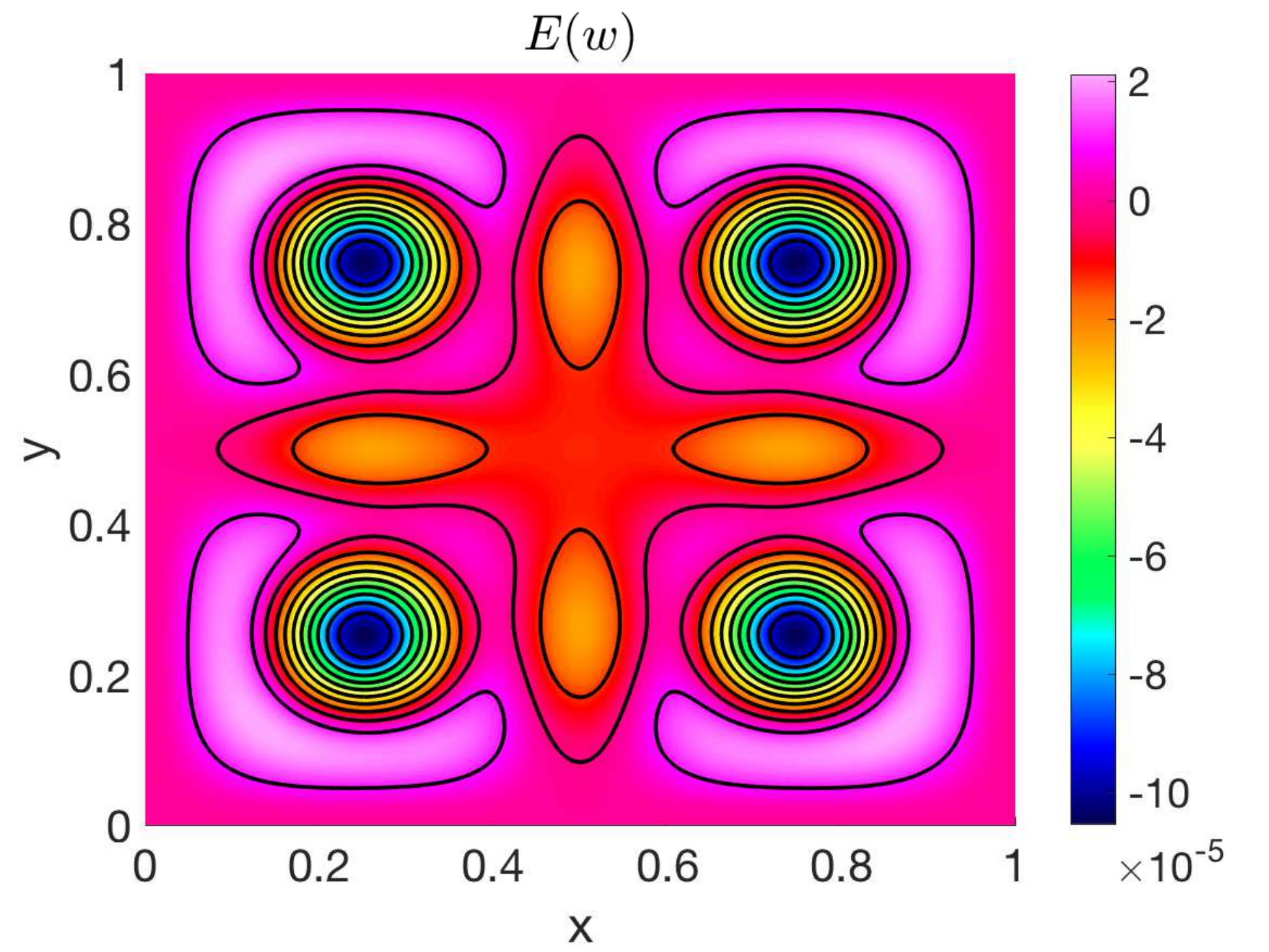}{\figWidth}};
\draw(8.0,8.0) node[anchor=south west,xshift=0pt,yshift=0pt] {\trimfig{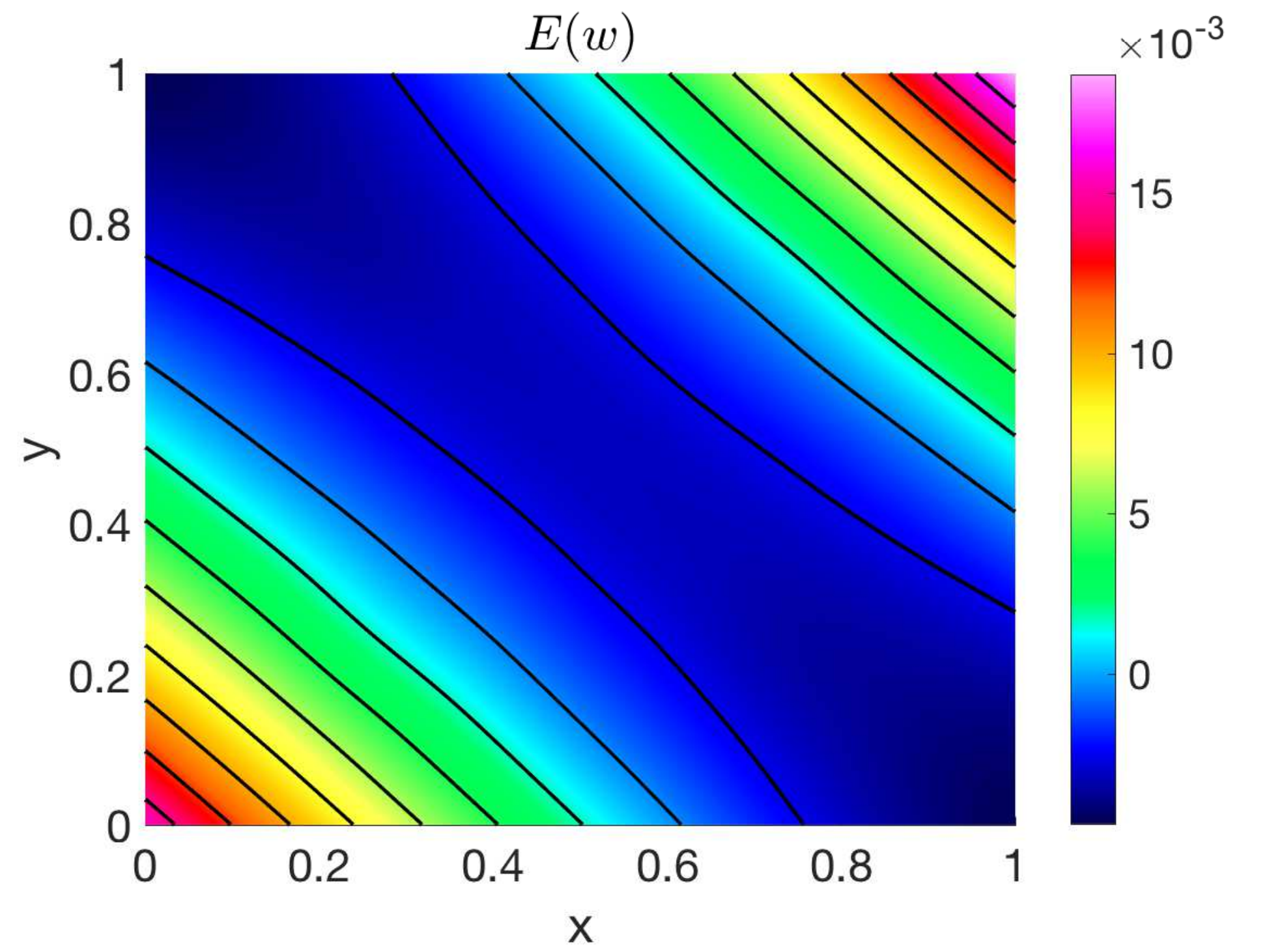}{\figWidth}};
\draw(-0.5,1.0) node[anchor=south west,xshift=0pt,yshift=0pt] {\trimfig{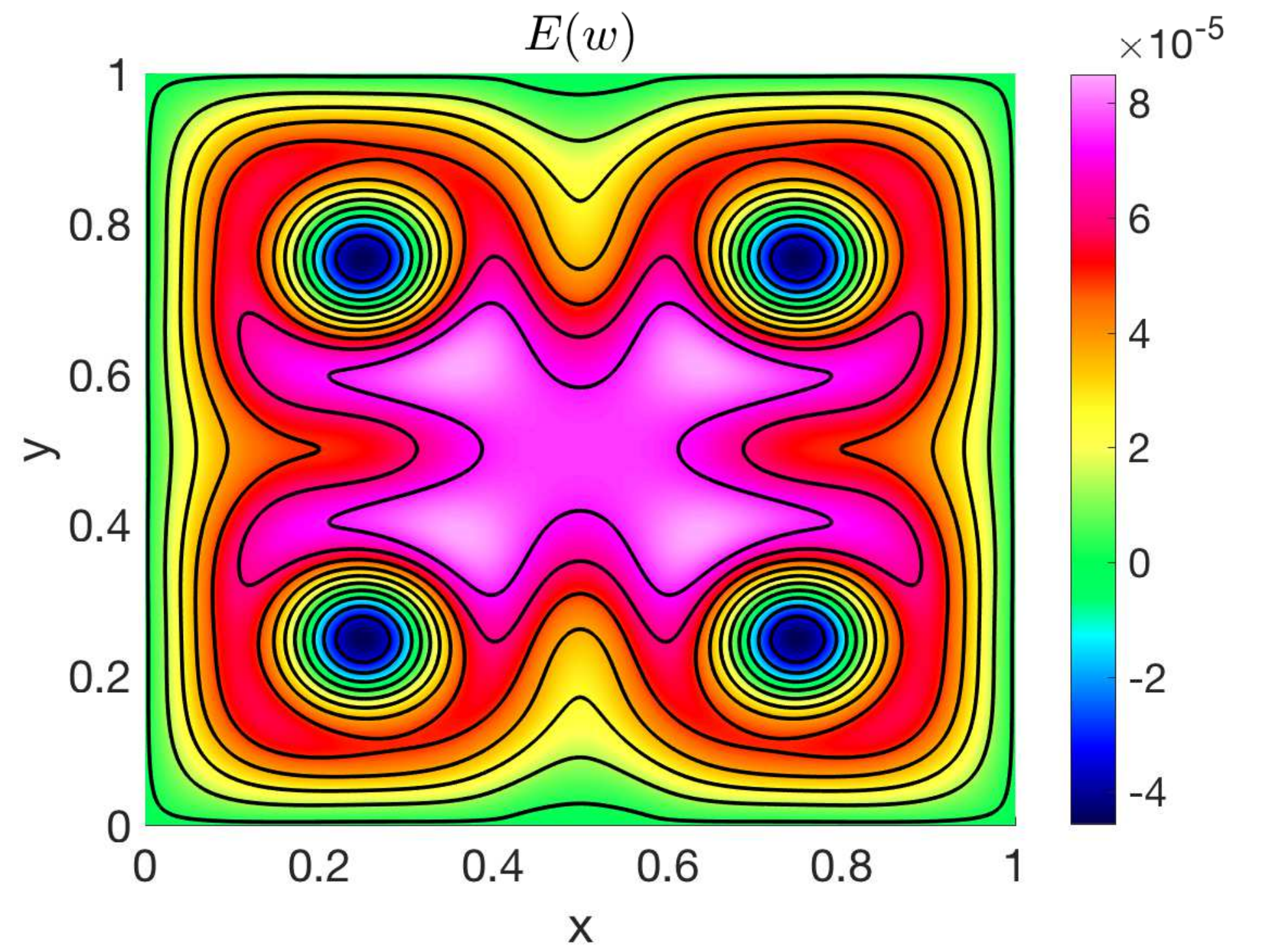}{\figWidth}};
\draw(8.0,1.0) node[anchor=south west,xshift=0pt,yshift=0pt] {\trimfig{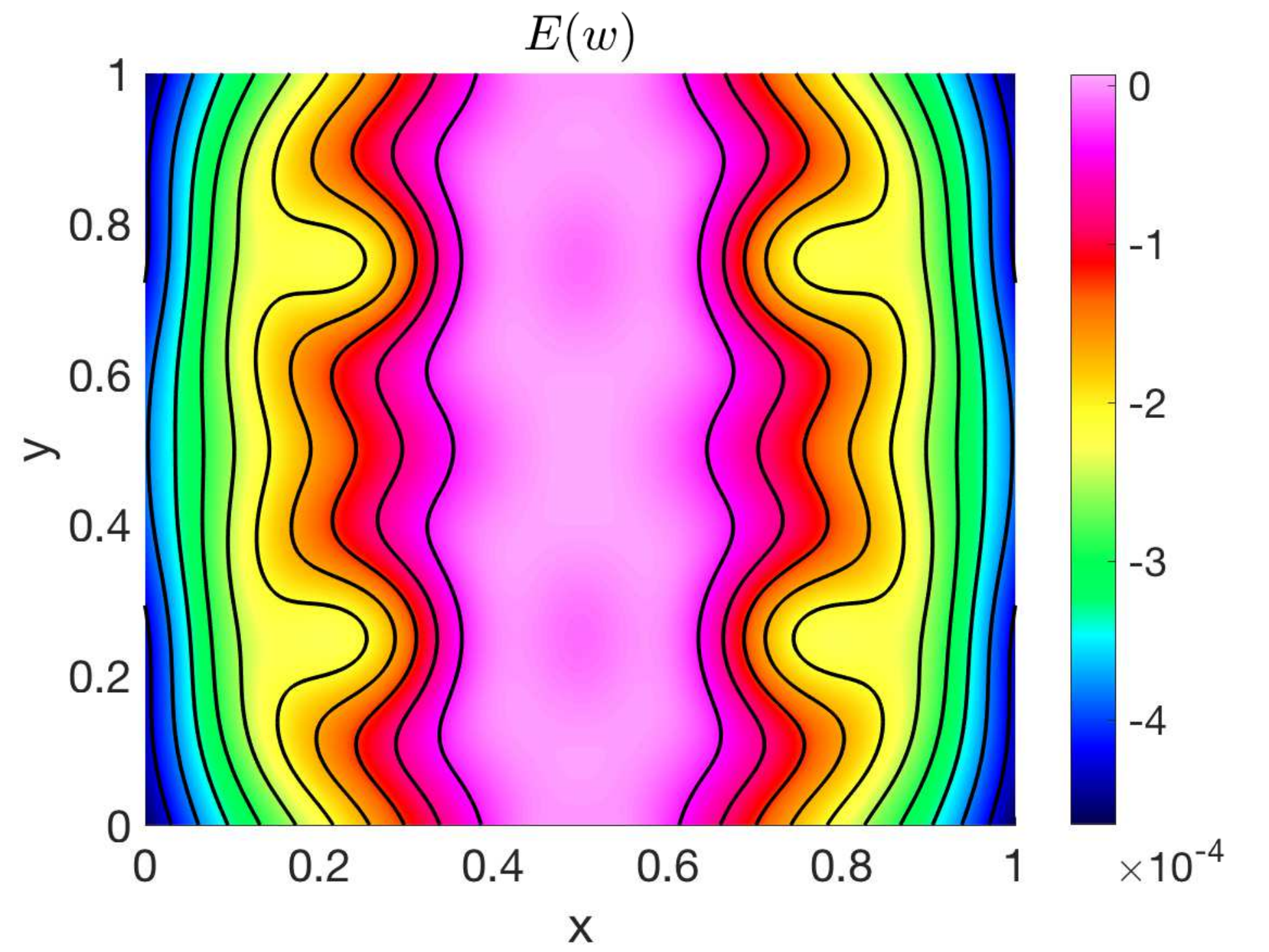}{\figWidth}};

\draw(4,15) node[anchor=north]{\footnotesize(a) Solution};
\draw(13,15) node[anchor=north]{\footnotesize(b) Simply supported BC};
\draw(4,8) node[anchor=north]{\footnotesize(c) Clamped BC};
\draw(13,8) node[anchor=north]{\footnotesize(d) Free BC};
\draw(4,1) node[anchor=north]{\footnotesize(e) Clamped-Supported (CS) BC};
\draw(13,1) node[anchor=north]{\footnotesize(f) Clamped-Free (CF) BC};
%
\end{tikzpicture}
\caption{Contour plots showing the solution and errors of the $w$ component  of the linear coupled system with various boundary conditions  on grid $\mathcal{G}_{640}$. Tolerance for this simulation is $tol=10^{-6}$. Results obtained from  the Picard method are shown here; those of the other two algorithms are similar.  }\label{fig:linearCoupledTestResultContour_w}
\end{center}
\end{figure}
}

It is also observed that the Picard method both explicit and implicit are more efficient than the Newton and trust-region-dogleg methods (fsolve).  At each iteration step,  we need to solve two $N\times N$ matrix equations using the Picard method, while  a $2N\times 2N$ (Jacobian) matrix equation is solved with the other two methods. 
 For large $N$, the Newton's method and the  trust-region-dogleg method   are more computationally expensive at each step, which may result in a longer overall time to solve the system even though the Newton method converges in a faster rate than the Picard method. In addition, computer memory can also be an issue using Newton solve and fsolve on a high resolution grid.  For example, both Newton solve and fsolve encountered out-of-memory  error when solving this problem on grid $\G_{640}$ (the finest grid considered for the mesh refinement study)  using a single processor of a linux desktop computer equipped with 64GB memory; however, the Picard method handles this resolution with no problem. For this reason,  the data points for $G_{640}$ are absent in plots (c) and (d) of Fig.~\ref{fig:convRateLinearCoupledSystem}.

{
\newcommand{\figWidth}{8cm}
\newcommand{\trimfig}[2]{\trimw{#1}{#2}{0.}{0.}{0.}{0.0}}
\begin{figure}[h!]
\begin{center}
\begin{tikzpicture}[scale=1]
\useasboundingbox (0.0,0.0) rectangle (17.,14);  
\draw(-0.5,7.5) node[anchor=south west,xshift=0pt,yshift=0pt] {\trimfig{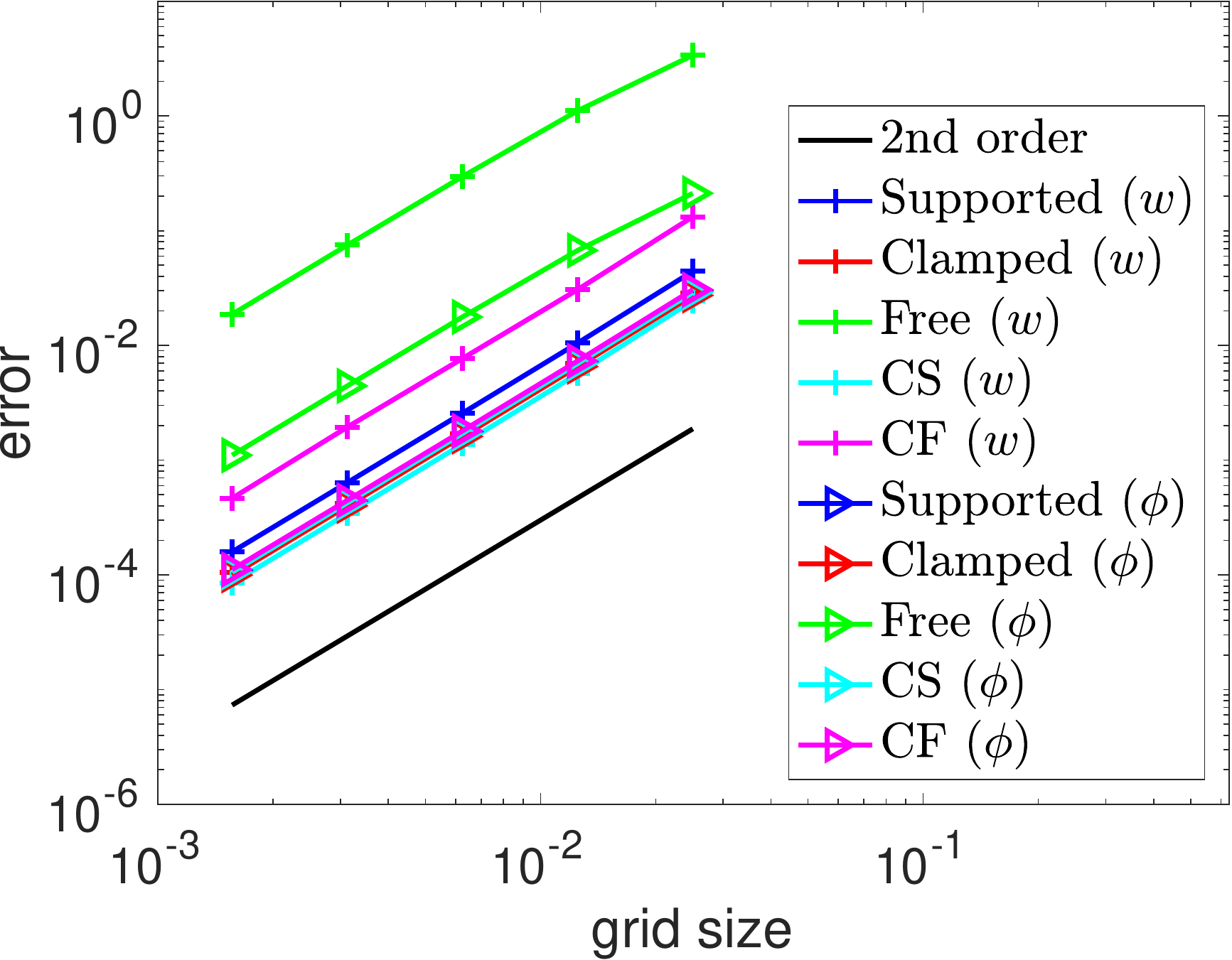}{\figWidth}};
\draw(8.0,7.5) node[anchor=south west,xshift=0pt,yshift=0pt] {\trimfig{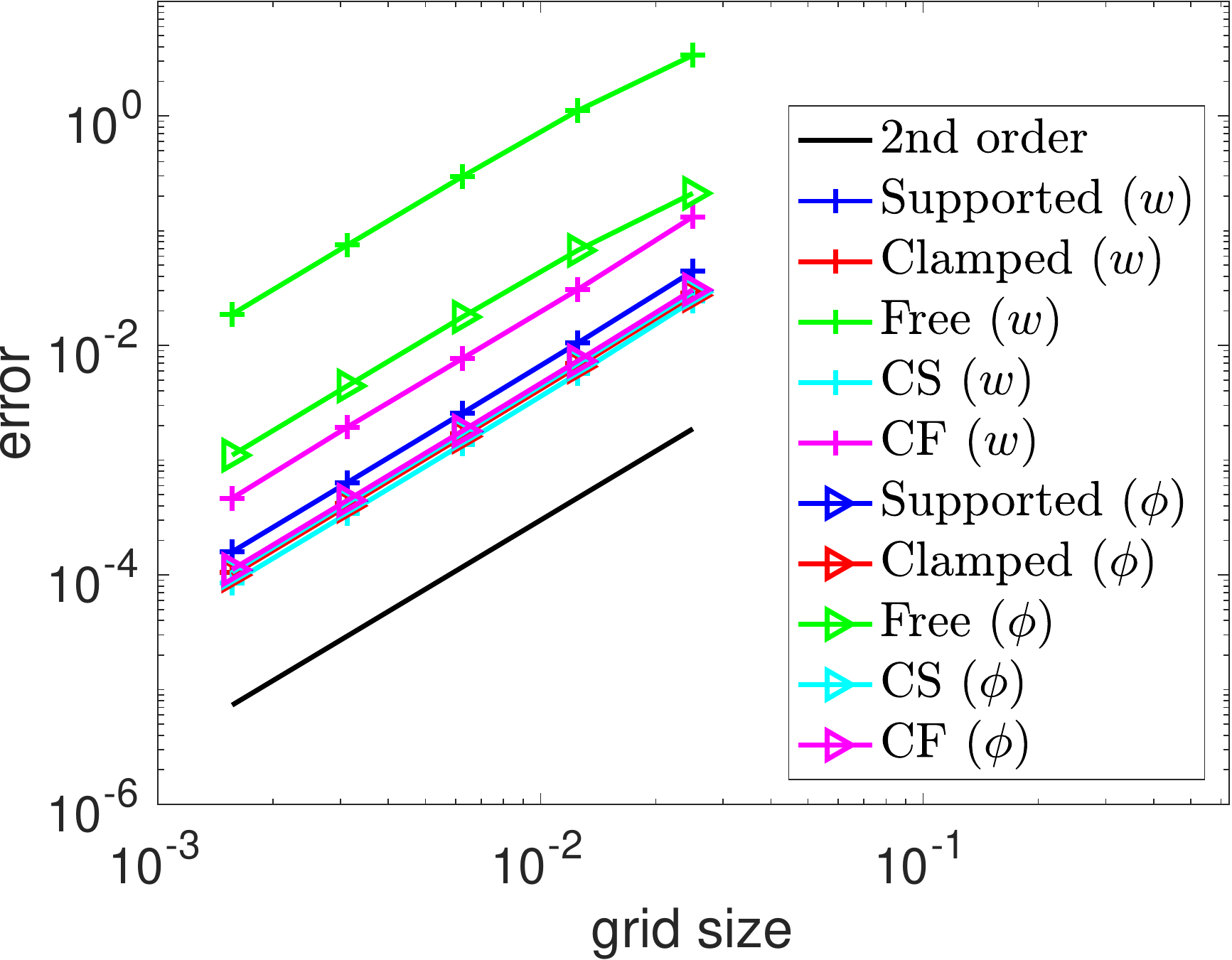}{\figWidth}};
\draw(-0.5,0.5) node[anchor=south west,xshift=0pt,yshift=0pt] {\trimfig{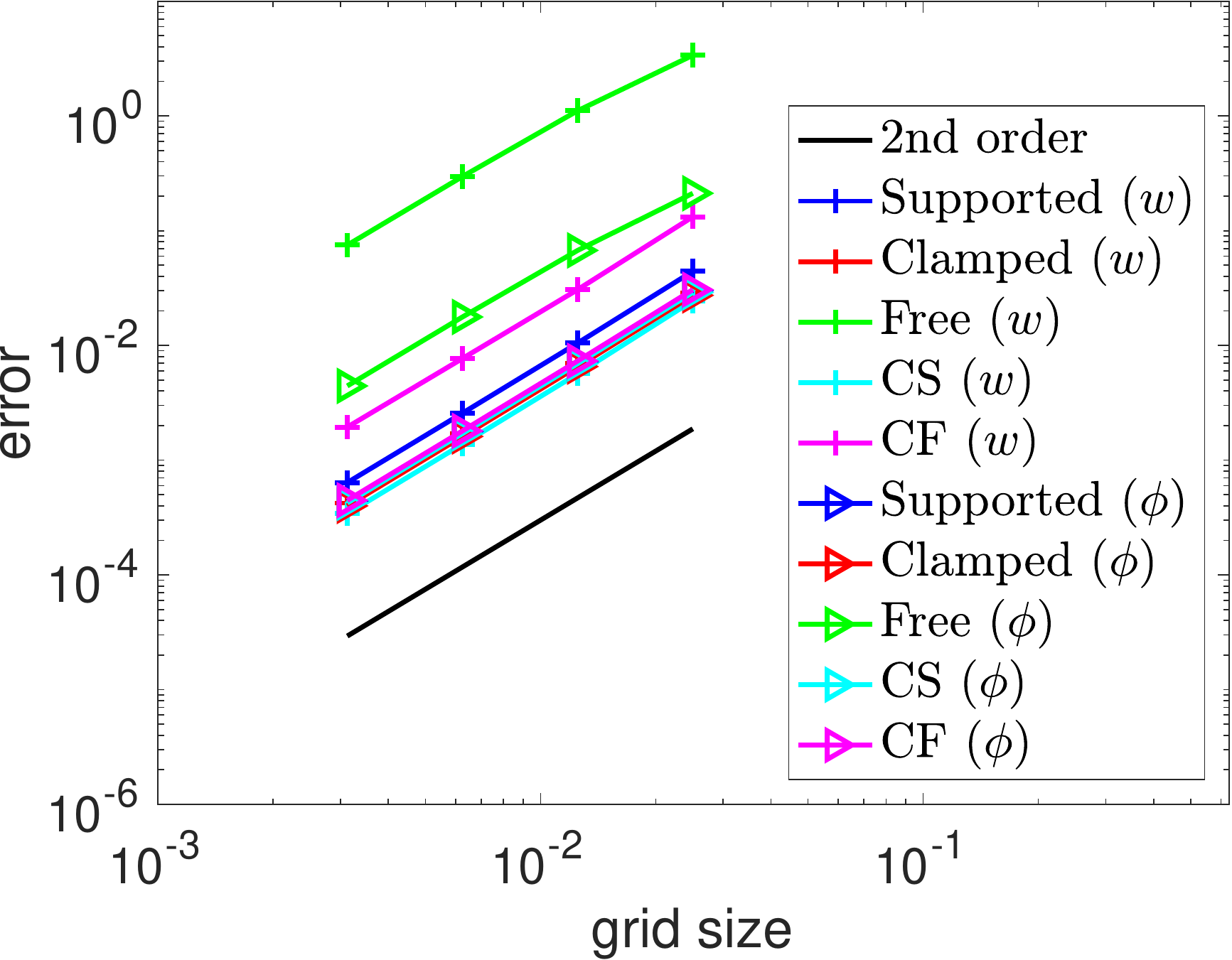}{\figWidth}};
\draw(8.0,0.5) node[anchor=south west,xshift=0pt,yshift=0pt] {\trimfig{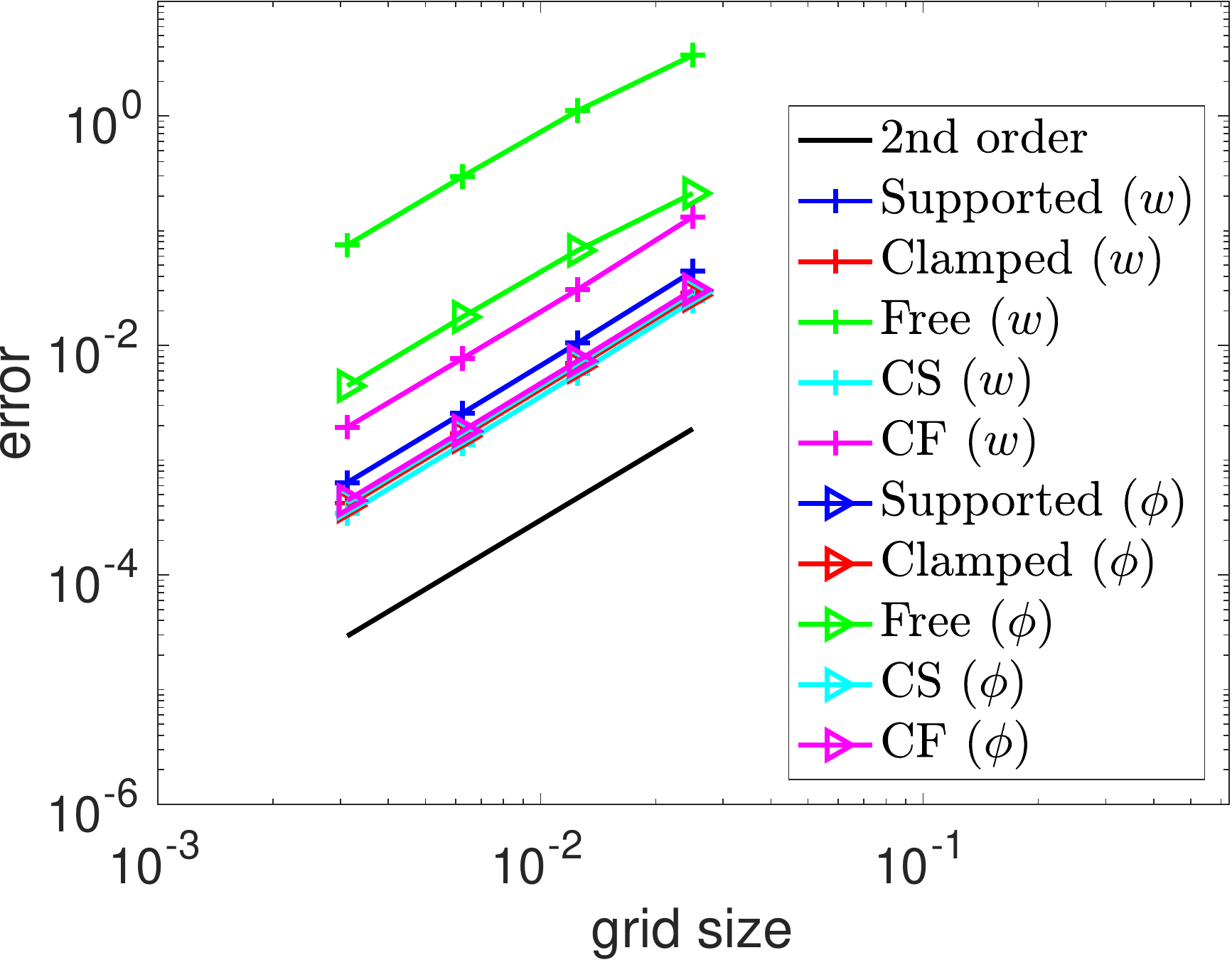}{\figWidth}};

\draw(4,7.5) node[anchor=north]{\footnotesize(a) Picard method (explicit; i.e., $\delta=0$)};
\draw(13,7.5) node[anchor=north]{\footnotesize(b) Picard method (implicit; i.e., $\delta=1$)};
\draw(4,0.5) node[anchor=north]{\footnotesize(c) Newton's method};
\draw(13,0.5) node[anchor=north]{\footnotesize(d) fsolve};
%
\end{tikzpicture}
\caption{A mesh refinement study for the numerical solutions of the linear coupled system with various methods and boundary conditions. Errors are in  maximum norm $L_{\infty}$. Tolerance for this simulation is $tol=10^{-6}$.}\label{fig:convRateLinearCoupledSystem}
\end{center}
\end{figure}
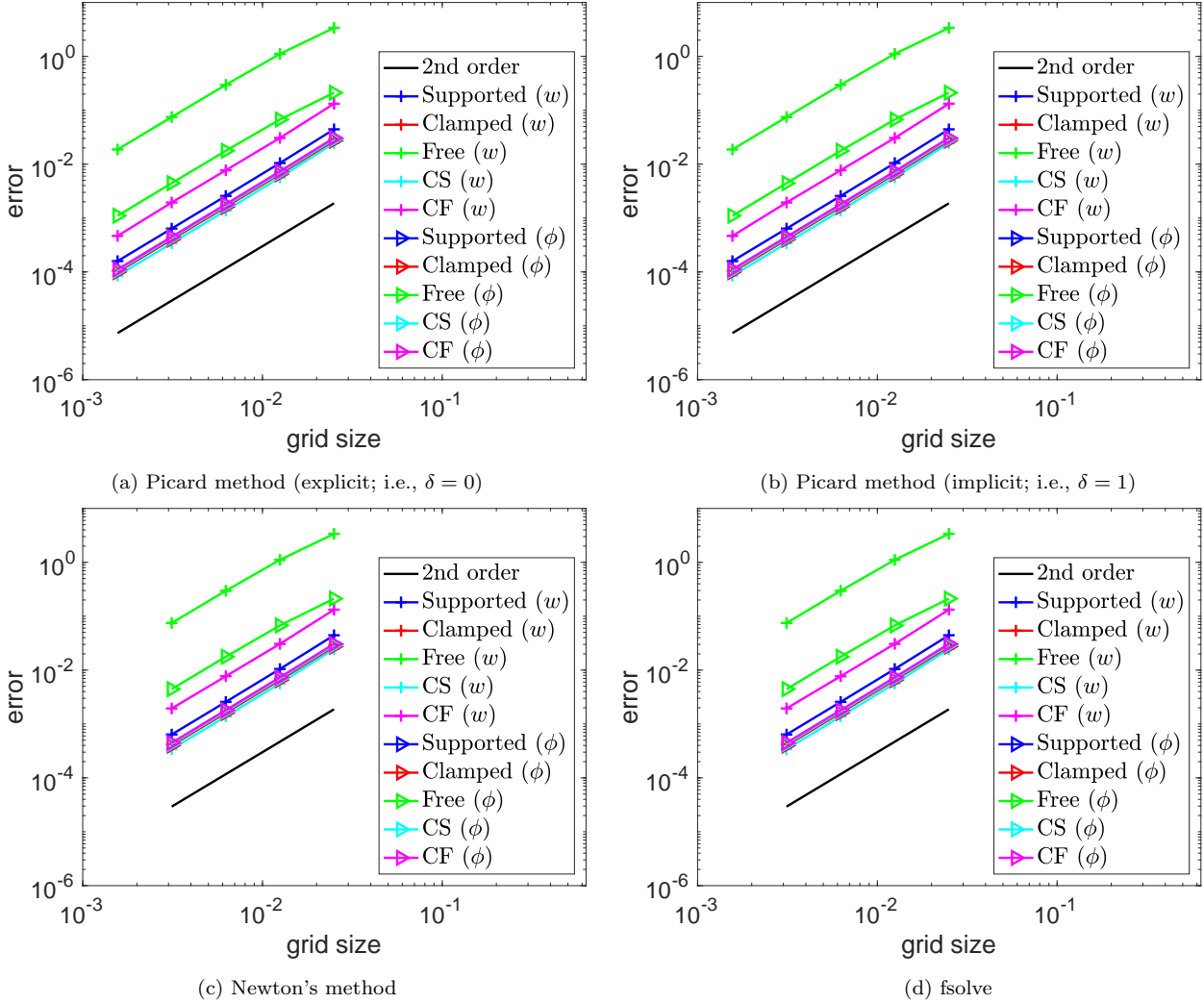
}

\subsubsection{Nonlinear coupled system}
As a final mesh refinement study, we test our numerical methods by solving  the nonlinear shallow shell equations \eqref{eq:coupledSystemNonlinear}.  The exact solutions  and  the  precast shell shape are specified to be the same as the linear coupled system test, which are given in equations \eqref{eq:exactSolutionCoupled} and \eqref{eq:precastShapeCoupled}.   The forcing terms are different due to the nonlinear terms, and they are visualized in Fig.~\ref{fig:nonlinearCoupledSystemForcingsContour}.

{
\newcommand{\figWidth}{8cm}
\newcommand{\trimfig}[2]{\trimw{#1}{#2}{0.}{0.}{0.}{0.0}}
\begin{figure}[h!]
\begin{center}
\begin{tikzpicture}[scale=1]
\useasboundingbox (0.0,0.0) rectangle (17.,7);  
\draw(-0.5,0.0) node[anchor=south west,xshift=0pt,yshift=0pt] {\trimfig{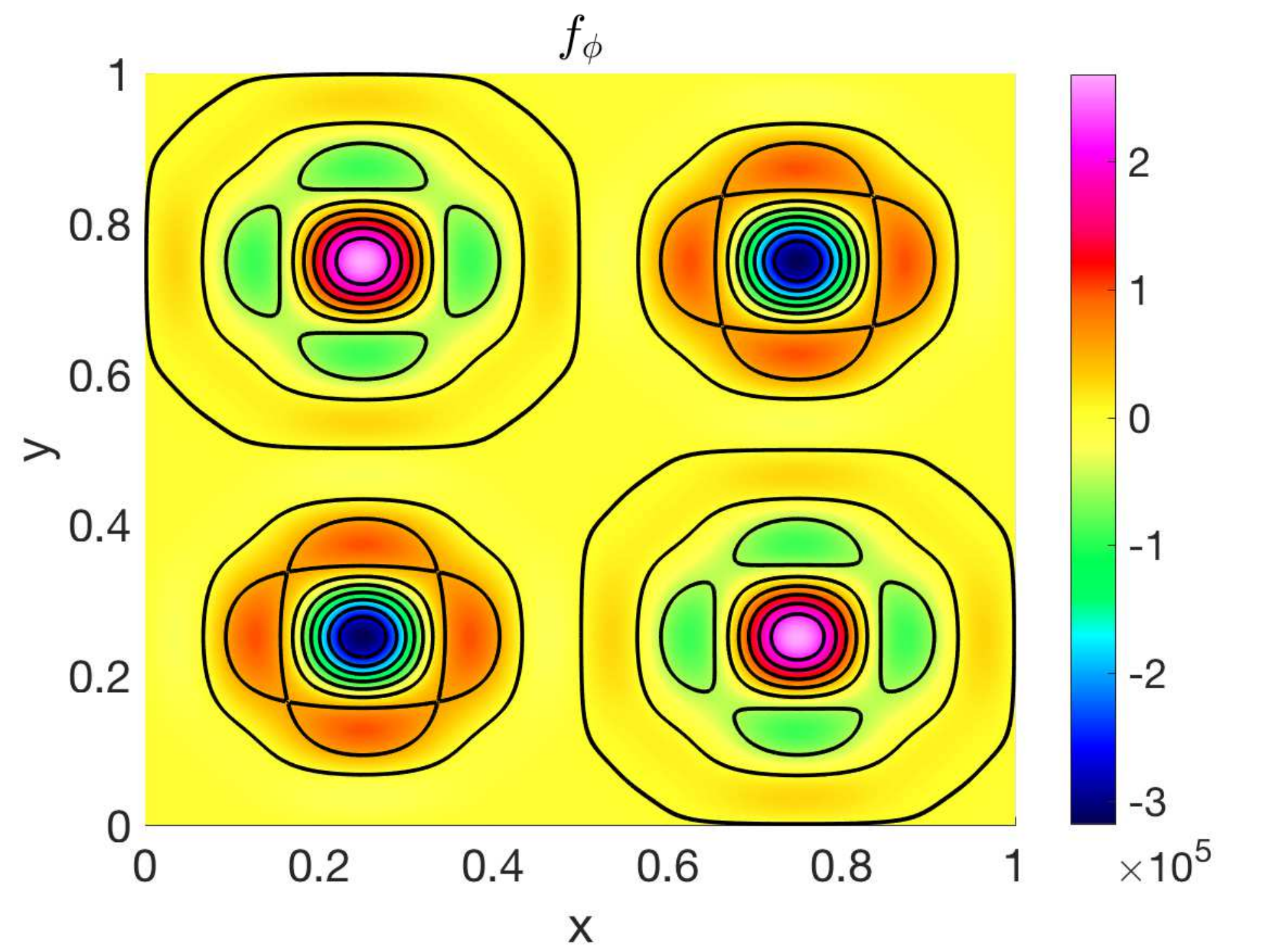}{\figWidth}};
\draw(8.0,0.0) node[anchor=south west,xshift=0pt,yshift=0pt] {\trimfig{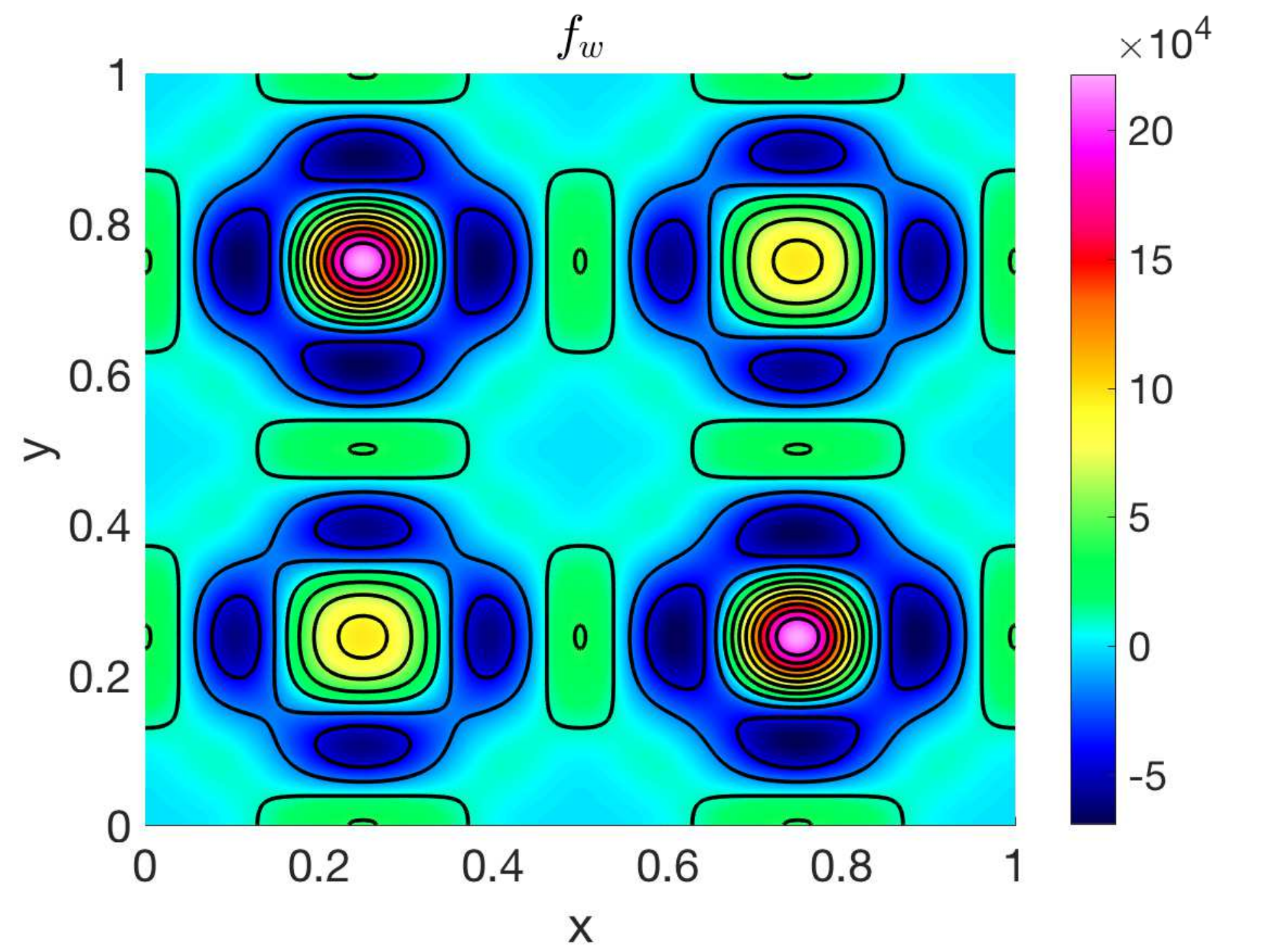}{\figWidth}};
%
\end{tikzpicture}
\caption{Contour plots showing the forcing terms in the nonlinear coupled system \eqref{eq:coupledSystemNonlinear}.}\label{fig:nonlinearCoupledSystemForcingsContour}
\end{center}
\end{figure}
}

The numerical solutions for the nonlinear case are   accurate  for all the boundary conditions and for all numerical methods that are considered in this paper.   As are illustrated in  Fig.~\ref{fig:nonlinearCoupledTestResultContour_phi} and Fig.~\ref{fig:nonlinearCoupledTestResultContour_w}, the numerical errors are well behaved in the sense that   their magnitudes are small and smooth throughout the domain including the boundaries. Importantly, the technique employed to  regularize the displacement equation with free boundary conditions performs well in the context of  a nonlinear coupled system, too. The 2nd-order spatial accuracy for all of the numerical schemes are again confirmed by the mesh refinement results shown in Fig.~\ref{fig:convRateNonlinearCoupledSystem}. We note that  the Picard method with free boundary conditions demands the grid to be fine enough to converge; thus,  the mesh refinement study  starts from grid $\G_{40}$ for $\delta=0$ and from $\G_{80}$ for $\delta=1$.

\begin{table}[h!]
\begin{center}
\begin{tabular}{|c|c|c|c||c|c|c|} \hline
\multicolumn{7}{|c|}{Run-time performance} \\ \hline\hline  
               & \multicolumn{3}{|c||}{Picard ($\delta=0$)}            & \multicolumn{3}{|c|}{Newton}     \\ \cline{2-7}
               & s/step            & steps  &  rate  &   s/step        & steps  &  rate    \\ \cline{2-7}
  $\G_{160}$      & 2.48 & 24 &1.00 & 5.38 &5  &1.83                            \\ \cline{2-7}
 $\G_{320}$     &  15.87   &  19  &    1.01     &    80.23  & 5  & 1.64   \\ \cline{2-7}
 $\G_{640}$   &   255.30  &  17  &    1.02      &    \multicolumn{3}{|c|}{out-of-memory}       \\ \hline 
\end{tabular}
\end{center}
\caption{
Comparison of the run-time performance of the explicit Picard  method versus the Newton's method for the nonlinear shallow shell equations with free boundary conditions. The column labeled ``s/step'' gives the CPU time in seconds per iteration step; the column labeled ``steps'' gives the number of steps taken; and the column labeled ``rate'' gives the estimated rate of convergence of the corresponding method.
}\label{tab:performance}
\end{table}

{
\newcommand{\figWidth}{8cm}
\newcommand{\trimfig}[2]{\trimw{#1}{#2}{0.}{0.}{0.}{0.0}}
\begin{figure}[hp!]
\begin{center}
\begin{tikzpicture}[scale=1]
\useasboundingbox (0.0,0.0) rectangle (17.,21);  
\draw(-0.5,15.0) node[anchor=south west,xshift=0pt,yshift=0pt] {\trimfig{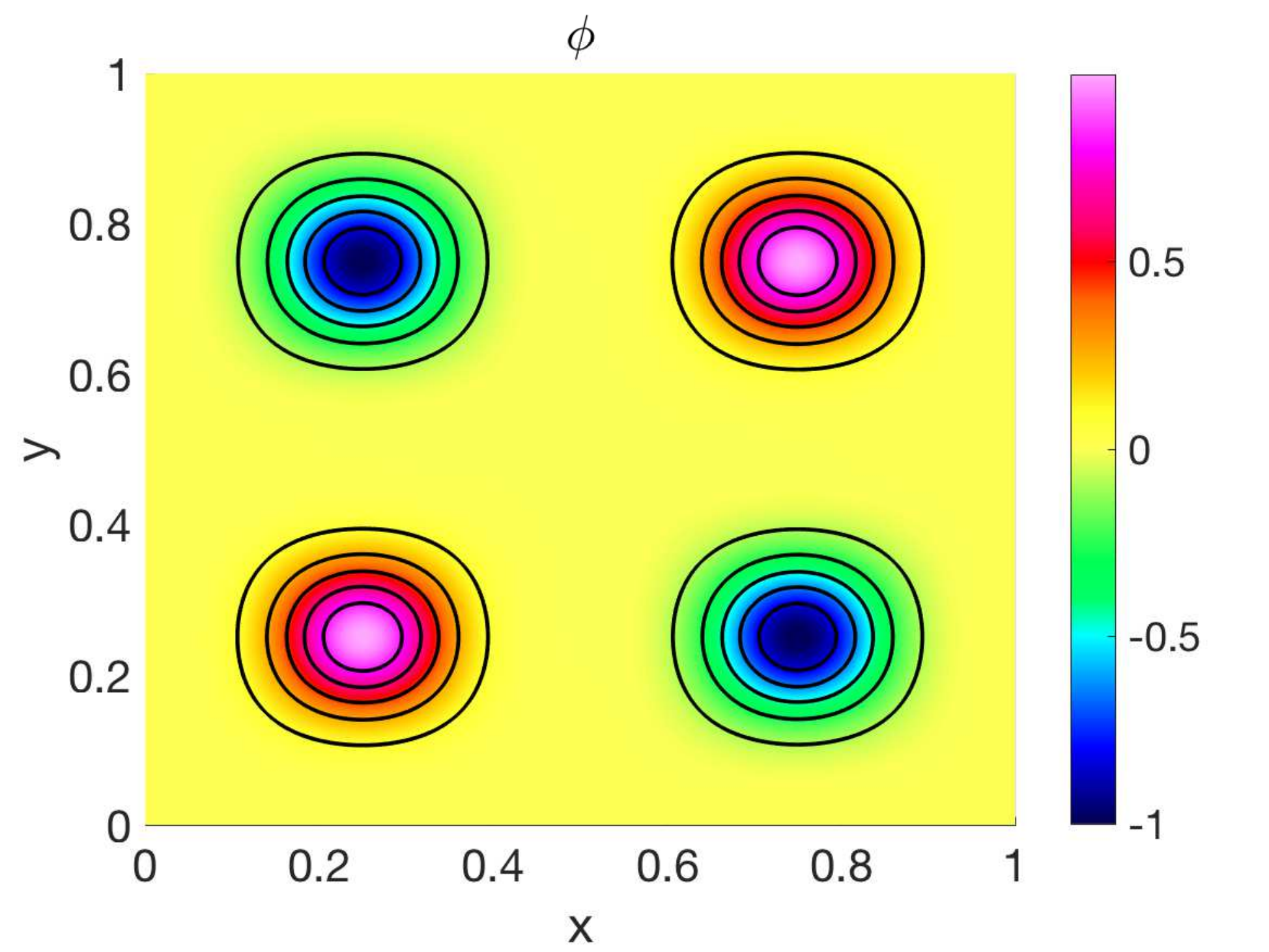}{\figWidth}};
\draw(8.0,15.0) node[anchor=south west,xshift=0pt,yshift=0pt] {\trimfig{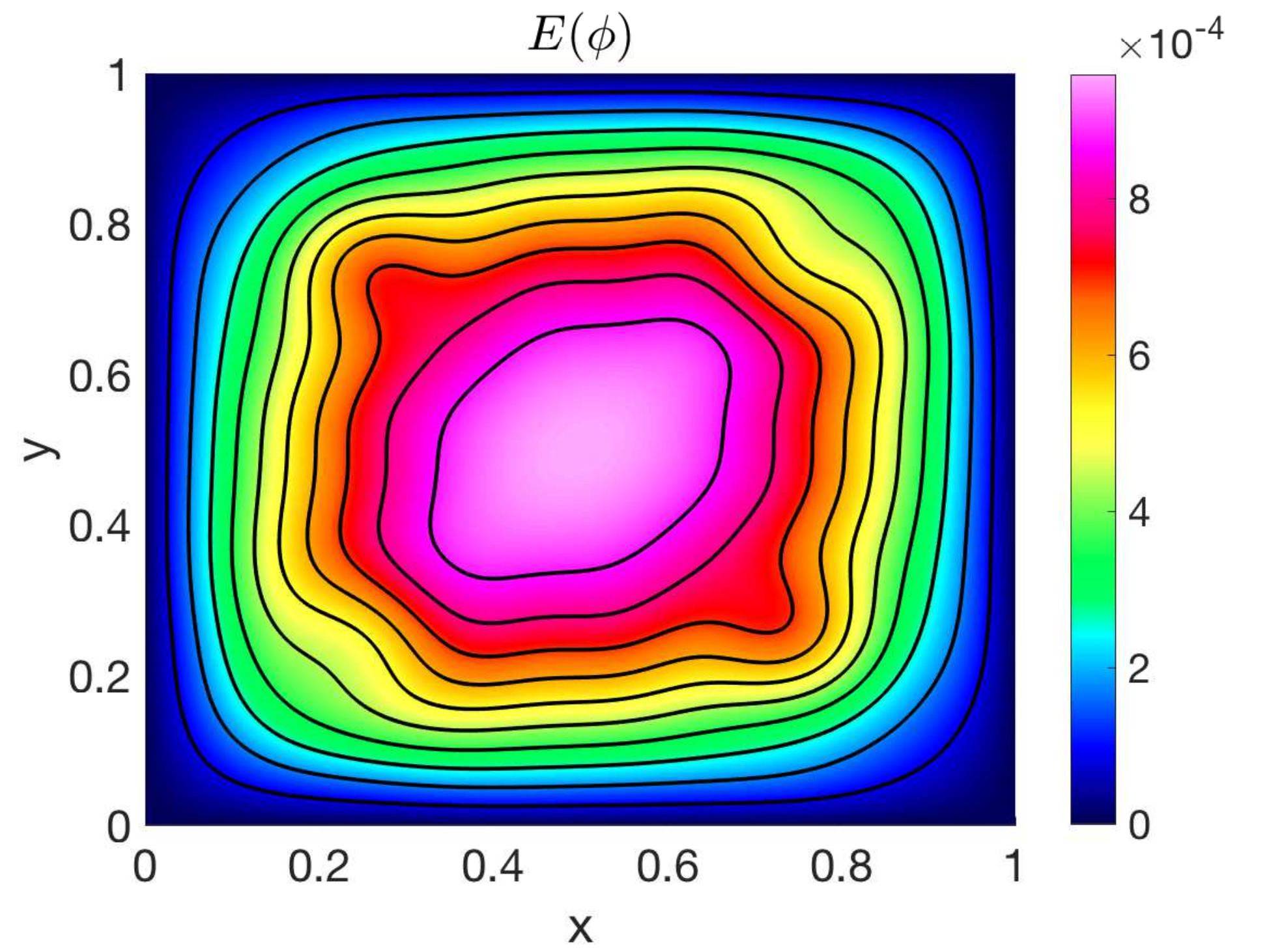}{\figWidth}};
\draw(-0.5,8.0) node[anchor=south west,xshift=0pt,yshift=0pt] {\trimfig{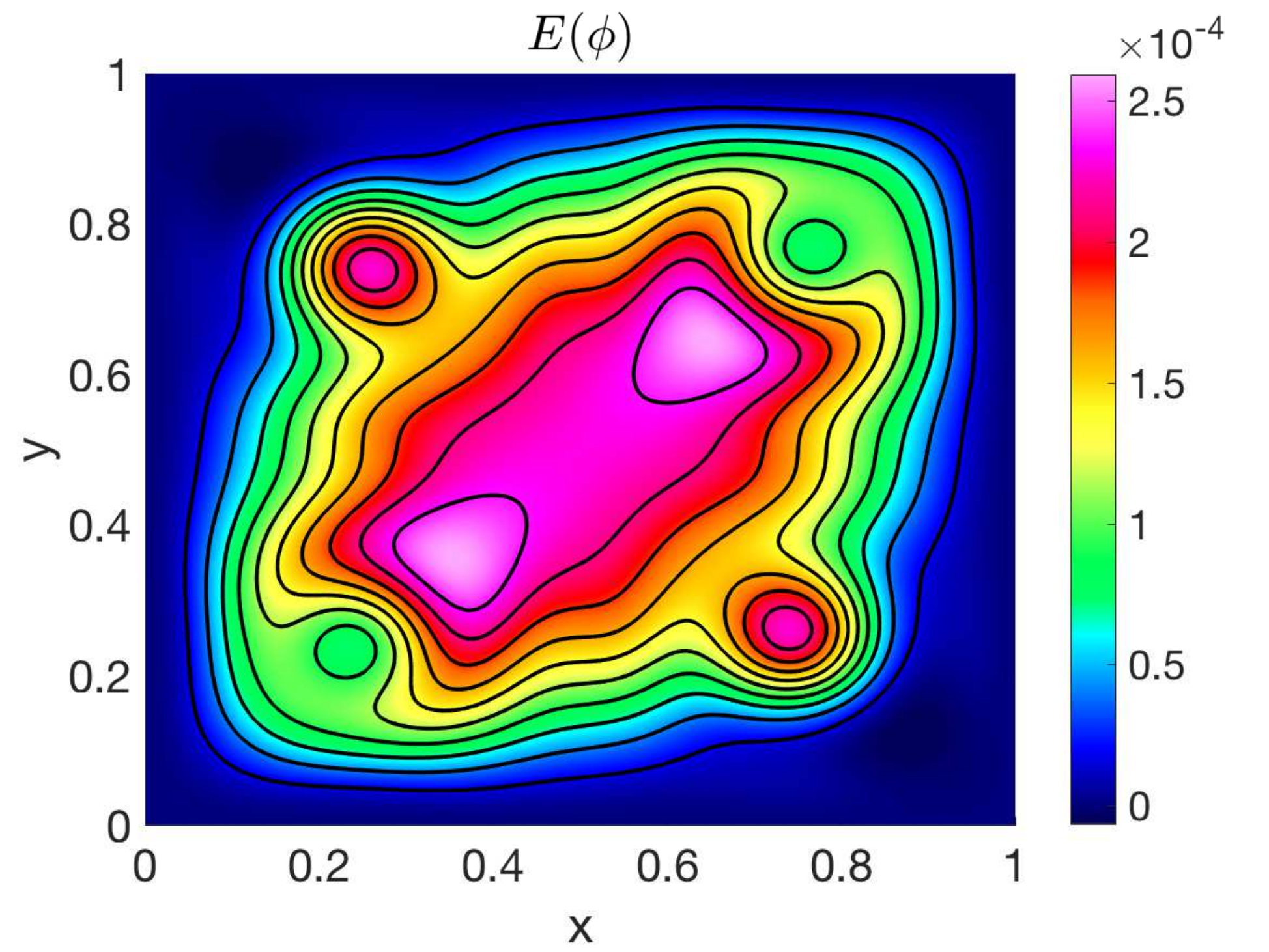}{\figWidth}};
\draw(8.0,8.0) node[anchor=south west,xshift=0pt,yshift=0pt] {\trimfig{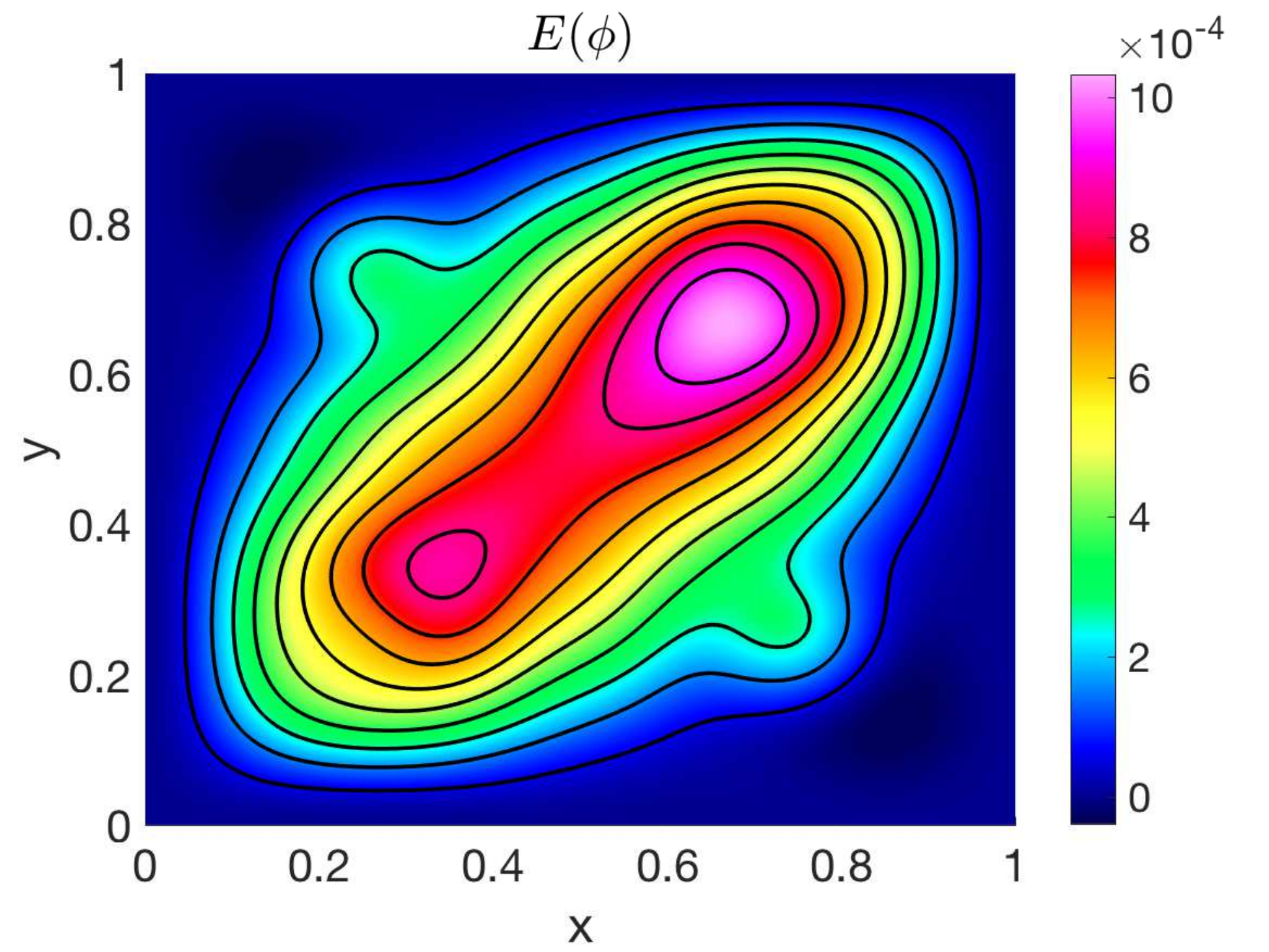}{\figWidth}};
\draw(-0.5,1.0) node[anchor=south west,xshift=0pt,yshift=0pt] {\trimfig{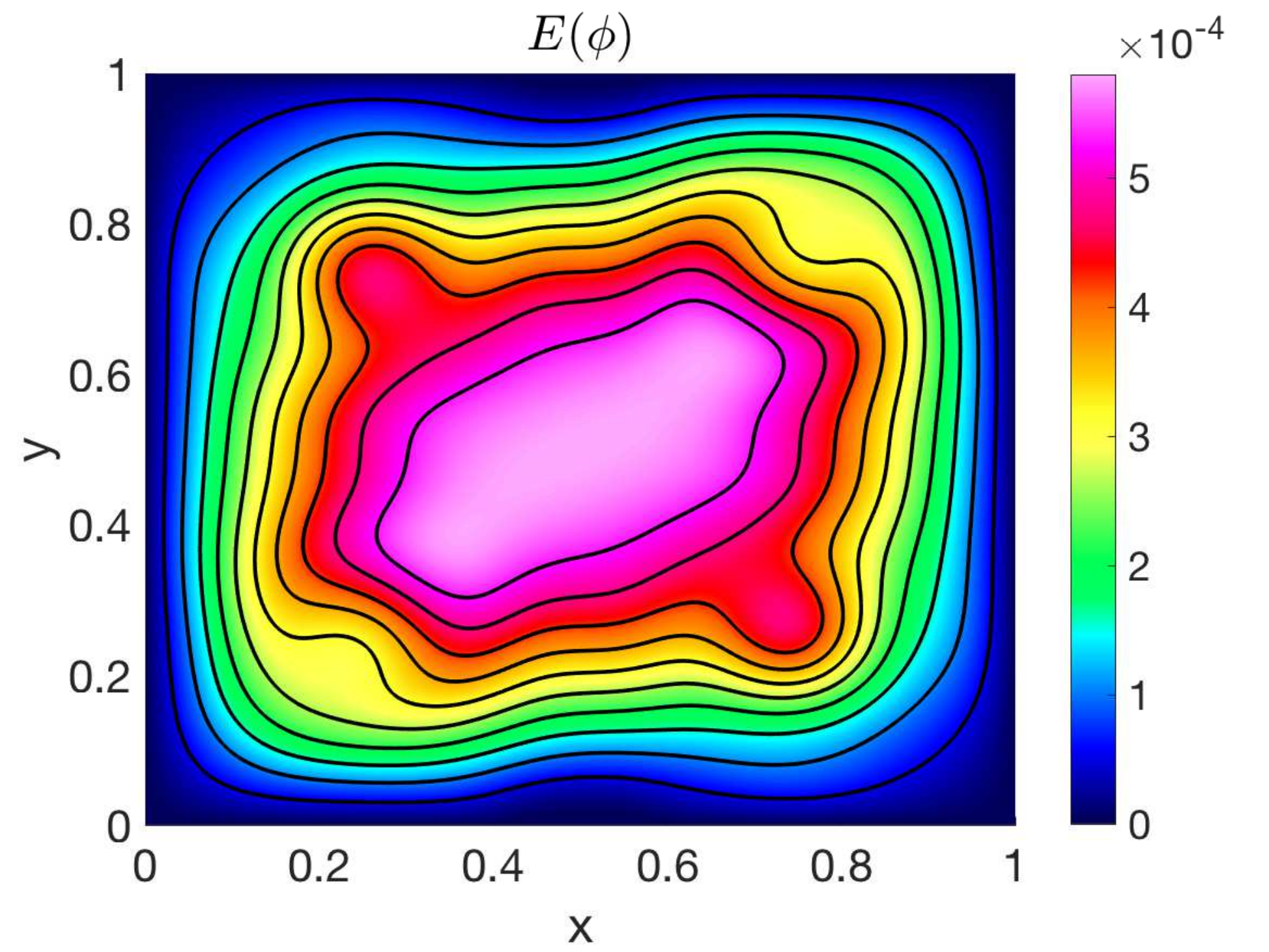}{\figWidth}};
\draw(8.0,1.0) node[anchor=south west,xshift=0pt,yshift=0pt] {\trimfig{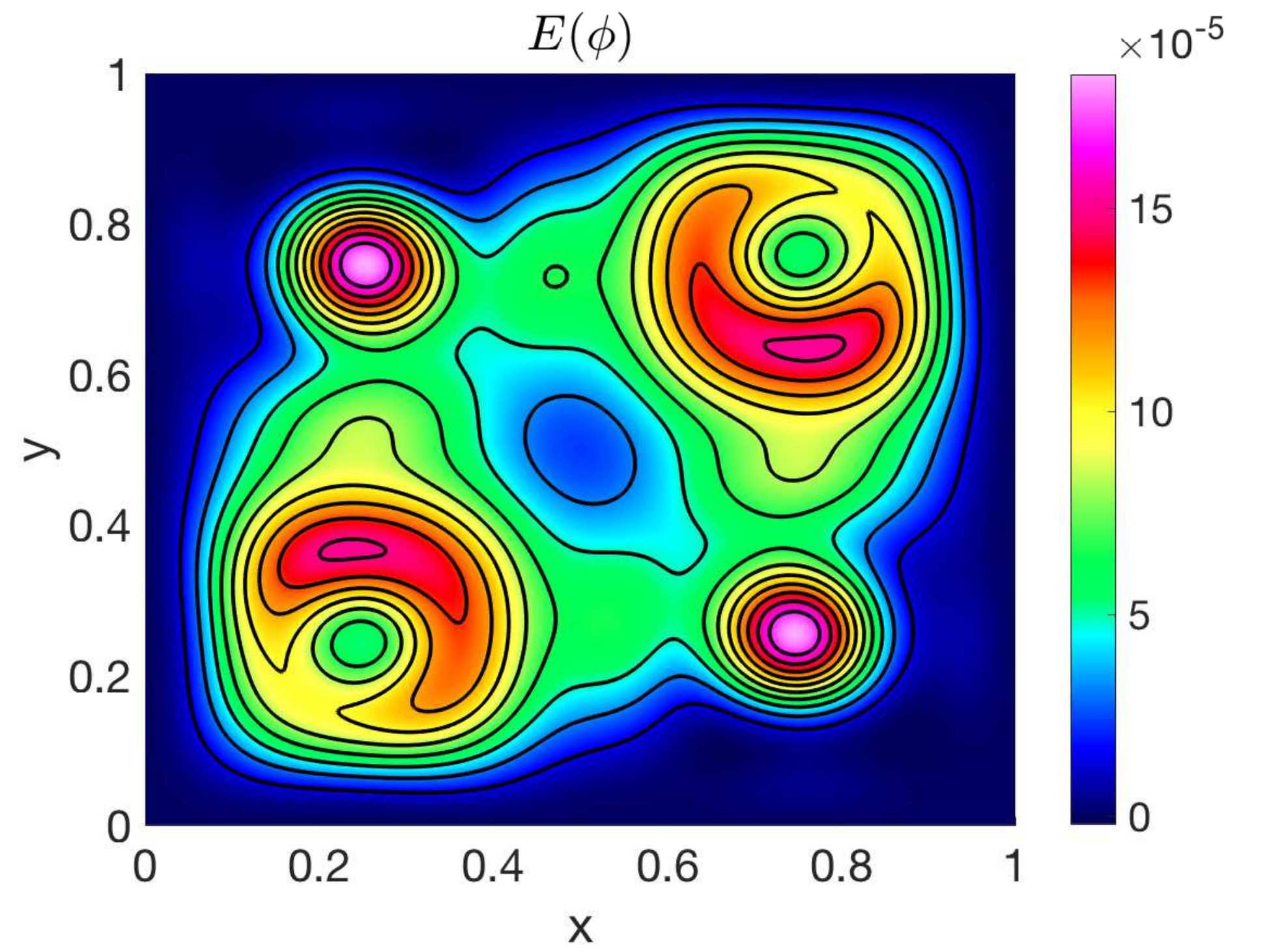}{\figWidth}};

\draw(4,15) node[anchor=north]{\footnotesize(a) Solution};
\draw(13,15) node[anchor=north]{\footnotesize(b) Simply supported BC};
\draw(4,8) node[anchor=north]{\footnotesize(c) Clamped BC};
\draw(13,8) node[anchor=north]{\footnotesize(d) Free BC};
\draw(4,1) node[anchor=north]{\footnotesize(e) Clamped-Supported (CS) BC};
\draw(13,1) node[anchor=north]{\footnotesize(f) Clamped-Free (CF) BC};
%
\end{tikzpicture}
\caption{Contour plots showing the $\phi$ component of the solution and errors of the nonlinear coupled system with various boundary conditions  on grid $\mathcal{G}_{640}$. Tolerance for this simulation is $tol=10^{-6}$. Picard method is used here. Results of the other two algorithms are similar.  
}\label{fig:nonlinearCoupledTestResultContour_phi}
\end{center}
\end{figure}
}

{
\newcommand{\figWidth}{8cm}
\newcommand{\trimfig}[2]{\trimw{#1}{#2}{0.}{0.}{0.}{0.0}}
\begin{figure}[hp!]
\begin{center}
\begin{tikzpicture}[scale=1]
\useasboundingbox (0.0,0.0) rectangle (17.,21);  
\draw(-0.5,15.0) node[anchor=south west,xshift=0pt,yshift=0pt] {\trimfig{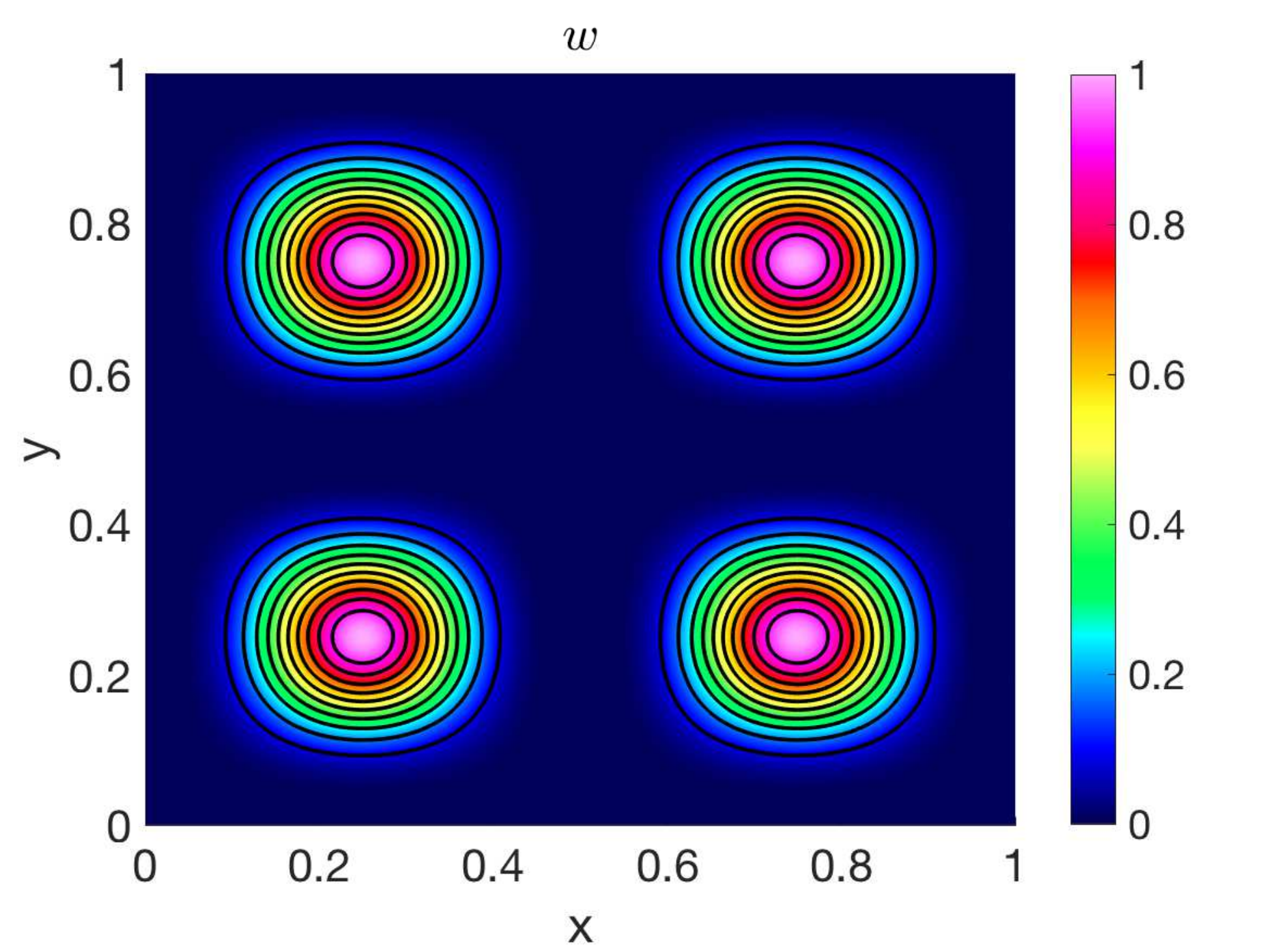}{\figWidth}};
\draw(8.0,15.0) node[anchor=south west,xshift=0pt,yshift=0pt] {\trimfig{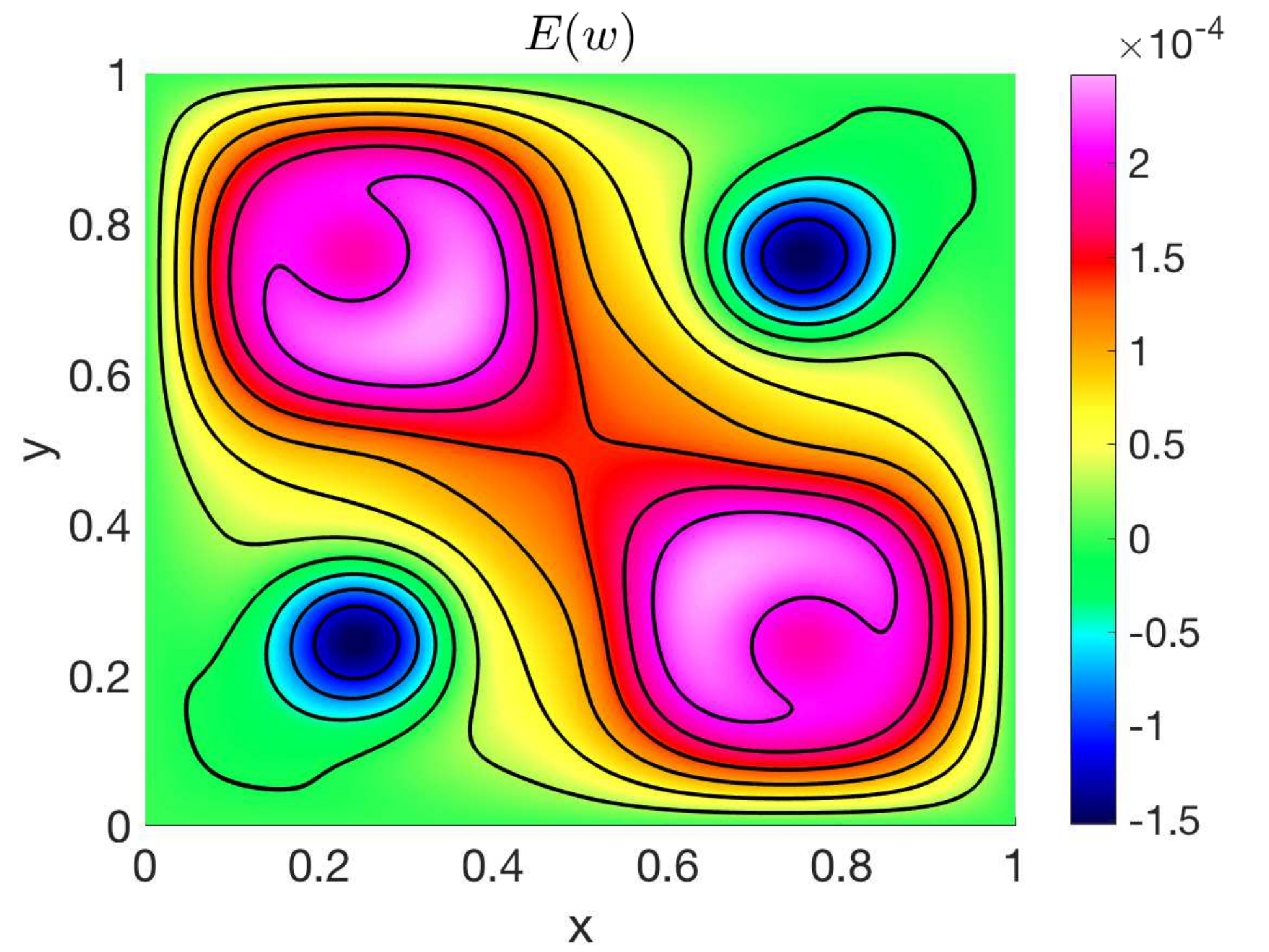}{\figWidth}};
\draw(-0.5,8.0) node[anchor=south west,xshift=0pt,yshift=0pt] {\trimfig{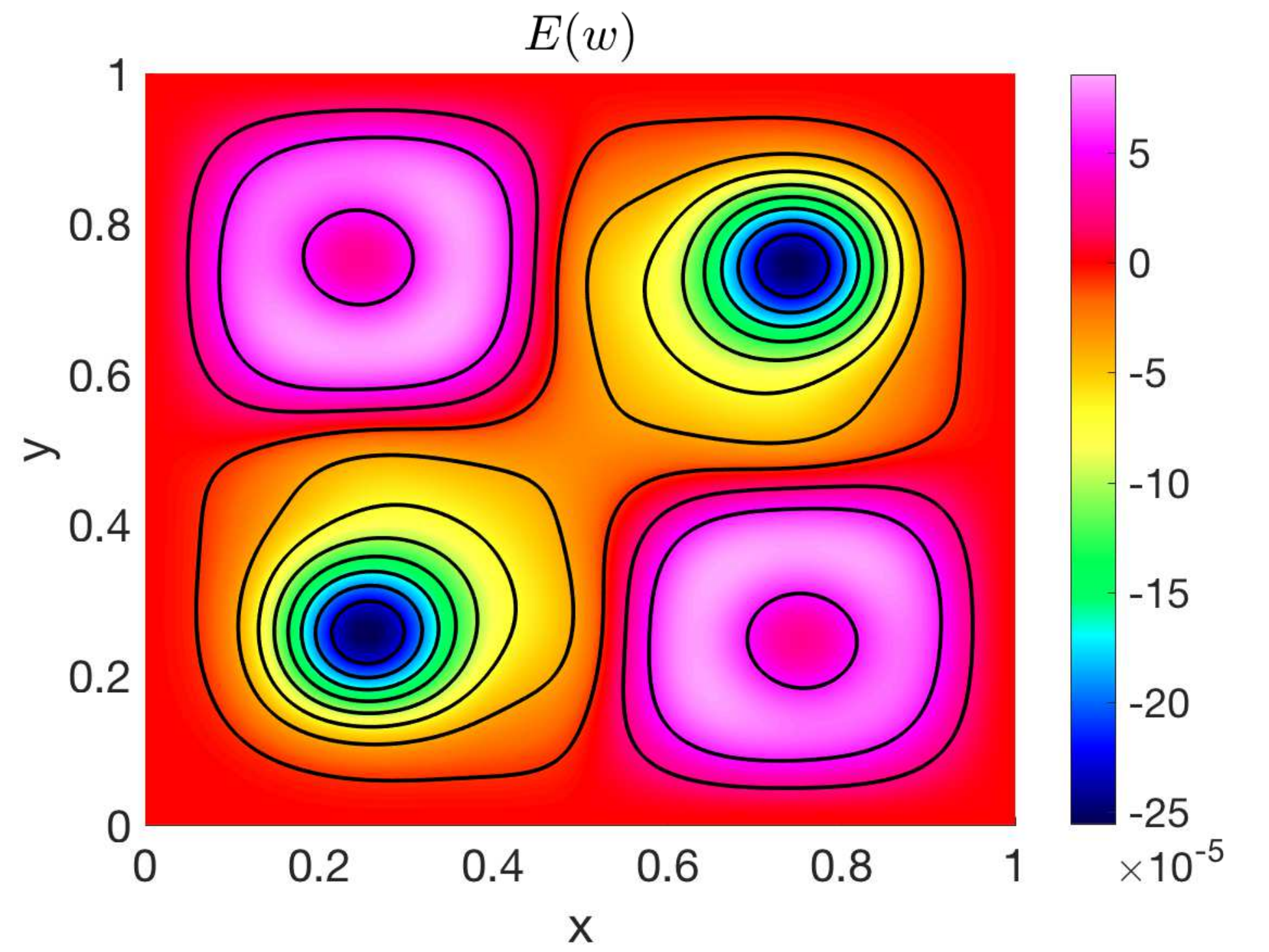}{\figWidth}};
\draw(8.0,8.0) node[anchor=south west,xshift=0pt,yshift=0pt] {\trimfig{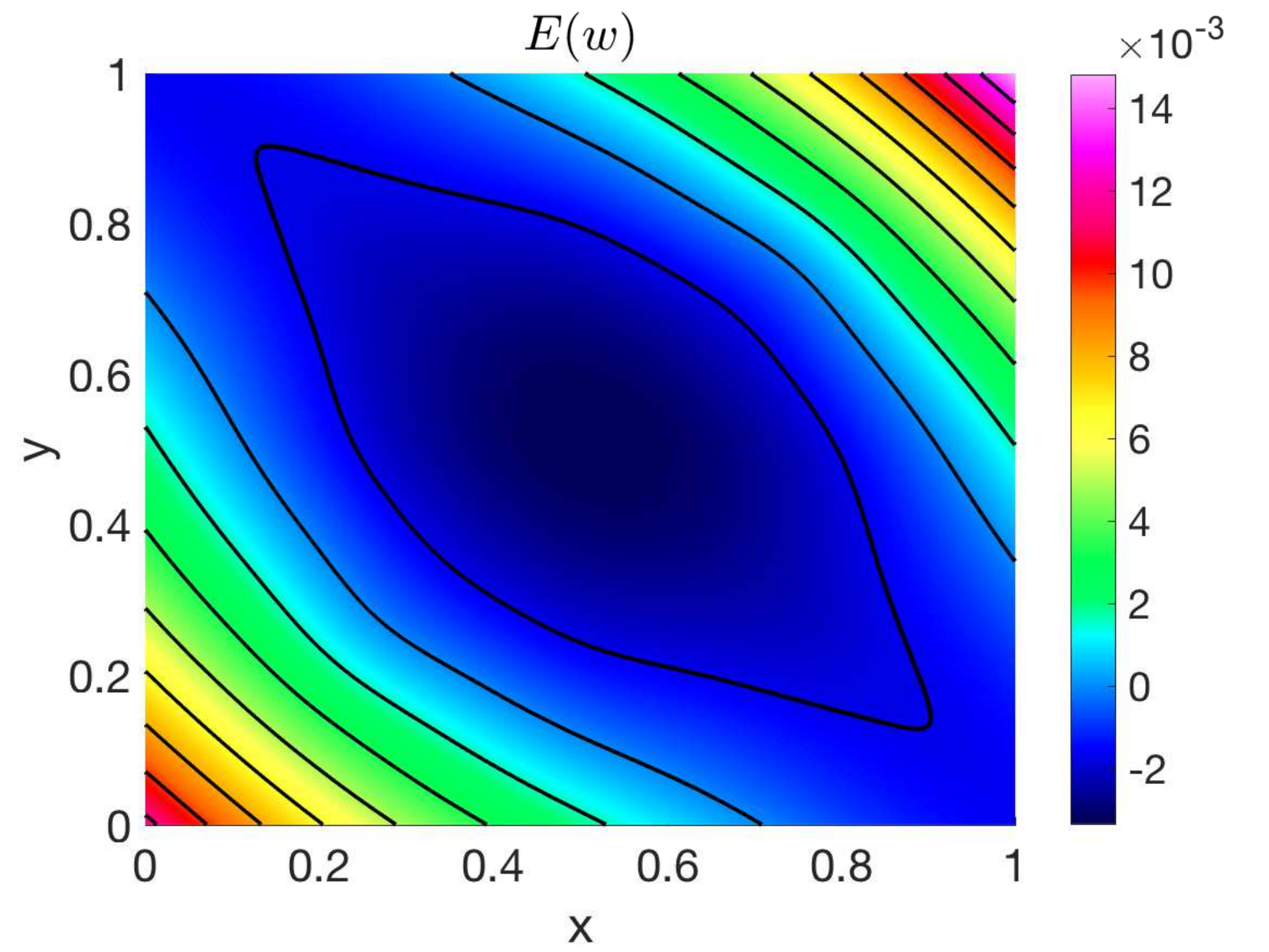}{\figWidth}};
\draw(-0.5,1.0) node[anchor=south west,xshift=0pt,yshift=0pt] {\trimfig{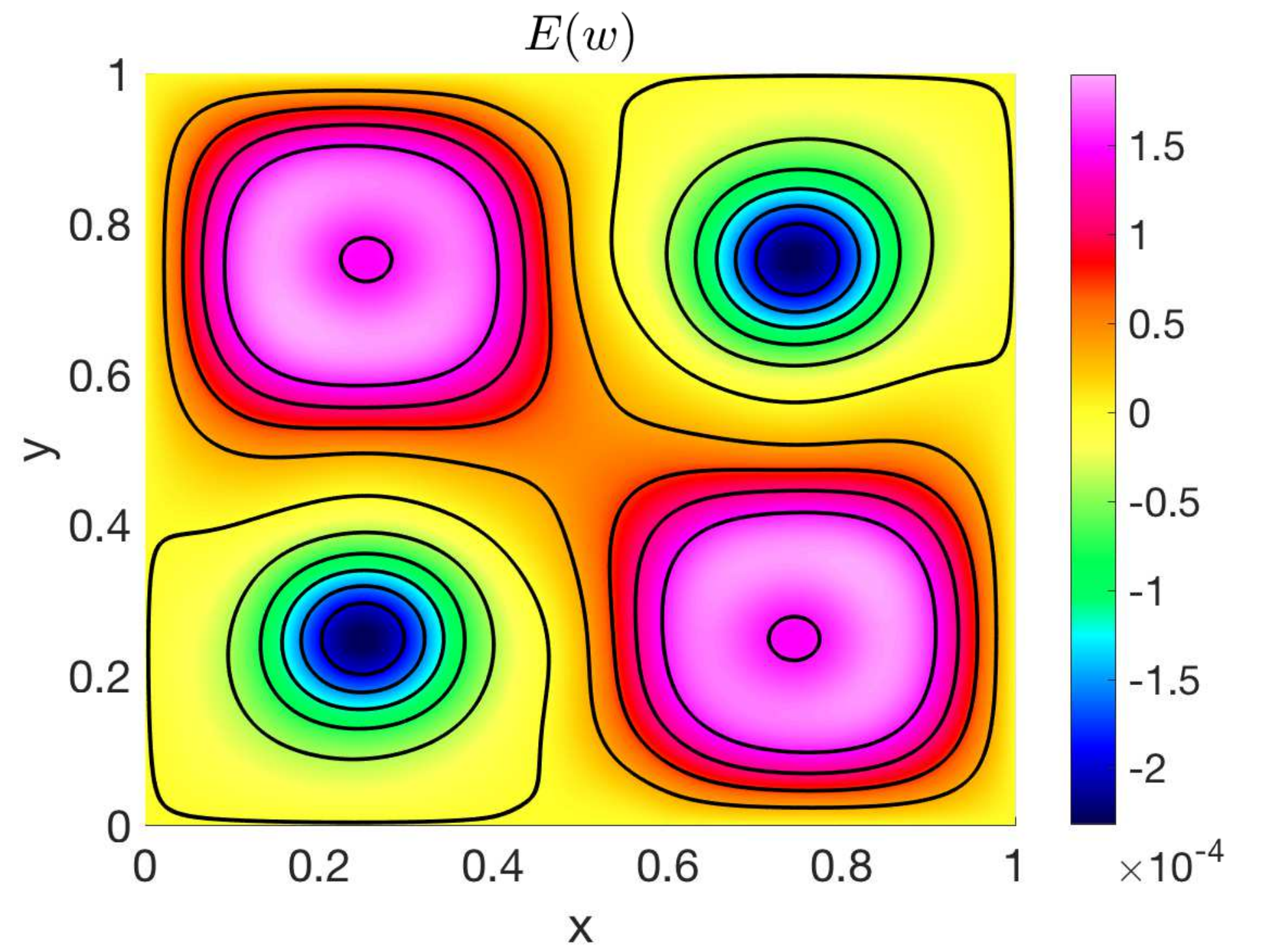}{\figWidth}};
\draw(8.0,1.0) node[anchor=south west,xshift=0pt,yshift=0pt] {\trimfig{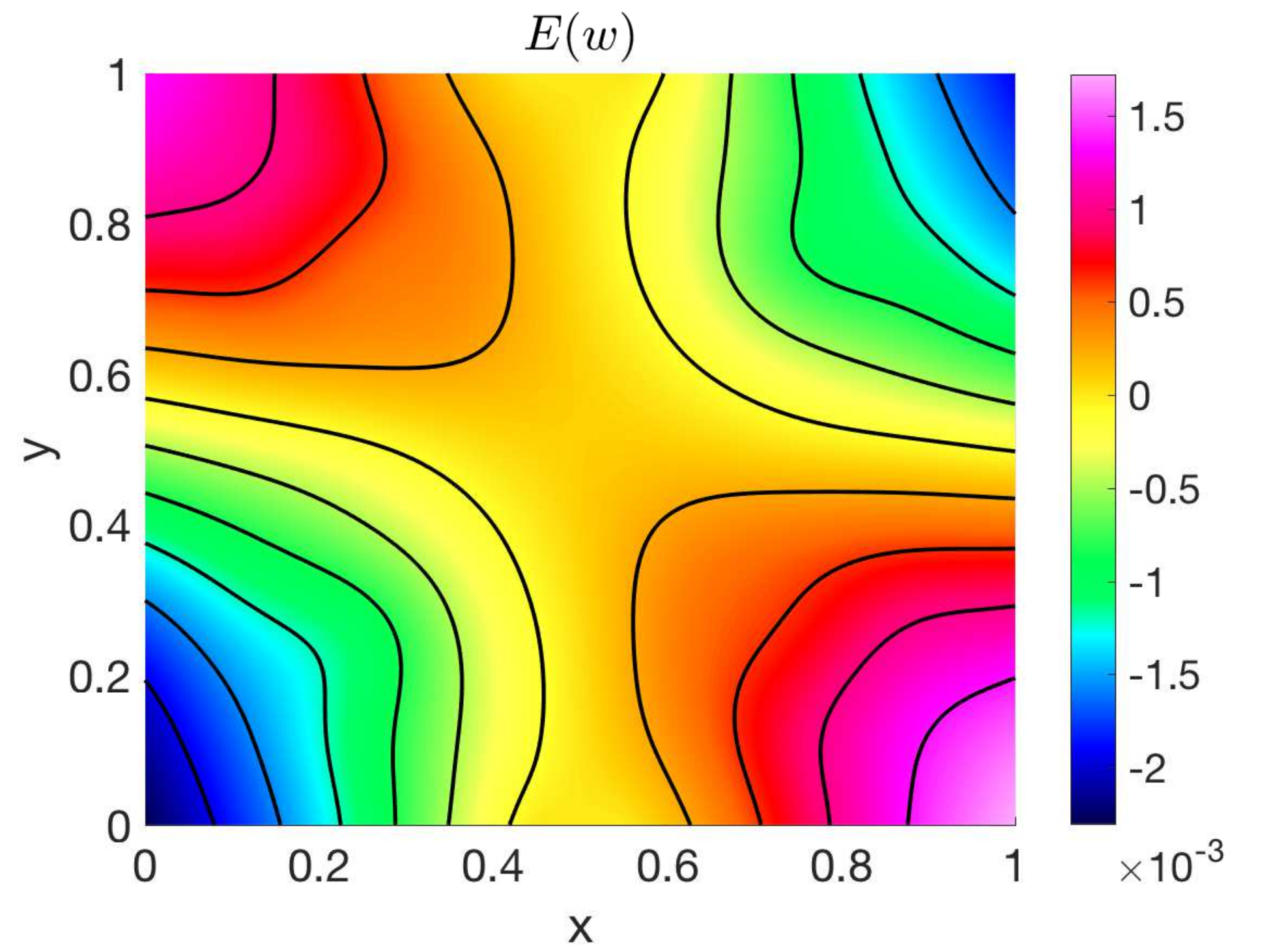}{\figWidth}};

\draw(4,15) node[anchor=north]{\footnotesize(a) Solution};
\draw(13,15) node[anchor=north]{\footnotesize(b) Simply supported BC};
\draw(4,8) node[anchor=north]{\footnotesize(c) Clamped BC};
\draw(13,8) node[anchor=north]{\footnotesize(d) Free BC};
\draw(4,1) node[anchor=north]{\footnotesize(e) Clamped-Supported (CS) BC};
\draw(13,1) node[anchor=north]{\footnotesize(f) Clamped-Free (CF) BC};
%
\end{tikzpicture}
\caption{Contour plots showing $w$ component of the solution and errors of the nonlinear coupled system \eqref{eq:coupledSystemNonlinear} with various boundary conditions  on grid $\mathcal{G}_{640}$. Tolerance for this simulation is $tol=10^{-6}$. Picard method is used here. Results of the other two algorithms are similar.  }\label{fig:nonlinearCoupledTestResultContour_w}
\end{center}
\end{figure}
}

{
\newcommand{\figWidth}{8cm}
\newcommand{\trimfig}[2]{\trimw{#1}{#2}{0.}{0.}{0.}{0.0}}
\begin{figure}[h!]
\begin{center}
\begin{tikzpicture}[scale=1]
\useasboundingbox (0.0,0.0) rectangle (17.,14);  
\draw(-0.5,7.5) node[anchor=south west,xshift=0pt,yshift=0pt] {\trimfig{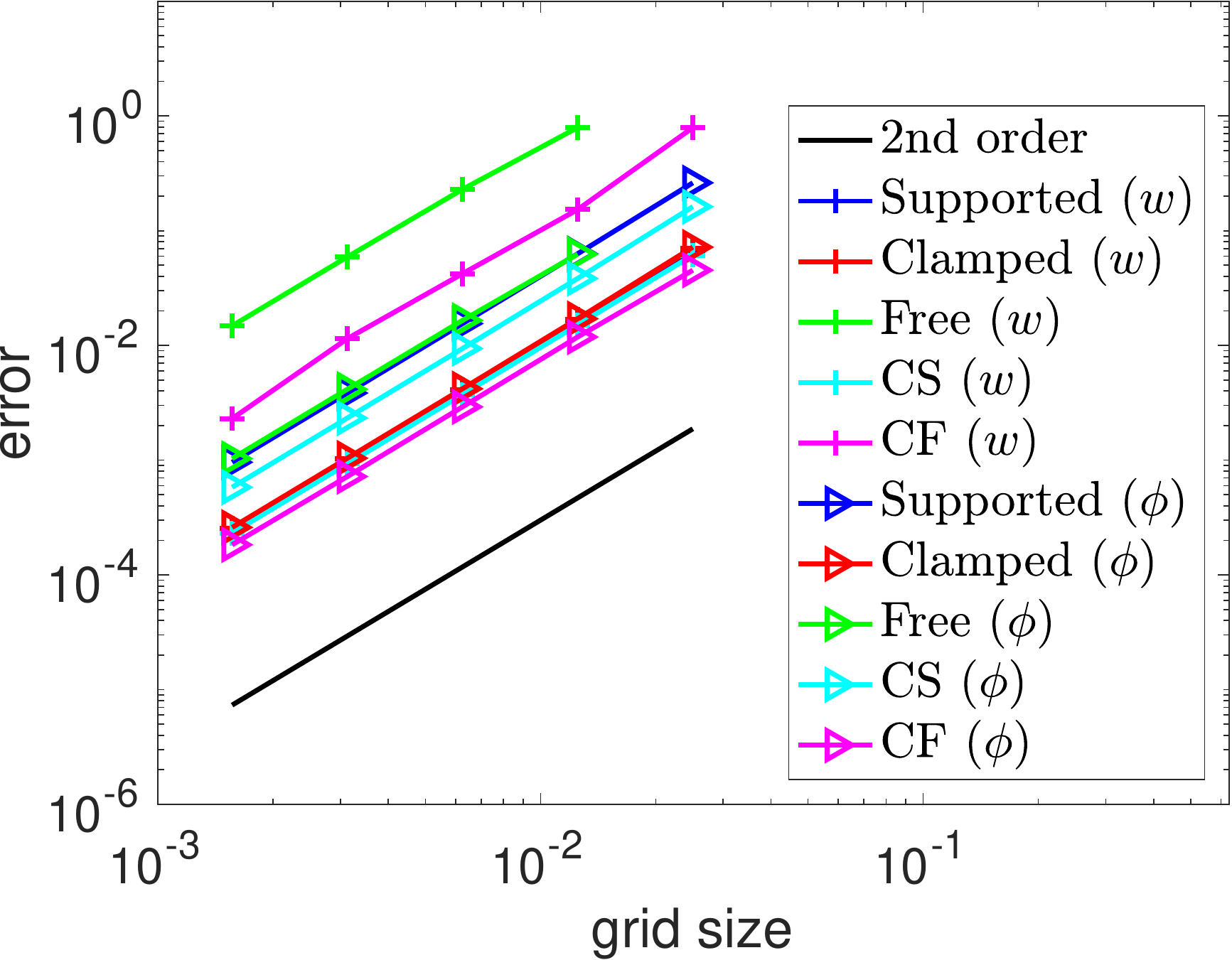}{\figWidth}};
\draw(8.0,7.5) node[anchor=south west,xshift=0pt,yshift=0pt] {\trimfig{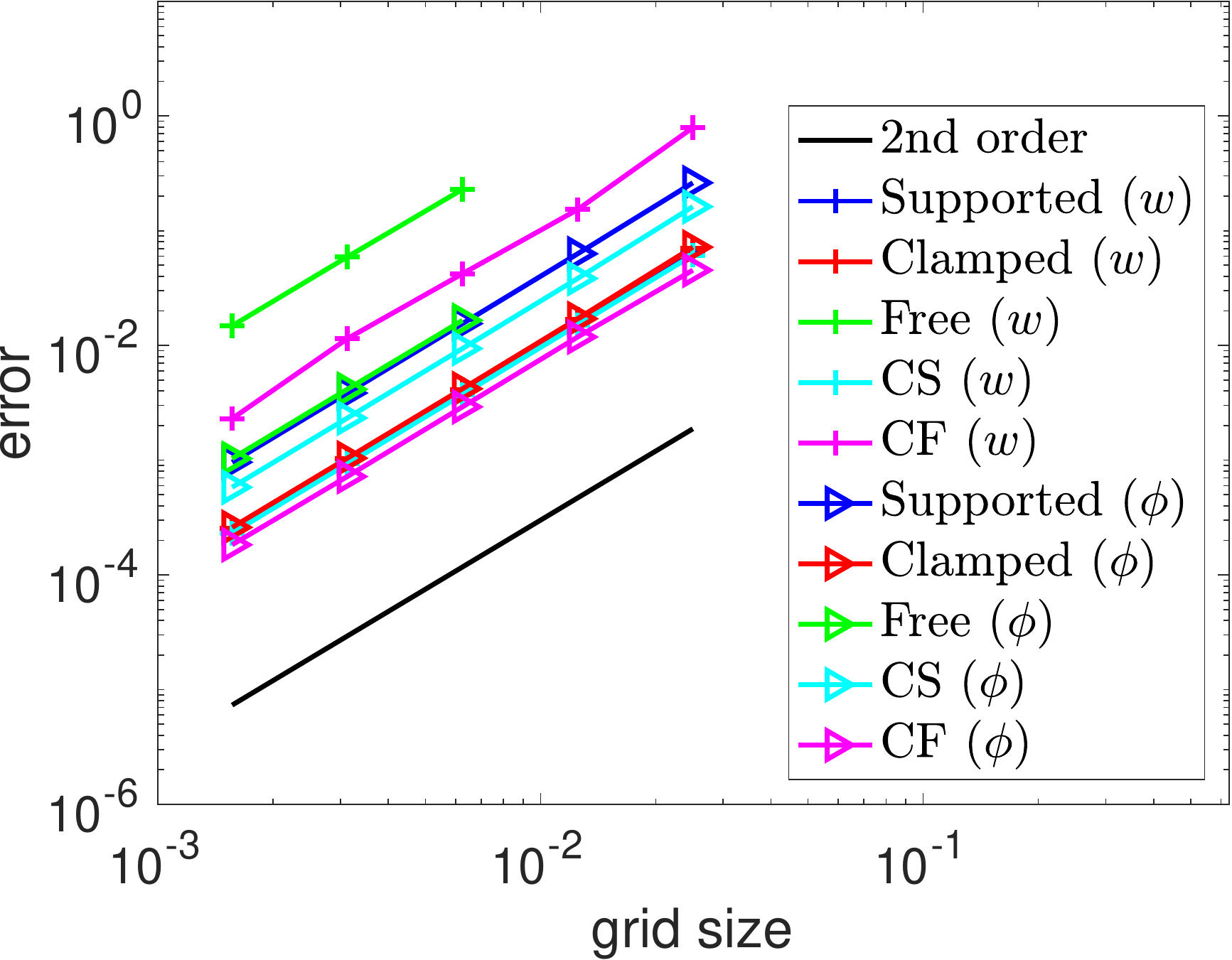}{\figWidth}};
\draw(-0.5,0.5) node[anchor=south west,xshift=0pt,yshift=0pt] {\trimfig{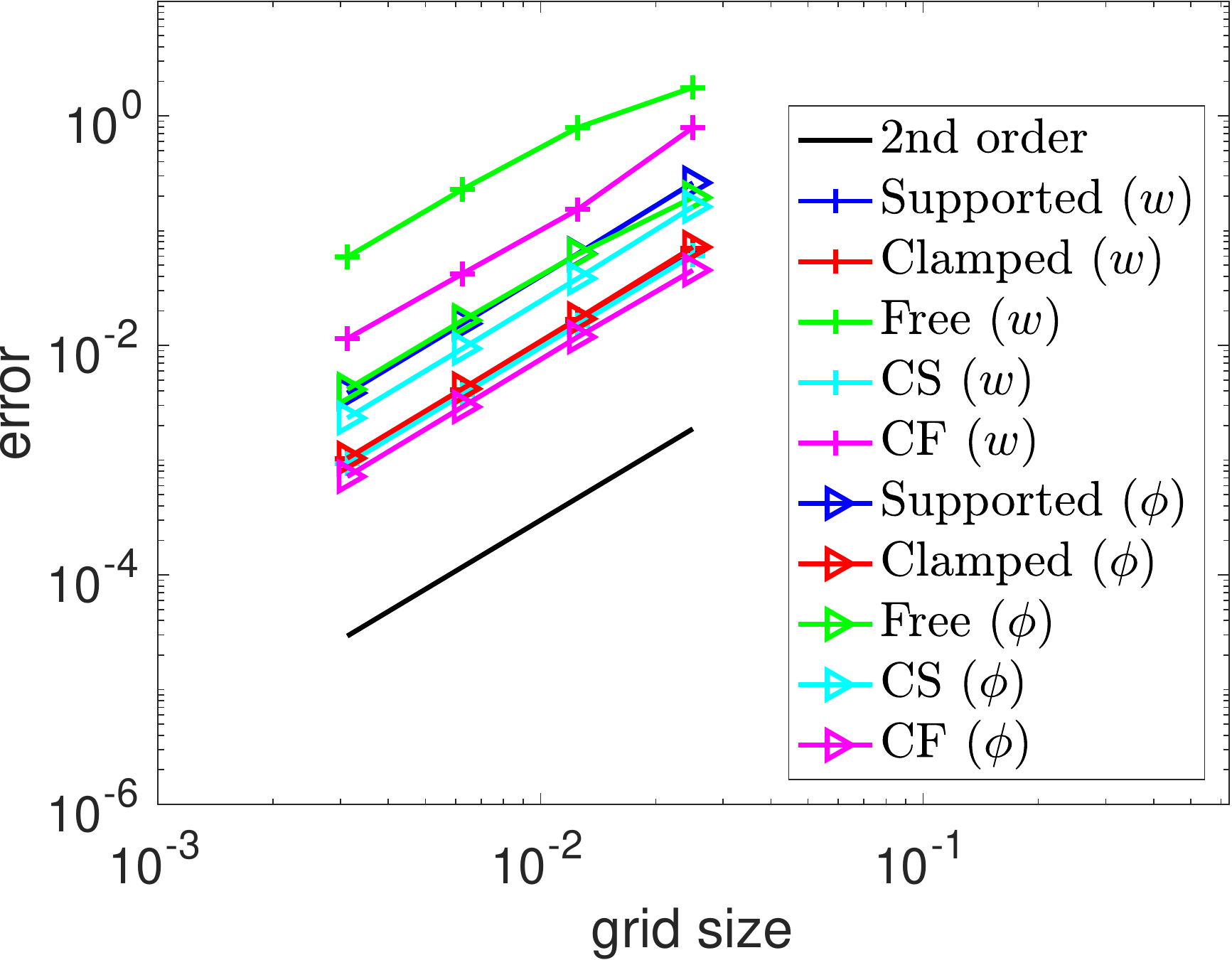}{\figWidth}};
\draw(8.0,0.5) node[anchor=south west,xshift=0pt,yshift=0pt] {\trimfig{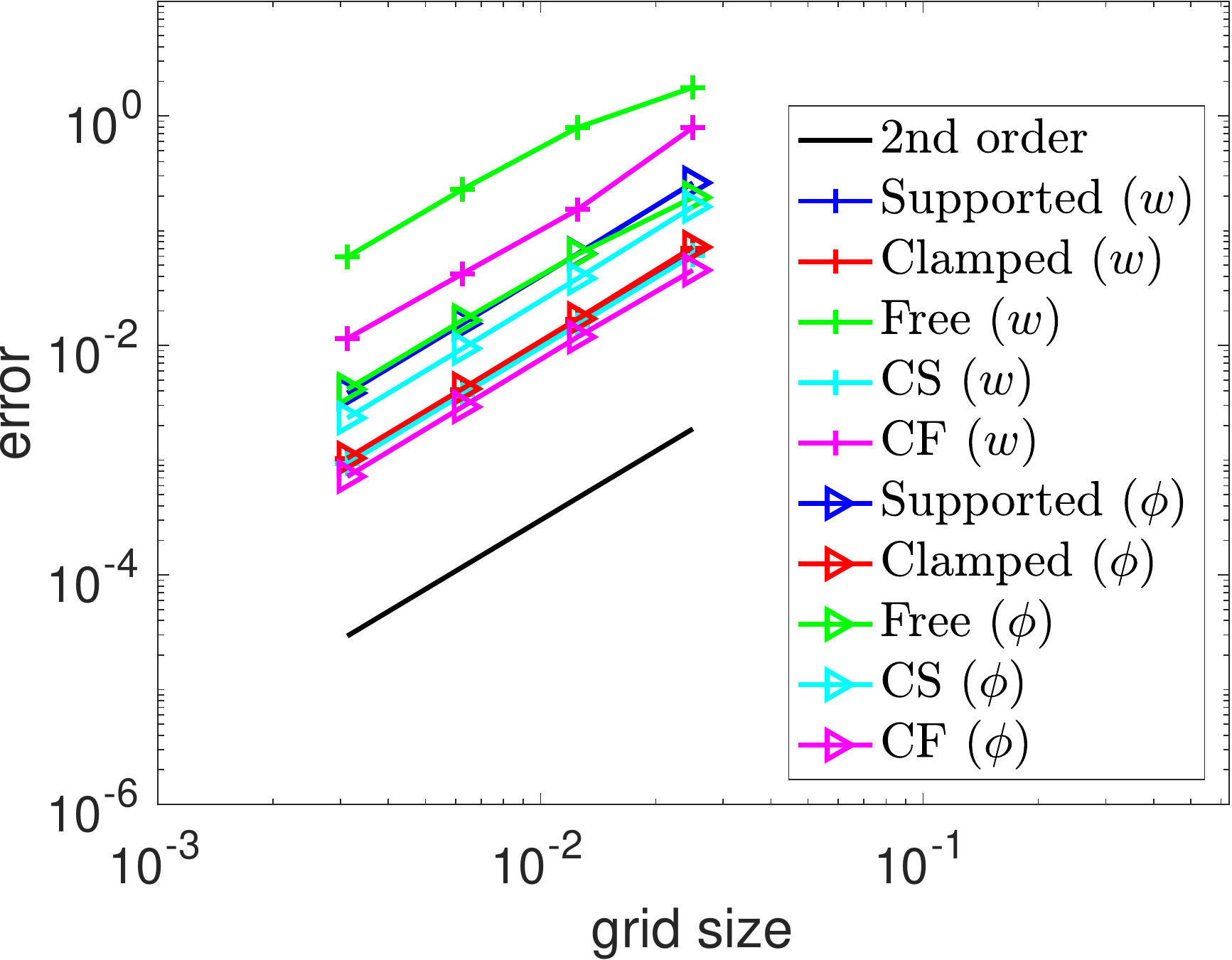}{\figWidth}};

\draw(4,7.5) node[anchor=north]{\footnotesize(a) Picard method (explicit; i.e., $\delta=0$)};
\draw(13,7.5) node[anchor=north]{\footnotesize(b) Picard method (implicit; i.e., $\delta=1$)};
\draw(4,0.5) node[anchor=north]{\footnotesize(c) Newton's method};
\draw(13,0.5) node[anchor=north]{\footnotesize(d) fsolve};
%
\end{tikzpicture}
\caption{A mesh refinement study for the numerical solutions of the nonlinear coupled system with various methods and boundary conditions. Errors are in  maximum norm $L_{\infty}$. Tolerance for this simulation is $tol=10^{-6}$.}\label{fig:convRateNonlinearCoupledSystem}
\end{center}
\end{figure}
}

We end the mesh refinement study by providing a rough comparison between the run-time CPU costs for the explicit Picard method ($\delta=0$) and  the Newton's method on grids with increasing resolutions.  In this comparison, the solving process of both methods are profiled and  their performance for solving the nonlinear test problem with free boundary conditions is summarized in Table~\ref{tab:performance}.   In this table, we list the average CPU time  in seconds-per-step (s/step), the number of steps taken and the estimated rate of convergence; the convergence rate is approximated following \cite{WeerakoonFernando00, CorderoTorregrosa07} 
   by the average  of 
\begin{equation}
p^{k+1} = \frac{\ln\left( ||\Xv^{k+1}-\Xv^{k}||_\infty\right)}{\ln\left( ||\Xv^{k}-\Xv^{k-1}||_\infty\right)}
\end{equation}
for all steps. Recall that, in contrast to the Newton's method which solves one matrix equation involving the Jacobian matrix, the Picard method solves 
two matrix equations that are much smaller in dimension. Therefore,  the Picard method is expected to be faster than the Newton's method per step; and this is observed in Table~\ref{tab:performance}.  However, since the Newton's method converges in 2nd-order rate and the Picard method converges  in  a slower 1st-order one, it takes more steps for the Picard method to converge. Overall speaking, 
the  Newton's method could still beat the Picard method despite being slower at each step; for example, the case with grid $\G_{160}$ in Table~\ref{tab:performance}.  As the grid gets more refined, the Newton's method takes much more time per step  since the size of the Jacobian matrix  becomes too big and eventually surpasses the capacity of our  computational resources.  For the case of $\G_{320}$, the Picard method is faster than the Newton's method both  per step  and  in total time (s/step $\times$ steps); and for the case of $\G_{640}$, the Newton's method encounters out-of-memory issue while the Picard method still works fine.

It is important to note that the performance comparison conducted here is just a rough one, which can be affected by many factors,  the quality of initial guess for instance. In addition, many improvements can be employed to speed up the iteration, such as  Anderson acceleration \cite{WalkerNi11}. The memory issue of the Newton's method can also be alleviated  by switching the solver of the linear system; solvers based on iterative schemes such as the biconjugate gradient stabilized method can be more suitable for problems with large matrix dimension than the direct QR solver used in this paper.   But all those numerical techniques are topics beyond the scope of this paper.

\subsection{Effects of boundary conditions and localized thermal source}
In this section, we solve a realistic problem from industrial application to demonstrate
the effectiveness of our scheme in capturing the influences of thermal stresses, precast shell shape, and various boundary supports to the final shell shape.  
In the process of manufacturing curved glass sheets, the  ``frozen-in'' thermal strain due to a non-ideal cooling history can cause small non-uniformities in the final configuration of the glass sheets. How the thermal stress interplays with the precast shape and boundary conditions is of great interest. 
Motivated by this application, we consider a thin shallow shell with a nonuniform precast  shell shape and various boundary conditions. In order to separate the impact of the thermal stress from the impact of the geometry of the precast shape, we focus on the case where the shell is subjected to a localized thermal loading. The influence of localized thermal heating on plate thermal buckling and post buckling has been investigated in \cite{kumar2017semi} when the plate is either simply supported or clamped. Here, we can do a more thorough study with our new numerical methods.

We assume that there is no external forcing to the shell displacement ($f_w=0$), and prescribe the 
 precast shell shape $w_0$ and the thermal loading $f_{\phi}$ using the following given functions,
\begin{equation}\label{eq:givenFunctionsForNonuniformTest}
w_0= 0.1 -0.4(y -0.5)^2,\quad
f_\phi=32634.2\max\{-100.0((x -0.75)^2 + (y -0.25)^2) + 1,0\}.
\end{equation}
The contour plots of $w_0$ and $f_\phi$ are shown in Fig.~\ref{fig:nonuniformw0andForcing}.
{
\newcommand{\figWidth}{8cm}
\newcommand{\trimfig}[2]{\trimw{#1}{#2}{0.}{0.}{0.}{0.0}}
\begin{figure}[h!]
\begin{center}
\begin{tikzpicture}[scale=1]
\useasboundingbox (0.0,0.0) rectangle (17.,7);  
\draw(-0.5,0) node[anchor=south west,xshift=0pt,yshift=0pt] {\trimfig{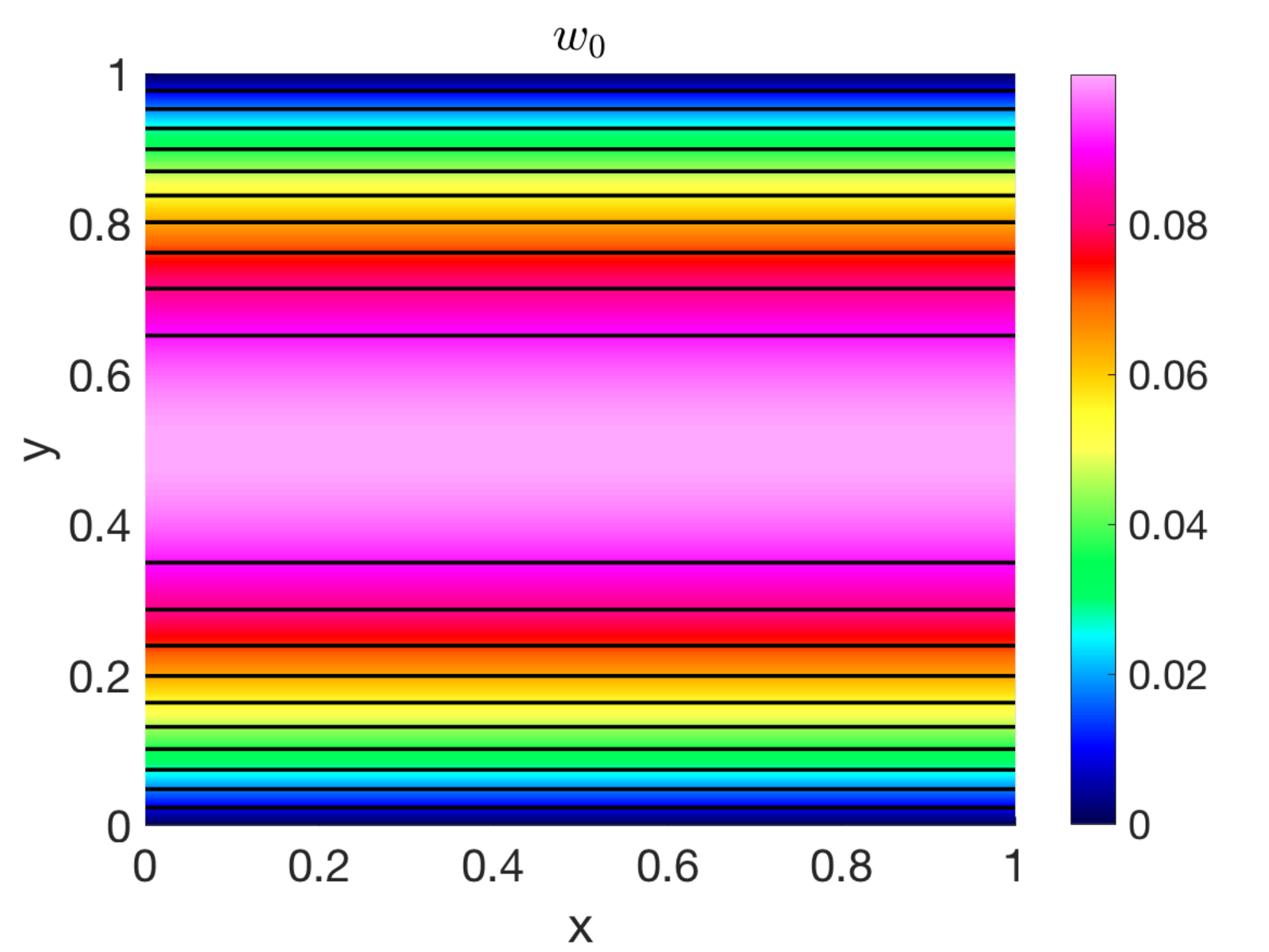}{\figWidth}};
\draw(8.0,0) node[anchor=south west,xshift=0pt,yshift=0pt] {\trimfig{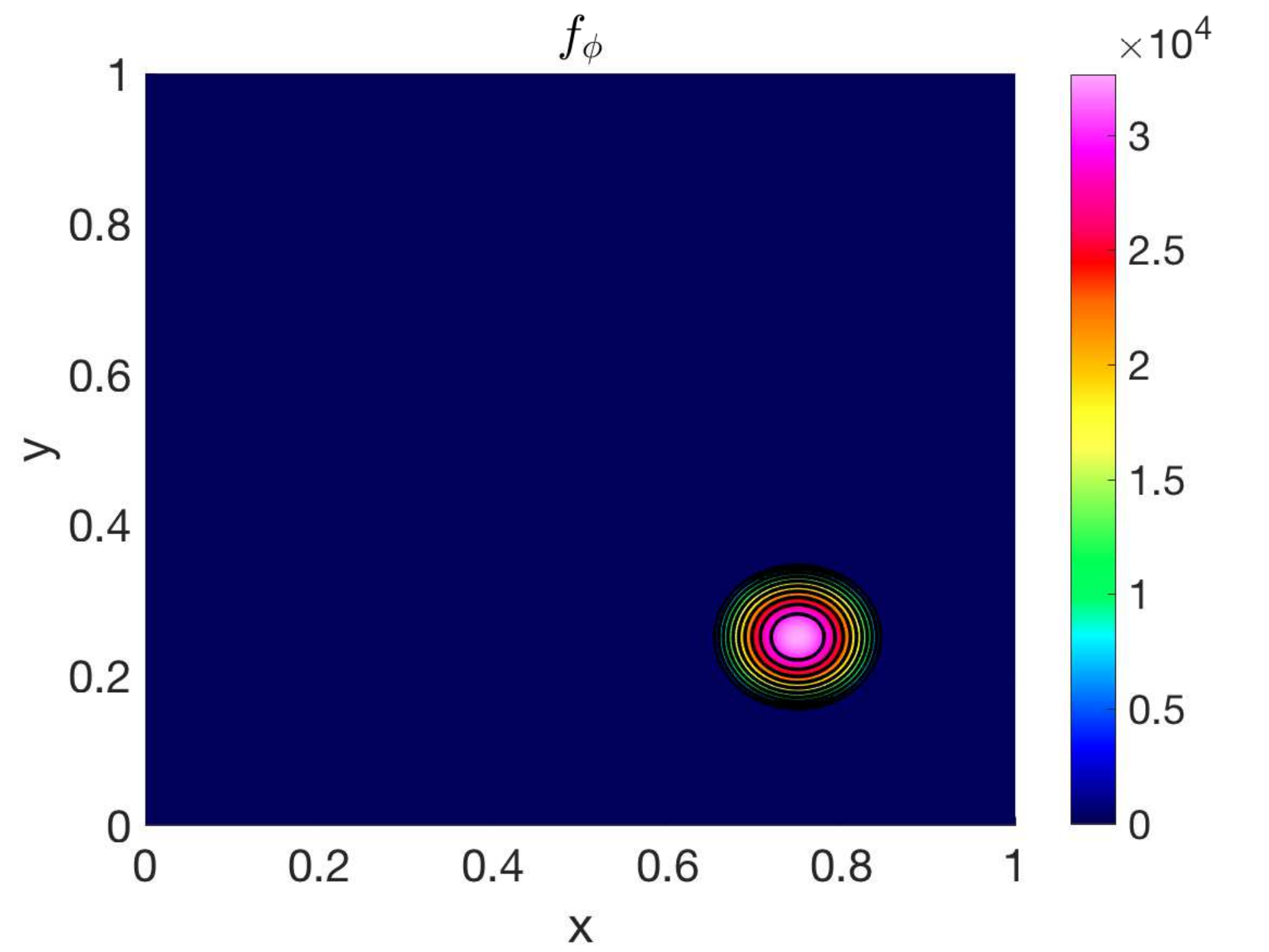}{\figWidth}};

%
\end{tikzpicture}
\caption{Contour plots of (left) the nonuniform precast shell shape and (right) the localized thermal loading.}\label{fig:nonuniformw0andForcing}
\end{center}
\end{figure}
}

This problem is solved again using all three proposed numerical methods subject to all five boundary conditions. From the numerical results, we observe that
while the unstressed shell shape $w_0$  defined in \eqref{eq:givenFunctionsForNonuniformTest} is cylindrically symmetric, both the displacement $w$ and Airy stress function $\phi$ associated with all the boundary conditions are asymmetric due to the effects of the localized thermal forcing at the bottom right corner of the domain (see the right image of  Fig.~\ref{fig:nonuniformw0andForcing}). This suggests that thermal effects can  cause small  non-uniformity in the  shell.

For this problem, we are interested in  how the final shell shape  ($w+w_0$) is affected by  various boundary conditions. In Fig.~\ref{fig:nonuniformWPlusW0}, we collect our numerical solutions for  the final shape subject to these boundary conditions.  We note that,  since the results of all the numerical methods are similar,  only the ones from the Newton's method are presented here.   From  Fig.~\ref{fig:nonuniformWPlusW0}, it can be seen that boundary conditions have a significant impact on the final shell shapes. In particular, we   observe that the free boundary conditions (Fig.~\ref{fig:nonuniformWPlusW0}(c)), as well as the CF boundary conditions  (Fig.~\ref{fig:nonuniformWPlusW0}(e)), introduce the largest deflections to the shell shape; nonetheless, the clamped boundary conditions (Fig.~\ref{fig:nonuniformWPlusW0}(b)) preserve the precast shell shape $w_0$ the best. Our numerical results for the clamped and supported boundary conditions agree qualitatively with the results reported in  \cite{kumar2017semi} where a similar problem with rectangular plates under localized thermal stresses subject to these boundary conditions is investigated. The results for the other three boundary conditions, which are made possible by our numerical methods,  are not available in literature for comparison.

{
\newcommand{\figWidth}{6cm}
\newcommand{\trimfig}[2]{\trimw{#1}{#2}{0.}{0.}{0.}{0.0}}
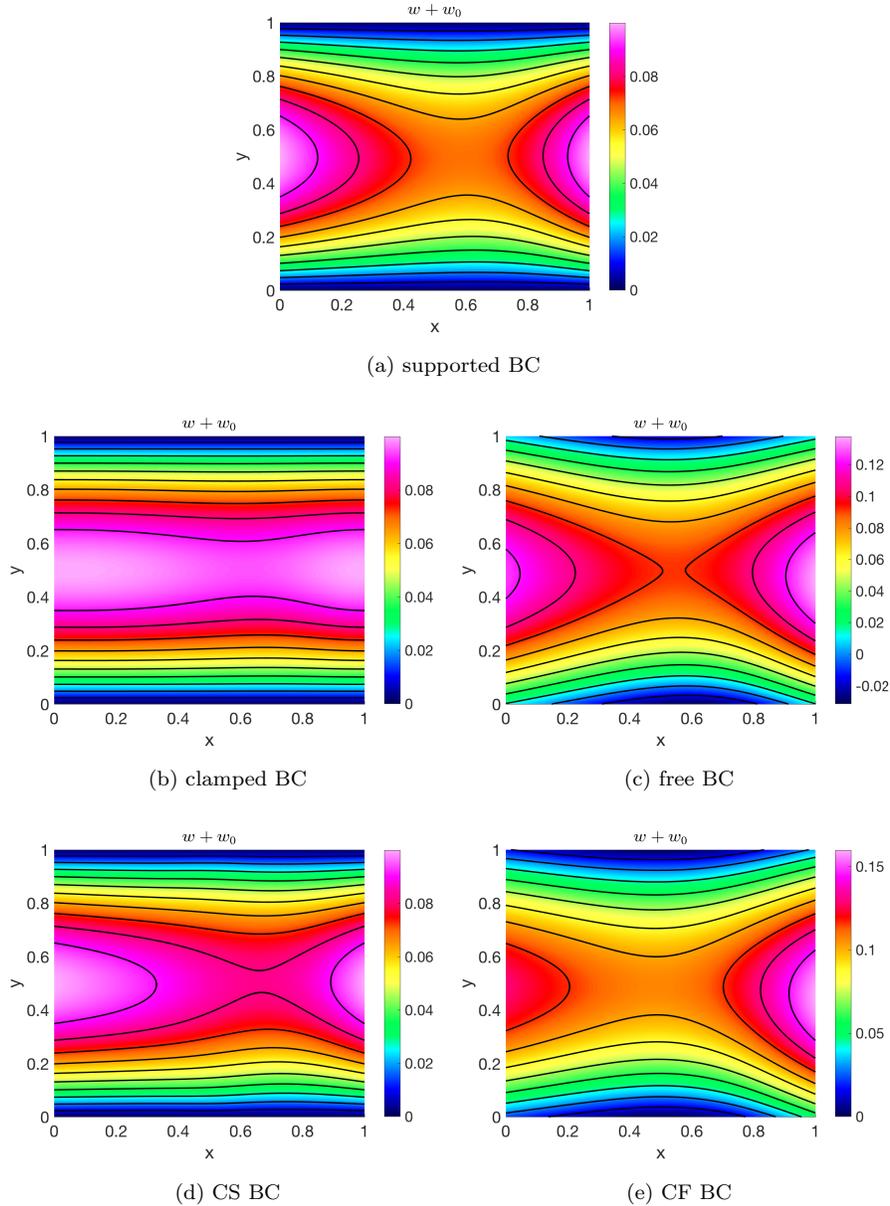
\begin{figure}[h!]
\begin{center}
\begin{tikzpicture}[scale=1]
\useasboundingbox (0.0,0.0) rectangle (12.,16);  

\begin{scope}[yshift=0.5cm]
\draw(2.5,11) node[anchor=south west,xshift=0pt,yshift=0pt] {\trimfig{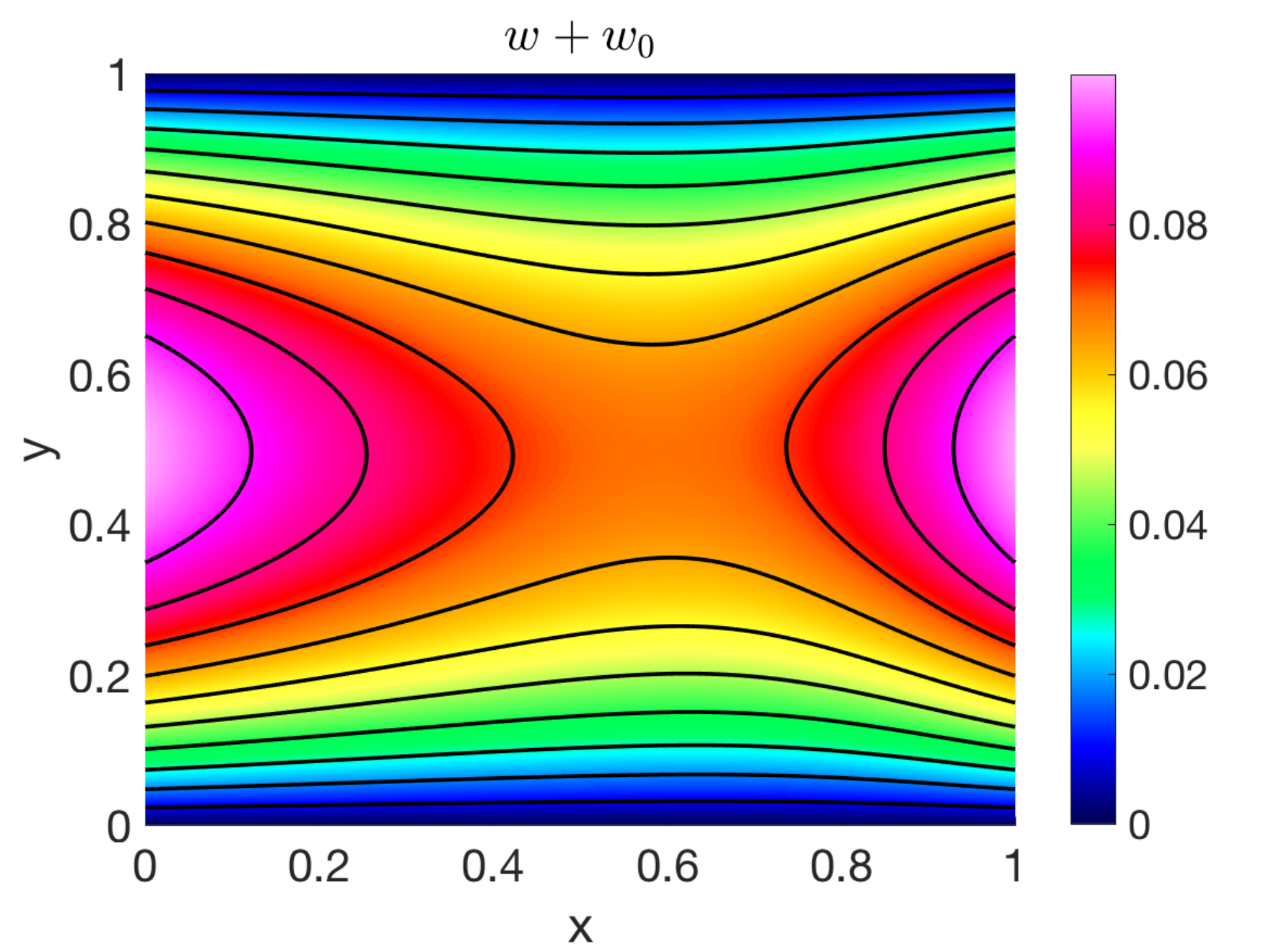}{\figWidth}};
\draw(-0.5,5.5) node[anchor=south west,xshift=0pt,yshift=0pt] {\trimfig{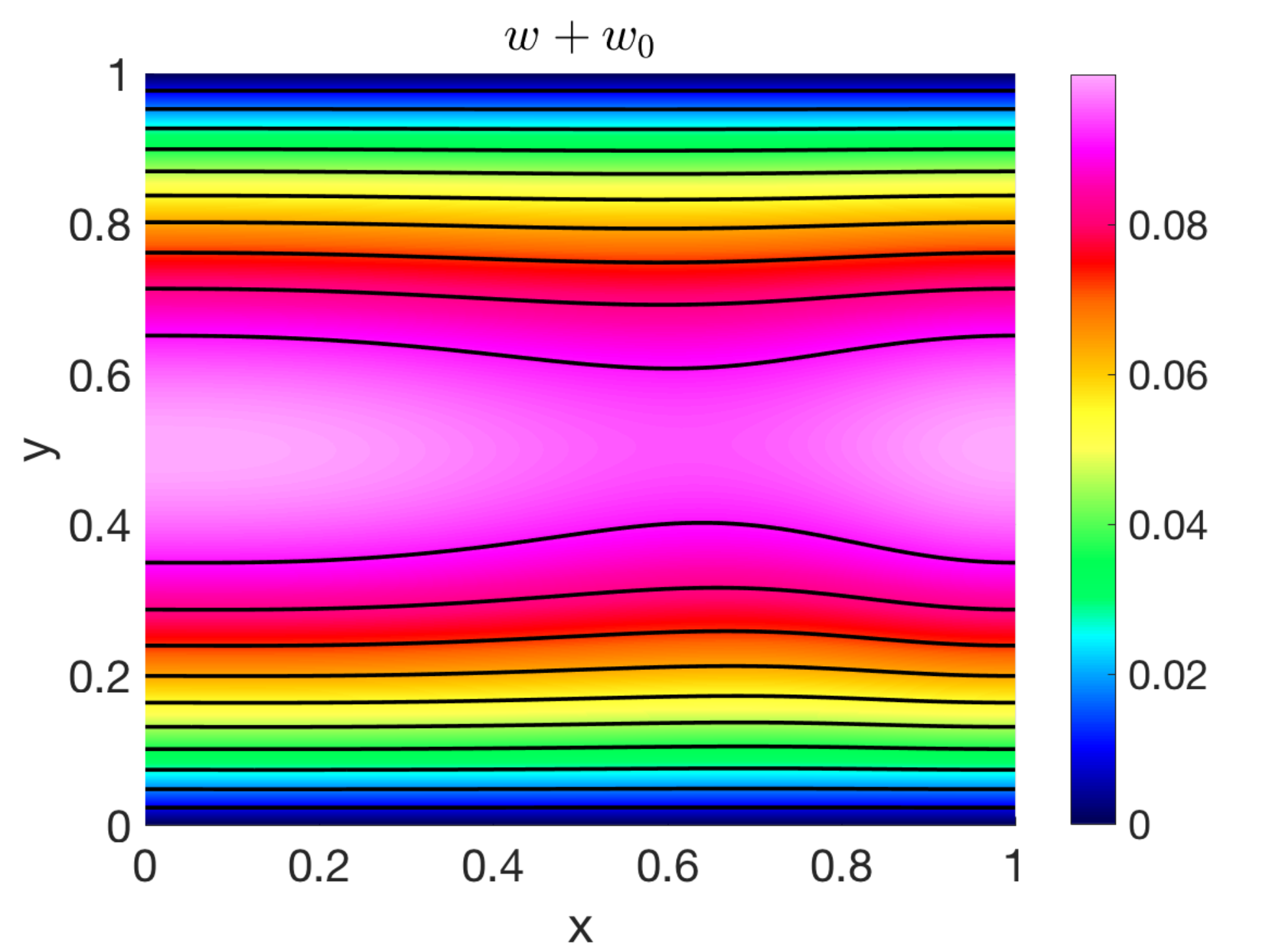}{\figWidth}};
\draw(5.5,5.5) node[anchor=south west,xshift=0pt,yshift=0pt] {\trimfig{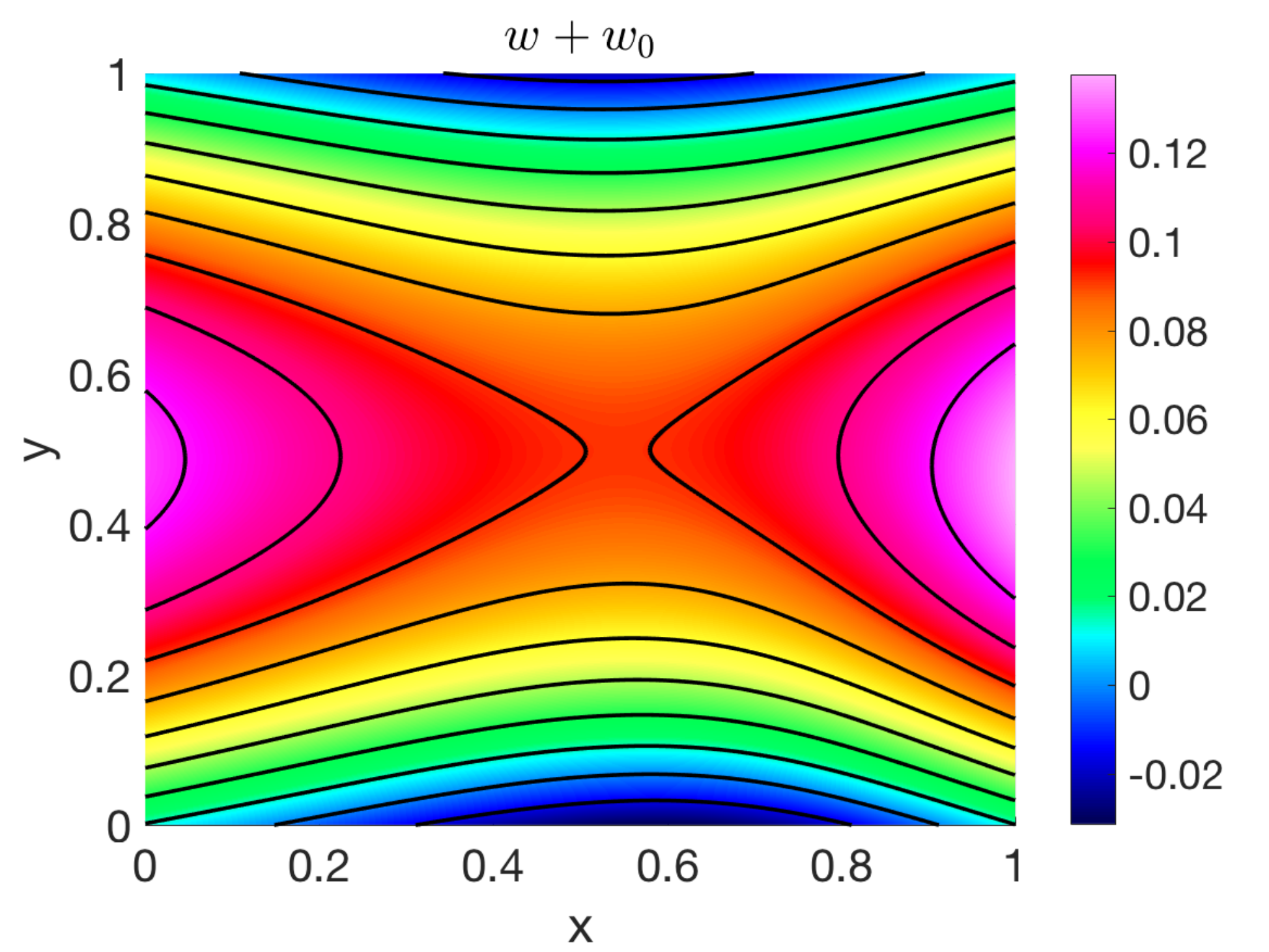}{\figWidth}};
\draw(-0.5,0) node[anchor=south west,xshift=0pt,yshift=0pt] {\trimfig{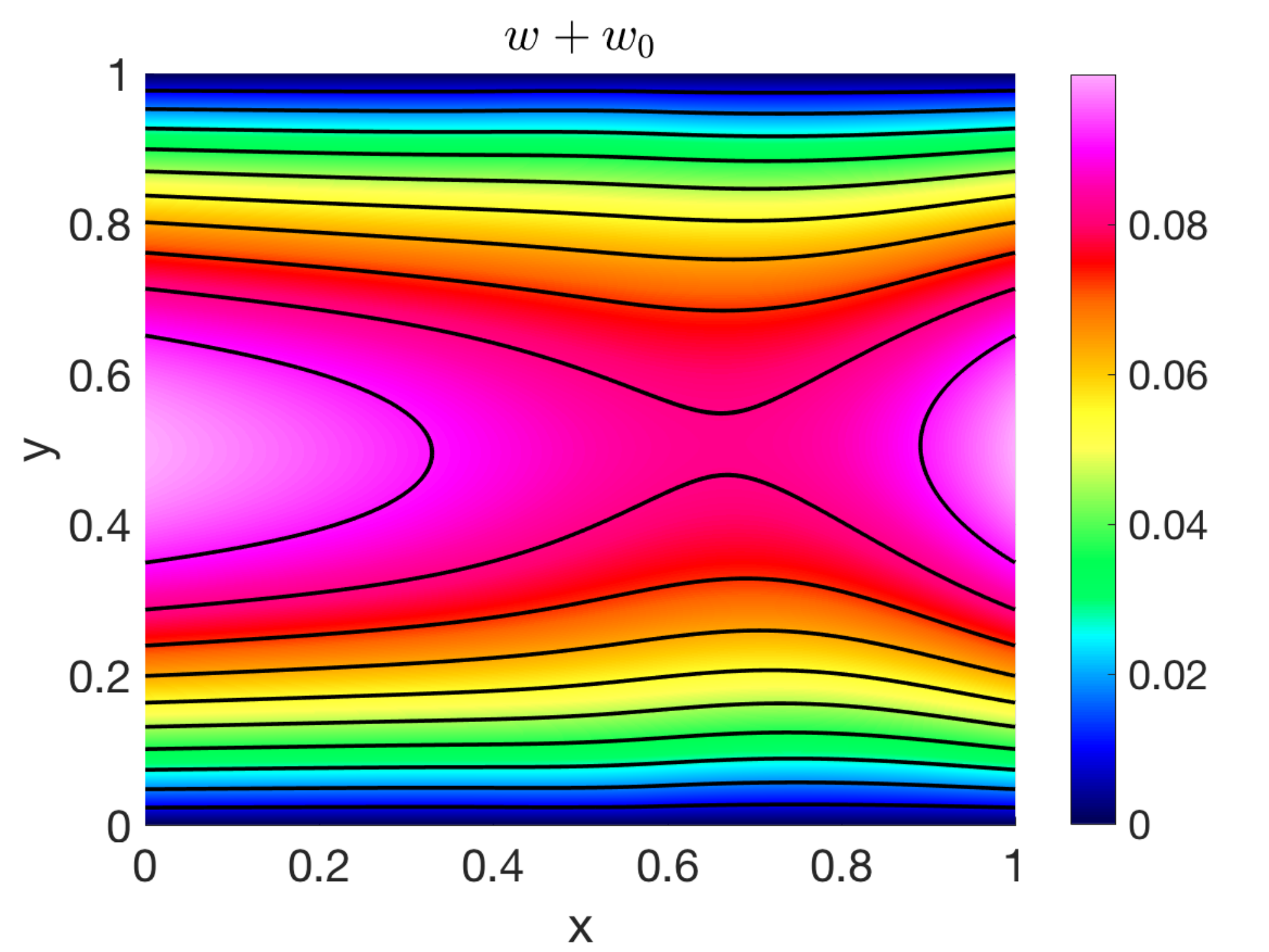}{\figWidth}};
\draw(5.5,0) node[anchor=south west,xshift=0pt,yshift=0pt] {\trimfig{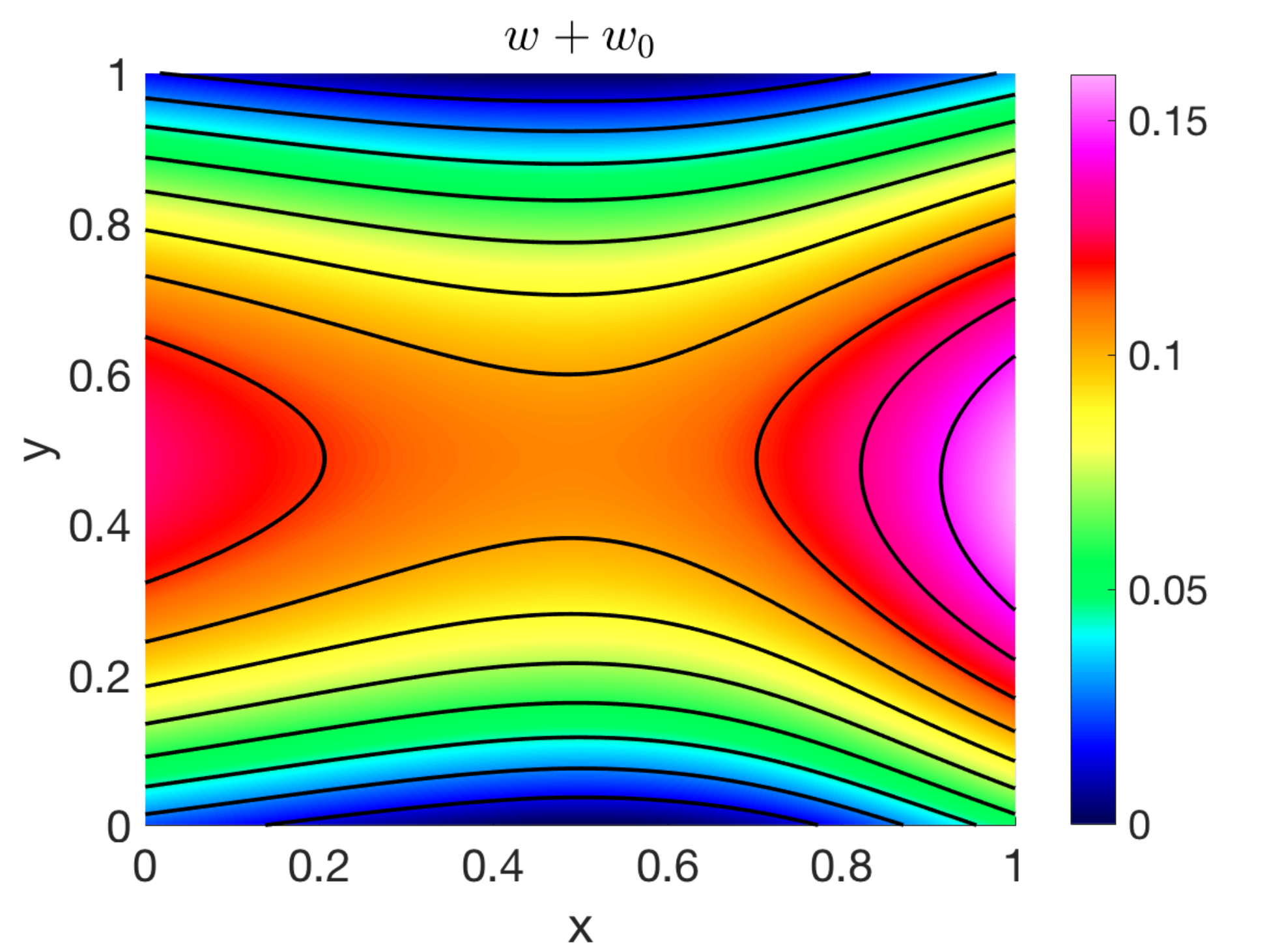}{\figWidth}};

\draw(6,11) node[anchor=north]{\footnotesize(a)  supported BC};

\draw(3,5.5) node[anchor=north]{\footnotesize(b) clamped BC};
\draw(9,5.5) node[anchor=north]{\footnotesize(c) free BC};
\draw(3,0) node[anchor=north]{\footnotesize(d) CS  BC};
\draw(9,0) node[anchor=north]{\footnotesize(e) CF BC};

\end{scope}
%
\end{tikzpicture}
\caption{Contour plots showing the final shape $w+w_0$ governed by the nonlinear shell equations \eqref{eq:coupledSystemNonlinear} with the precast shell shape $w_0$ and thermal loading $f_{\phi}$ in \eqref{eq:givenFunctionsForNonuniformTest} subject to various boundary conditions.}
\label{fig:nonuniformWPlusW0}
\end{center}
\end{figure}
}

\subsection{Snap-through bifurcations}
The critical thermal loading for the snap-through bifurcation is important for the understanding of the maximum allowed temperature for a shell or plate structure. In \cite{tauchert1991thermally, thornton1993thermal},   thermal buckling  with various types of shell shapes and boundary supports has been studied. It has been shown in  \cite{murphy2001thermal, mahayni1966thermal} that the deflection of a perfectly flat plate develops a symmetric pitchfork bifurcation associated with elevated temperatures, while a shallow shell undergoes an asymmetric saddle-node bifurcation at a relatively high critical temperature.  
As an application of our proposed  numerical methods for the nonlinear shallow shell equations,  the snap-through thermal buckling problem is studied  numerically for each of the boundary conditions  so as to  demonstrate the effectiveness and accuracy of our numerical methods.

Specifically,   we solve the nonlinear shallow shell equations \eqref{eq:coupledSystemNonlinear}   with the following specifications
\begin{align}\label{eq:bifurcationProblemSpecification}
f_\phi(x,y)=\xi,\quad
f_w(x,y)=0,\quad
w_0(x,y) = 0.3(1-(x-0.5)^2-(y-0.5)^2),
\end{align}
where $\xi$ is a spatially-uniform thermal loading.  The snap-through bifurcation can be obtained numerically using path following (parameter continuation) techniques. The natural choice of parameter for continuation in this problem would be the constant thermal loading $\xi$. The idea of parameter
continuation is to find a solution to the governing equations at $\xi_0+\Delta\xi$ for a small perturbation $\Delta\xi$ given the solutions at $\xi_0$, and then we  proceed step by step to get a global solution path; solutions at each step are solved using our iterative methods with the solutions from the previous step as the initial guess.  However, it is well-known that
the natural parameter continuation method may fail at some step due to the existence of singularities on the curve (e.g., folds or bifurcation points) \cite{Keller87}.  Even though an estimate of the locations of the folding points can be made from bifurcation branches away from the  singularities that are obtained during the natural parameter continuation prior to failure, we would like to have  a more accurate estimation  since  the folding points represent the critical thermal loading of our problem. In order to capture the critical thermal loading, the so-called Pseudo-Arclength Continuation (PAC) method is used here to circumvent the simple fold difficulties \cite{Keller87}. 

The main idea in PAC is to drop the natural parametrization by $\xi$ and use some other parameterization.  Detail discussion about the PAC method can be found in the lecture notes \cite{Keller87} by Keller.  Here, to be self-contained, we briefly describe the PAC method used for our problem. To simplify notation,
we denote the  shallow shell equation \eqref{eq:coupledSystemNonlinear}  together with the specifications given in \eqref{eq:bifurcationProblemSpecification} as
\begin{equation}\label{eq:PACMainEqn}
G(w,\phi,\xi)=0.
\end{equation}
Instead of tracing out a solution path from the incrementation of the natural parameter,  we treat $\xi$ as an unknown and solve \eqref{eq:PACMainEqn} together with a  scalar
normalization equation,
\begin{equation}\label{eq:PACNormalEqn}
N(w,\phi,\xi;\Delta s) \equiv  \dot{w}_p(w-w_p)+\dot{\phi}_p(\phi-\phi_p)+\dot{\xi}_p(\xi-\xi_p)-\Delta s=0,
\end{equation}
where $N(w,\phi,\xi;\Delta s)=0$ is the equation of a plane perpendicular to the tangent $(\dot{w}_p,\dot{\phi}_p, \dot{\xi}_p)$ at a distance $\Delta s$ from a solution
 $({w}_p,{\phi}_p, {\xi}_p)$. 
This plane will intersect the solution path  if $\Delta s$ and the curvature of the path are not too large. Here $\Delta s$ can be regarded as  the increment of the pseudo-arclength of the path curve.  So  by solving   \eqref{eq:PACMainEqn} \& \eqref{eq:PACNormalEqn} step by step with the solutions $({w}_p,{\phi}_p, {\xi}_p)$ from the previous step as initial guess, we are able to circumvent the simple fold difficulties and  obtain a complete bifurcation branch.  We note that  at each step the equations  \eqref{eq:PACMainEqn} \& \eqref{eq:PACNormalEqn} are solved using the iterative solvers proposed in this paper. In addition, in order to speed up the continuation process,  $\Delta s$ is set dynamically in our numerical implementation; that is, the size of $\Delta s$ is set to decrease when approaching a folding point and increase when leaving one automatically.

%
%

{
\newcommand{\figWidth}{8cm}
\newcommand{\trimfig}[2]{\trimw{#1}{#2}{0.}{0.}{0.}{0.0}}
\begin{figure}[h!]
\begin{center}
\begin{tikzpicture}[scale=1]
\useasboundingbox (0.0,0.0) rectangle (17.,7);  
\draw(-0.5,.5) node[anchor=south west,xshift=0pt,yshift=0pt] {\trimfig{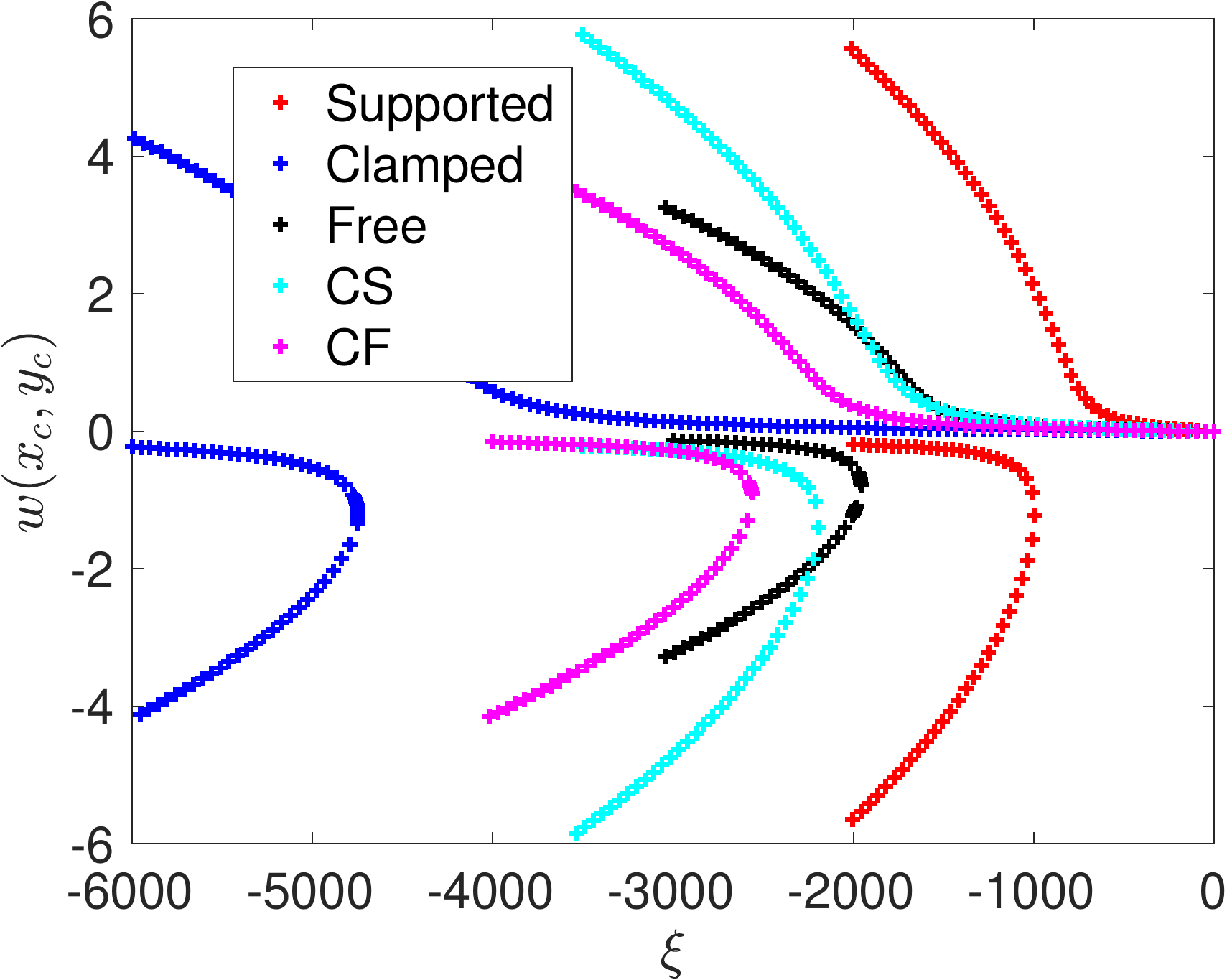}{\figWidth}};
\draw(8.0,.5) node[anchor=south west,xshift=0pt,yshift=0pt] {\trimfig{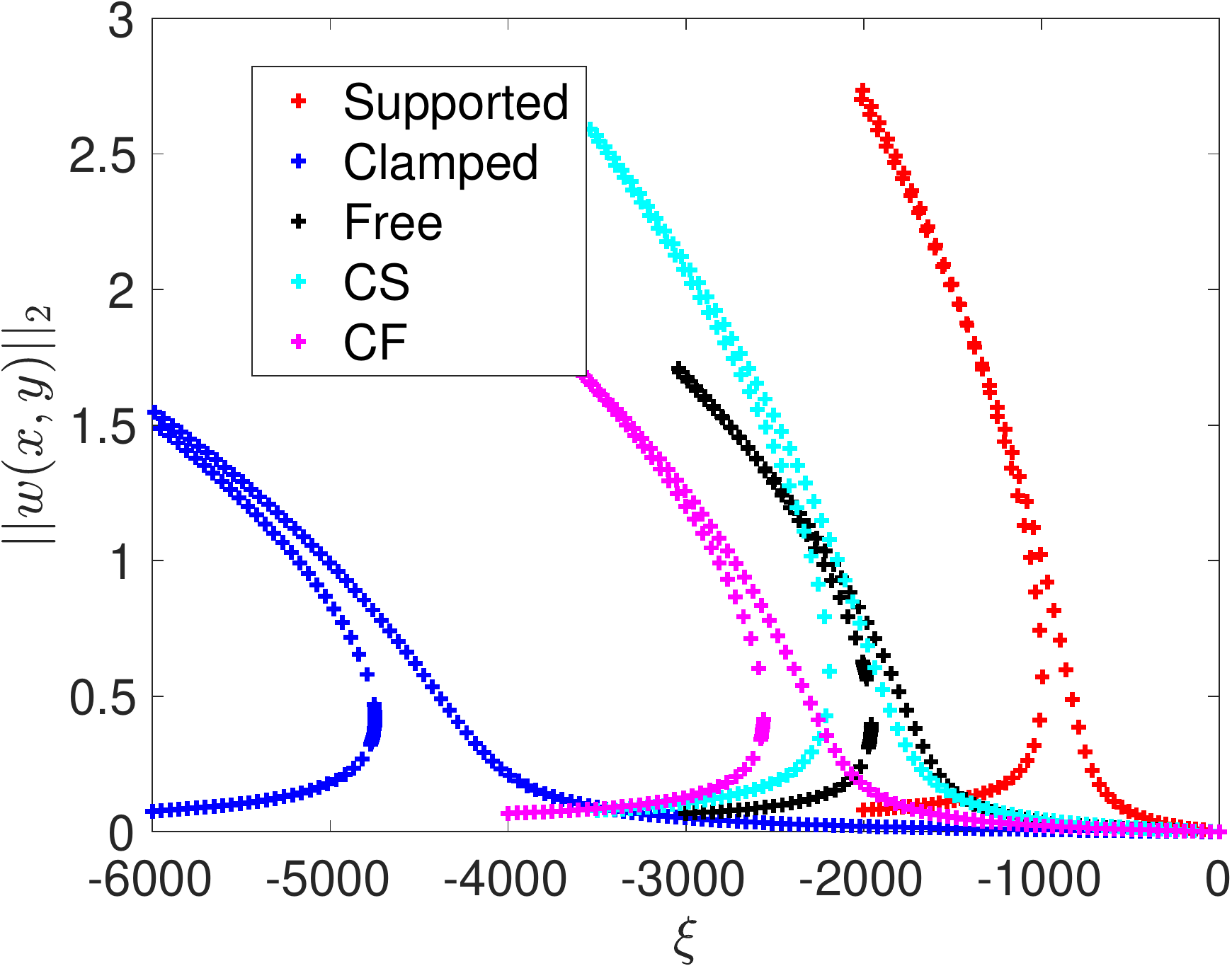}{\figWidth}};

\draw(4,.5) node[anchor=north]{\footnotesize(a) center point};
\draw(13,.5) node[anchor=north]{\footnotesize(b) $L_2$ norm };
%
\end{tikzpicture}
\caption{Snap-through bifurcations for the nonlinear shallow shell equations subject to constant thermal loading and various boundary conditions. Left: bifurcation diagram of the displacement  $w$ at the center of the domain against the thermal loading $\xi$ with various boundary conditions. Right: bifurcation diagram of the $L_2$ norm of the displacement  $w$  against the thermal loading $\xi$ with various boundary conditions.}\label{fig:bifurcation}
\end{center}
\end{figure}
}

The bifurcation results obtained using the PAC method and the Newton's solver for the shallow shell equations are shown in Fig.~\ref{fig:bifurcation}. It is interesting to notice that, by using PAC, there exhibits no issues related to the simple fold difficulties that had plagued the natural parameter continuation and notice that the data points  are clustered towards the folding points due to the  dynamical strategy of setting the value of $\Delta s$.  
 For all the five boundary conditions, we collect the bifurcation diagrams of the displacement at the center of the domain  in Fig.~\ref{fig:bifurcation} (left), and show the bifurcation diagrams in terms of the $L_2$ norm of the shell displacement  in Fig.~\ref{fig:bifurcation} (right). The saddle-node bifurcation curves in Fig.~\ref{fig:bifurcation} (left) all qualitatively agree with the load-deflection curves for thermal buckling of shallow shells in the literature \cite{murphy2001thermal, mahayni1966thermal}.
 
The results reveal the effects of boundary conditions on critical thermal loadings of snap-through buckling.  It is clearly seen from the numerical results that the critical thermal loading for fully clamped boundary supports is much smaller than the other ones, while the clamped  boundary conditions possess  the largest  critical thermal loading. Let $\xi_{\text{bc}}$ denote the critical thermal loading of various boundary conditions, the locations of the folding points in the left plot of Fig.~\ref{fig:bifurcation} indicate the following relation,
$$
\xi_{\text{clamped}} < \xi_{\text{cf}} <\xi_{\text{cs}} < \xi_{\text{free}}< \xi_{\text{supported}}.
$$

\section{Conclusions}

We have developed novel finite difference based iterative  schemes to solve a von-Karman type nonlinear shallow shell model \eqref{eq:coupledSystemNonlinear} that incorporates the thermal stresses. The boundary conditions considered for the system are three simple boundary conditions and two application-motivated mixed boundary conditions. 
To deal with the boundary singularities introduced by the  mixed boundary conditions and maintain the second order accuracy, a transition function approach and a local asymptotic solution approach  are proposed.
All  proposed numerical  methods  for solving the shallow shell equations  are  verified as second order accurate by numerical mesh refinement studies.
 
 As a  demonstration of the efficiency and accuracy of our numerical schemes for engineering applications, we also solve two realistic shallow shell problems; namely the localized thermal source problem and the the snap-through thermal buckling problem. Our numerical results directly reveal  
the combined effects of unstressed shell shape, thermal stresses and boundary conditions on the shallow shell system. In addition,   for the snap-through thermal buckling problem, we are able numerically obtain the snap-through bifurcations using the  pseudo-arclength continuation method with the equations at each continuation step being solved by one of our proposed numerical methods. Our results for both problems are consistent with existing studies.

A number of interesting questions remain to be answered.
While this paper is devoted to the static von Karman shell equations, we are also interested in extending our methods to the corresponding dynamical systems. For instance, in \cite{bilbao2008family} Bilbao studied numerical stability of a family of finite difference schemes for the dynamical plate equations. To obtain numerical stability, special treatments to our methods will be needed when applying to simulate dynamic evolution of shell structures with mixed boundary conditions.

For the study on the influences of thermal stresses and the precast shell shape, our investigation has been focused on the forward problem proposed in \cite{abbottmethods}  which involves numerically solving the governing equations to obtain the overall deflection. Related inverse problems would also be interesting and challenging. For instance, we may ask: is it possible to recover the precast shape with given deflection and thermal stresses, or is it possible to obtain the thermal stresses with given final deflection and precast shape? 
Whether these inverse problems are well-posed or not is still unclear, and some regularization may be needed in order to form an optimization problem to solve for these inverse problems.

\section*{Acknowledgement}
H. Ji's research is supported by the 2015 Mathematical Problems in Industry Workshop (MPI) Graduate Fellowship. L. Li's research is supported by the  National Science Foundation under Grant No. DMS-1261596, and the Margaret A. Darrin Postdoctoral Fellowship at Rensselaer Polytechnic Institute.
The authors thank John Abbott and Leslie Button of Corning Corporation and Thomas Witelski of Duke University. L. Li would like to acknowledge Professors D.W. Schewendeman and W.D. Henshaw of RPI for helpful conversations.

\bibliographystyle{elsart-num}
\bibliography{MPI}

\end{document}